\documentclass[10pt,a4paper]{article} 

\usepackage[latin1,utf8]{inputenc}
\usepackage[T1]{fontenc} 
\usepackage{textcomp} 
\usepackage[english]{babel}
\usepackage[autolanguage]{numprint} 
\usepackage{hyperref} 
\usepackage{graphicx} 
\usepackage{verbatim} 
\usepackage[all]{xy}
\usepackage{float}

\usepackage{lmodern}
\usepackage{mathrsfs}

\usepackage{amsmath,amssymb,amsthm} 
\usepackage{tikz,tkz-tab} 

\usepackage{mathdots}
\usepackage{makeidx}
\usepackage{stmaryrd} 

\usepackage{fancyhdr} 
\usepackage{enumitem}
\setenumerate[0]{label=(\alph*)}

\usepackage{a4wide}

\renewcommand\epsilon\varespilon 
\renewcommand\phi\varphi 

\newcommand\Adj{\mathrm{Adj}}
\newcommand\bP{\mathbf{P}} 
\newcommand\bL{\mathbf{L}} 
\newcommand\ba{\mathbf{a}}  
\newcommand\bb{\mathbf{b}}  
\newcommand\bx{\mathbf{x}}  
\newcommand\by{\mathbf{y}}  
\newcommand\bz{\mathbf{z}}  
\newcommand\bu{\mathbf{u}}  
\newcommand\bv{\mathbf{v}}  
\newcommand\bw{\mathbf{w}}  
\newcommand\CC{\mathbb{C}} 
\newcommand\CP{\mathrm{P}} 
\newcommand\CL{\mathrm{L}} 
\newcommand{\dist}[2]{\textrm{ {\rm dist}}\left(#1,#2\right)}
\newcommand\ee{\varepsilon}
\newcommand\espc{\,}
\newcommand\FF{\mathcal{F_{\textsl{stur}}}} 
\newcommand\GL{\mathrm{GL}_2}
\newcommand\GrO{\mathcal{O}} 
\newcommand\hCL{\widehat{\mathrm{L}}} 
\newcommand\hW{\widehat{W}} 
\newcommand\hY{\widehat{Y}} 
\newcommand\hZ{\widehat{Z}} 
\newcommand\hlambda{\widehat{\lambda}} 
\newcommand\homega{\widehat{\omega}} 

\newcommand\Id{\mathrm{Id}_{2}}
\newcommand\ie{\textsl{i.e. }} 
\newcommand{\indiceGauche}[2]{{\vphantom{#2}}^{#1}\!\!\!#2} 
\newcommand\Li[2]{\mathrm{L}_{#1}^{(#2)}} 

\newcommand\Mat{\mathrm{Mat}_{2\times2}}
\newcommand\MM{\mathcal{M}} 
\newcommand\NN{\mathbb{N}} 
\newcommand{\norm}[1]{\|#1\|} 

\newcommand\PeO{o} 
\newcommand\PP{\mathcal{P}} 
\newcommand\Proj{\mathbb{P}} 
\newcommand{\psc}[2]{\left\langle#1,#2\right\rangle} 
\newcommand\QQ{\mathbb{Q}} 
\newcommand\RR{\mathbb{R}} 
\newcommand\Sp{\textrm{spec}} 
\newcommand\Sturm{\mathcal{S}} 
\newcommand\SSS{\mathcal{A}} 
\newcommand\Tr{\mathrm{Tr}}
\newcommand{\transpose}[1]{\indiceGauche{t\;}{#1}}

\newcommand\ZZ{\mathbb{Z}} 

\newcommand\pu{\underline{\psi}} 
\newcommand\po{\overline{\psi}} 

\newcommand{\intervalle}[4]{\mathopen{#1}#2
\mathclose{}\mathpunct{};#3
\mathclose{#4}}

\newcommand{\intervalleEfo}[2]{\intervalle{\llbracket}{#1}{#2}{\llbracket}}

\theoremstyle{definition} 
\newtheorem{Def}{Definition}[section]
\newtheorem*{Nota}{Notation}

\theoremstyle{plain} 
\newtheorem{Prop}[Def]{Proposition} 
\newtheorem{Lem}[Def]{Lemma} 
\newtheorem{Thm}{Theorem}[section] 
\newtheorem{Cor}[Def]{Corollary} 

\theoremstyle{remark} 
\newtheorem{Exe}[Def]{Example} 
\newtheorem{Rem}[Def]{Remark} 
\newtheorem*{Rem2}{Remark} 

\newenvironment*{Dem}{\noindent{\bf Proof}}{\hfill$\square$} 

\numberwithin{equation}{section} 

\newcommand{\Addresses}{{
  \bigskip
  \footnotesize
  A.~Poëls, \textsc{Laboratoire de Math\'ematiques d'Orsay, Univ. Paris-Sud, CNRS,
Universit\'e Paris-Saclay, 91405 Orsay, France}\par\nopagebreak
  \textit{E-mail}: \texttt{anthony.poels@math.u-psud.fr}
}}

\newcommand{\MSC}{{
  \footnotesize
  \textbf{MSC~2010}: 11J13(Primary), 11H06 (Secondary), 11J82.
}}

\newcommand{\keysW}{{
  \footnotesize
  \textbf{Keywords}: Diophantine approximation, geometry of numbers, Sturmian sequence, simultaneous approximation
}}

\title{Exponents of diophantine approximation in dimension $2$ for numbers of Sturmian type}  
\author{Anthony Poëls} 

\begin{document} 
\maketitle

\begin{abstract}
We generalize the construction of Roy's Fibonacci type numbers to the case of a Sturmian recurrence and we determine the classical exponents of approximation $\omega_2(\xi)$, $\homega_2(\xi)$, $\lambda_2(\xi)$, $\hlambda_2(\xi)$ associated with these real numbers. This also extends similar results established by Bugeaud and Laurent in the case of Sturmian continued fractions. More generally we provide an almost complete description of the combined graph of parametric successive minima functions defined by Schmidt and Summerer in dimension two for such \textsl{Sturmian type} numbers. As a side result we obtain new information on the joint spectra of the above exponents as well as a new family of numbers for which it is possible to construct the sequence of the best rational approximations.
\end{abstract}

\MSC

\keysW

\section{Introduction}

If $\xi\in\RR$ is not an algebraic number of degree $\leq 2$ we may study the following two problems which are standard in diophantine approximation:\\

\textbf{Problem $E_{\omega,X}$:} We search for non-zero integer points $\bx = (x_0,x_1,x_2)\in\ZZ^3$ solutions of the system
\begin{align*}
\left\{ \begin{array}{ll}
|x_0+x_1\xi+x_2\xi^2|\leq X^{-\omega}\\
\max\big(|x_1|,|x_2|\big) \leq X.
\end{array} \right.
\end{align*}

\textbf{Problem $E_{\lambda,X}'$:} We search for non-zero integer points $\bx = (x_0,x_1,x_2)\in\ZZ^3$ solutions of the system
\begin{align*}
\left\{ \begin{array}{ll}
|x_0|\leq X\\
\max\big(|x_0\xi-x_1|,|x_0\xi^2-x_2|\big) \leq X^{-\lambda}.
\end{array} \right.
\end{align*}

We denote by $\omega_2(\xi)$ (resp. $\homega_2(\xi)$) the supremum of real numbers $\omega$ for which the problem $E_{\omega,X}$ admits a non-zero integer solution for arbitrarily large values of $X$ (resp. for each sufficiently large value of $X$).\\
Similarly we denote by $\lambda_2(\xi)$ (resp. $\hlambda_2(\xi)$) the supremum of real numbers $\lambda$ for which the problem $E_{\lambda,X}'$ admits a non-zero integer solution for arbitrarily large values of $X$ (resp. for each sufficiently large value of $X$).

We call \textsl{spectrum} of an exponent the set of values taken by this exponent at real numbers which are not algebraic of degree $\leq 2$. We also call \textsl{joint spectrum} of a family of diophantine exponents $(\nu_1,\nu_2,\dots)$ the set of values taken by $(\nu_1(\xi),\nu_2(\xi),\dots)$ at real numbers which are not algebraic of degree $\leq 2$. We denote this set by $\Sp(\nu_1,\nu_2,\dots)$.\\

The previous four exponents have been studied a lot during the last decade. The reader may refer to \cite{bugeaud2015exponents} which contains a well supplied summary of the subject and to \cite{bugeaud2005exponentsSturmian} for a generalization to the approximation of a vector $(1,\xi_1,\xi_2)\in\RR^3$ (in this paper we study the special case $(1,\xi_1,\xi_2) = (1,\xi,\xi^2)$). Let us recall briefly some important properties they satisfy. We shall first start with Jarn\'ik's identity \cite[Theorem 1]{jarnik1938khintchineschen} which binds $\hlambda_2$ and $\homega_2$ together:
\begin{Prop}[Jarn\'ik's identity]
For each real number $\xi$ which is not algebraic of degree $\leq 2$, we have
\begin{equation}
\label{Eq Jarnik}
    \hlambda_2(\xi) = 1 - \frac{1}{\homega_2(\xi)},
\end{equation}
with the convention that the right hand side of this equality is $1$ if $\homega_2(\xi)=+\infty$.
\end{Prop}
In this paper $\gamma$ denotes the golden ratio $\frac{1+\sqrt{5}}{2}$. It is known that these exponents satisfy
\begin{equation}
\label{eq intro encadrement hlambda et homega}
    \frac{1}{2} \leq \hlambda_2(\xi)\leq \frac{1}{\gamma}=0.618... \quad \textrm{and}\quad 2 \leq \homega_2(\xi)\leq \gamma^2=2.618...
\end{equation}
for any $\xi$ which is not algebraic of degree $\leq 2$. The lower bounds are obtained by the Dirichlet box principle (or equivalently, Minkowski's Theorem), the upper bounds follow respectively from Theorem $1a$ of
\cite{davenport1969approximation} and from \cite{arbourRoyCriterionDegreTwo}. Roy proved \cite{roy2003approximation} the existence of real numbers $\xi$ called \textsl{extremal numbers} for which $\hlambda_2(\xi)=\frac{1}{\gamma}$ and $\homega_2(\xi)= \gamma^2$. Then, generalizing his construction he proved \cite{roy2007two} that the spectrum of $\hlambda_2$ (resp. $\homega_2$) forms a dense subset of the interval $[\frac{1}{2},\frac{1}{\gamma}]$ (resp. $[2,\gamma^2]$). The numbers given by his generalized construction are said to be of the \textsl{Fibonacci type}. It is still an open problem to describe precisely the spectra of the two exponents $\hlambda_2,\homega_2$.\\
On the contrary, the spectra of $\lambda_2$ and $\omega_2$ are well known:
\begin{equation}
\label{eq intro spec lambda_2 et omega_2}
    \Sp(\lambda_2) = [\frac{1}{2},+\infty]\quad\textrm{and}\quad \Sp(\omega_2) = [2,+\infty].
\end{equation}
The first equality of \eqref{eq intro spec lambda_2 et omega_2} results directly from the works of Beresnevich, Dickinson, Vaughan and Velani \cite{beresnevich2007diophantine}, \cite{vaughan2006diophantine}, who computed the Hausdorff dimension of the set of real numbers $\xi$ for which $\lambda_2(\xi)=\lambda$. The second equality of \eqref{eq intro spec lambda_2 et omega_2} was proved by Bernik in \cite{bernik1983use}, computing for all $\omega\geq 2$ the Hausdorff dimension of the set of real numbers $\xi$ for which $\omega_2(\xi) = \omega$. Finally, note that for almost all $\xi$ with respect to the Lebesgue's measure, one has
\[
    \omega_2(\xi)=\homega_2(\xi)=\frac{1}{\lambda_2(\xi)}=\frac{1}{\hlambda_2(\xi)} = 2.
\]
(Cf \cite[Ch. VI, Corollaries 1C, 1E]{schmidt1996diophantine} or \cite[Theorem 2.3]{bugeaud2005exponentsSturmian}). This equality is also achieved for each $\xi$ algebraic number of degree $\geq 3$ by virtue of Schmidt's subspace theorem \cite{schmidt1996diophantine}. Cf also \cite[Theorem 2.10]{bugeaud2015exponents}.\\

To construct extremal numbers, Roy \cite{roy2003approximation} associates to distinct positive integers $a,b$ the real number $\xi_{a,b}$ whose partial quotients are successively $0$, and the letters of the Fibonacci word $w=abaaba\dots$ constructed on $\{a,b\}$ (\ie $\xi_{a,b} = [0;a,b,a,a,b,\dots]$). The property $\hlambda_2(\xi_{a,b})=1/\gamma$ is related to the fact that the Fibonacci word has many palindromic prefixes, and each of them gives a good approximation of the vector $(1,\xi_{a,b},\xi_{a,b}^2)$. Bugeaud and Laurent generalized this approach and calculated (among others) the four exponents $\lambda_2(\xi),\hlambda_2(\xi),\omega_2(\xi)$ and $\homega_2(\xi)$ for a real number $\xi$ whose partial quotients are given by a Sturmian characteristic word (see \cite{bugeaud2005exponentsSturmian}; also note that Schleischitz deals with the problem of the cubic approximation of such numbers in \cite{schleischitz2017approxcubicextremal}, \cite{schleischitz2017cubic}). Again, the fact that Sturmian characteristic words have many palindromic prefixes that one knows how to describe precisely plays an essential role. Note that for all extremal numbers and all numbers $\xi$ constructed by Bugeaud and Laurent, we have $\lambda_2(\xi)=1$. On the other hand, for Fibonacci type numbers Roy generalized the construction of extremal numbers in a different way using Fibonacci sequences in $\GL(\RR)$. One crucial point is that he did not consider directly the continued fraction representations (and the $2\times2$ associated matrices) which allowed him to leave out the condition $\lambda_2(\xi)=1$ and show the density of the spectra of the exponents $\hlambda_2$ and $\homega_2$ in the intervals defined by \eqref{eq intro encadrement hlambda et homega}.\\

In this paper we construct numbers of \textsl{Sturmian type}, which generalize both Roy's Fibonacci type numbers and the numbers constructed by Bugeaud and Laurent. For such a number $\xi$ we compute the four exponents $\hlambda_2(\xi)$, $\homega_2(\xi)$, $\omega_2(\xi),\lambda_2(\xi)$ and (at least partially) the associated $3$-system coming from Schmidt and Summerer's parametric geometry of numbers (see the end of this introduction for more details, including the definition of Sturmian type numbers). As far as we know all real numbers $\xi$ satisfying $\homega_2(\xi)>2$ and for which one knows how to compute the associated diophantine exponents and the best rational approximations sequence of $(1,\xi,\xi^2)$ are of Sturmian type\footnote{Note that we also have a construction which gives for every $\omega\geq3$ a number $\xi$ for which $\omega_2(\xi)=\omega$, but for such $\xi$ we always have $\homega_2(\xi) = 2$ (see the proof of Theorem $5.1$ and Theorem $5.5$ of \cite{bugeaud2015exponents}).}.\\
As an application of the construction and study of Sturmian type numbers, we prove density results on the spectra. To state these results we need the following definition.
\begin{Def}
Let $\underline{s} = (s_k)_{k\geq 0}$ be a sequence of positive integers. We associate to $\underline{s}$ the quantity
\[
    \sigma(\underline{s}) = \frac{1}{\limsup_{k\rightarrow+\infty}[s_{k+1};s_k,\dots,s_1]},
\]
where $[a_0;a_1,a_2,\dots]$ denotes the continued fraction whose partial quotients are $(a_0,a_1,\dots)$. We denote by $\SSS$ the set of all bounded sequences of positive integers, and as in  \cite{bugeaud2005exponentsSturmian}, $\Sturm$ denotes the set of $\sigma(\underline{s})$ for $\underline{s}\in\SSS$.
\end{Def}

Cassaigne studied the elements of $\Sturm$ in \cite{cassaigne1999limit}. This set is connected to the spectrum of Fischler's exponent $\beta_0(\xi)$ \cite{fischler2007palindromic}. To any $\underline{s}\in\SSS$ is associated a characteristic Sturmian word of angle $[0;s_1,s_2,\dots]$ and Sturmian type numbers $\xi$. In the following result $\xi$ is of Sturmian type and associated to $\underline{s}$.

\begin{Thm}
\label{Thm propre exposants intro}
For each $\underline{s}\in\SSS$, there exists a set $\Delta_{\underline{s}}$ which contains $0$ and is dense in the interval $[0,\frac{\sigma}{1+\sigma}]$ (with $\sigma = \sigma(\underline{s})$) and such that for each $\delta\in\Delta_{\underline{s}}$, there is a real number $\xi$ satisfying
\begin{align*}
&\omega_2(\xi) =  \frac{2-\delta}{\sigma}+1-\delta,\\
&\homega_2(\xi) =  1+(1-\delta)(1+\sigma),\\
&\hlambda_2(\xi) = \frac{(1-\delta)(1+\sigma)}{1+(1-\delta)(1+\sigma)},\\ &1-\delta\leq\lambda_2(\xi)\leq\max\Big(1-\delta,\frac{1}{1-\delta+\sigma}\Big).
\end{align*}
Moreover if we have $\delta\in[0,h(\sigma)]$ with
\[
    h(\sigma) = \frac{\sigma}{2}+1-\sqrt{\big(\frac{\sigma}{2}\big)^2+1}\leq \frac{\sigma}{1+\sigma},
\]
then
\[
    \lambda_2(\xi) = 1-\delta.
\]

\end{Thm}

Note that the case $\delta = 0$ of Theorem \ref{Thm propre exposants intro} corresponds to the values of the four exponents for the number $\xi_{\phi}$ constructed by Bugeaud and Laurent (see Theorem~$3.1$ of \cite{bugeaud2005exponentsSturmian}); for a Fibonacci type number we have $s_k=1$ for all $k\geq 1$ and the associated $\sigma$ is $\sigma = 1/\gamma$.\\
We may deduce from Theorem \ref{Thm propre exposants intro} the following corollary about the spectrum of $\omega_2$:

\begin{Cor}
\label{Thm densité omega_2 sturm dans [c,+infty]}
We have
\[
    \overline{\{\omega_2(\xi)\;|\;\xi \textrm{ of Sturmian type}\}} = [1+\sqrt{5},2+\sqrt 5]\cup[2+2\sqrt 2,3+2\sqrt 3]\cup[3+\sqrt{13}+\infty\big),
\]
where $\overline{A}$ denotes the topological closure of the subset $A\subset\RR$.
\end{Cor}
One may compare this result with the special case where $s_k=1$ for all $k$: $\overline{\{\omega_2(\xi)\;|\;\xi \textrm{ of Fibonacci type}\}} = [1+\sqrt{5},2+\sqrt 5]$.\\
These new constructions do not seem to give additional information on the individual spectra of each exponent. However they bring new information on their joint spectrum. In this respect, a direct corollary of Theorem \ref{Thm propre exposants intro} is the following.

\begin{Cor}\
\label{Cor 2 spectre}
For each $\sigma\in\Sturm$, there is $c=c(\sigma)$ satisfying $0\leq c < 1$, such that the joint spectrum of $(\lambda_2,\hlambda_2,\omega_2,\homega_2)$ is dense in the curve
\[
    \Big\{\big(x,1-\frac{1}{1+(1+\sigma)x},\frac{1+(1+\sigma)x}{\sigma},1+(1+\sigma)x\big)\;\Big|\;x\in[c,1] \Big\}.
\]

\end{Cor}

Using their powerful theory of parametric geometry of numbers (\cite{Schmidt2009}, \cite{Schmidt2013}) Schmidt and Summerer proved that the study of the previous diophantine approximation exponents may be reduced to the study of what they call a $(3,0)$-system (and that Roy calls a $3$-system in \cite{Roy_juin}, cf Definition \ref{Def n-système} below). In this paper, we give an almost complete description of the $3$-system associated to a number $\xi$ of Sturmian type (see Figure \ref{figure 3system_partiel2.png} in Section \ref{Section 3-système (partiel) d'un nombre psi sturmien}), and from it we compute the four exponents of Theorem \ref{Thm propre exposants intro}.\\

To conclude this introduction, let us summarize our strategy, parallel to that of Roy \cite{roy2007two}. Let $(s_k)_{k\geq 0}$ be a sequence of positive integers (except for $s_0$) with $s_0 = -1, s_1 = 1$. For $k\geq 0$, we set $t_k = s_0+s_1+\dots+s_k$, and we associate to $(s_k)_k$ a function $\psi$ (called a \textsl{Sturmian} function, see Definition \ref{Def fonction sturmienne}) which is defined on $\intervalleEfo{0}{+\infty}$ by $\psi(t_k) = t_{k-1} - 1$ for $k\geq 1$ and $\psi(i) = i-1$  if $i$ is not among the $t_k$. We denote by $\FF$ the set of all Sturmian functions for which the corresponding sequence $(s_k)_k$ is bounded. Such a function $\psi$ corresponds to a characteristic Sturmian word through the recurrence relation of its palindromic prefixes (see Section $3.1$ of \cite{fischler2006palindromic}).\\
In order to define a proper $\psi$-Sturmian number, we need the two following definitions.
A $\psi$\textsl{-Sturmian} sequence in $\MM=\GL(\RR)\cap\Mat(\ZZ)$ is a sequence $(\bw_i)_{i\geq 0}$ of matrices such that
\[
    \bw_{i+1} = \bw_i^{s_{i+1}}\bw_{i-1}\quad (i\geq 1).
\]
Such a sequence is \textsl{admissible} if there exists a matrix $N\in\MM$ such that $\bw_1N, \bw_0\transpose{N}$ and $\bw_1\bw_0\transpose{N}$ are symmetric matrices (see Definition \ref{Def monoide suite admissible}). We set $N_k = N$ if $k$ is even, $N_k = \transpose{N}$ if $k$ is odd.\\
We say that a $\psi$-Sturmian sequence $(\bw_i)_{i\geq 0}$ has \textsl{multiplicative growth} if
\[
    \norm{\bw_k^l\bw_{k-1}} \asymp \norm{\bw_k}\times\norm{\bw_k^{l-1}\bw_{k-1}}
\]
for $k \geq 1$, $1\leq l \leq s_{k+1}+1$ (see Definition \ref{Def croissance multiplicative}). Sequences with multiplicative growth are studied in Section \ref{section suites psi-sturm croissance multi}. In the following we identify $\RR^3$ (resp. $\ZZ^3$) with the space of $2\times 2$ symmetric matrices with real (resp. integer) coefficients under the map $(x_0,x_1,x_2) \longrightarrow \left( \begin{array}{cc} x_0 & x_1 \\ x_1 & x_2\end{array}\right)$. Accordingly, it makes sense to define the determinant $\det(\bx)$ of a point $\bx\in\RR^3$. Similarly, given symmetric matrices $\bx,\by$ we write $\bx\wedge\by$ to denote the inner product of the corresponding points.\\
From now on, we consider an admissible $\psi$-Sturmian sequence in $\MM$ with multiplicative growth such that the content of $\bw_i$ is bounded and $(\norm{\bw_i})_i$ tends to infinity. For instance, the $\psi$-Sturmian sequence defined by Roy's matrices \cite{roy2007two}
\[
    \bw_0 = \left(
    \begin{array}{cc}
        1 & b \\
        a & a(b+1)
    \end{array}
    \right),
    \quad
    \bw_1 = \left(
    \begin{array}{cc}
        1 & c \\
        a & a(c+1)
    \end{array}
    \right)
\]
satisfies these hypotheses (see Section \ref{Section Exemple de Roy}).\\
Then we define two sequences of symmetric matrices $(\by_i)_i$ and $(\bz_i)_i$ by the formulas
\[
    \by_{t_k+l} = \bw_k^{l+1}\bw_{k-1}N_k\quad\textrm{and} \quad \bz_{t_k+l} = \frac{1}{\det(\bw_k)}\by_{\psi(t_{k+1})}\wedge \by_{t_k+l}\quad (0\leq k, 0\leq l < s_{k+1})
\]
(see Definitions \ref{Def (y_i)} and \ref{Def (z_i)_i}). Combinatorial properties of these sequences are studied in Sections \ref{section suites psi-sturmiennes de GL2(C)} and \ref{section suites psi-sturmiennes admissibles}. By multiplicative growth (and because $(\norm{\bw_i})_i$ tends to infinity) we show in Proposition \ref{Prop existence delta sur le det} that there exists $\delta\geq 0$ such that $|\det(\bw_i)|\asymp \norm{\bw_i}^{\delta}$. If $\delta < \frac{\sigma}{1+\sigma}$, with $\sigma = 1/\limsup_{k\rightarrow\infty}[s_{k+1};s_k,\dots,s_1]$, then Proposition \ref{Prop existence y = (1,xi,xi^2)} ensures that the associated sequence $(\by_i)_i$ converges in $\Proj^2(\RR)$ to a point $\by = (1,\xi,\xi^2)$. The number $\xi$ therefore produced is called a proper $\psi$-Sturmian number (and a Sturmian type number is a proper $\psi$-Sturmian number for some $\psi\in\FF$; see Definition \ref{Def nombre de type sturmien}). Then we show that the sequences $(\by_i)_i$ and $(\bz_i)_i$ are ``good'' solutions of Problems $E'_{\lambda,X}$ and $E_{\omega,X}$ respectively. In Section \ref{section xi construction de Roy} (see Proposition \ref{Prop estim normes yi zj}) we prove the existence of $\xi$ and give precise estimates of $\norm{\by_i}$, $\norm{\bz_i}$, $|\psc{\bz_i}{\by}|$ and $\norm{\by_i\wedge\by}$. This allows us to build in Section \ref{Section 3-système (partiel) d'un nombre psi sturmien} (see Propositions \ref{Prop Dessin L_j, L_j*}, \ref{Prop 3-sytème partiel} and Figure \ref{figure 3system_partiel2.png}) a partial $3$-system associated to $\xi$. Finally, we deduce from this $3$-system the Diophantine exponents associated to $\xi$ (see Theorem \ref{Thm exposants geom param nb psi-sturmien}). For $\lambda_2(\xi)$ an additional hypothesis --which is discussed in Section \ref{Section zones d'incertitude} -- is required in order to ignore the gray areas of the partial $3$-system of Figure \ref{figure 3system_partiel2.png}.

\section{Notations}
\label{section Notations}

In this section we set some notations (in particular the sequence $(s_k)_{k\geq 0}$ and the function $\psi$ that we will constantly use throughout this paper) and we introduce the notion of Sturmian functions.

\begin{Nota}
Let $I$ be a set (the \textsl{set of index}, typically $I$ of the form $\NN^r$), $(a_{\underline{i}})_{\underline{i}\in I}$ and $(b_{\underline{i}})_{\underline{i}\in I}$ two sequences of non negative real numbers indexed by $I$. Let $J$  be a subset of $I$ (the \textsl{set of conditions}). We denote by ``$a_{\underline{i}} \ll b_{\underline{i}}$ for $\underline{i} \in J$'' or ``$ b_{\underline{i}} \gg a_{\underline{i}}$ for $\underline{i} \in J$'' if there is a constant $c>0$ such that for each $\underline{i}\in J$ we have $a_{\underline{i}} \leq cb_{\underline{i}}$. We write ``$a_{\underline{i}}\asymp b_{\underline{i}}$ for $\underline{i} \in J$'' if $a_{\underline{i}}\ll b_{\underline{i}}$ for $\underline{i} \in J$ and $b_{\underline{i}}\ll a_{\underline{i}}$ for $\underline{i} \in J$.\\
In the special case $I = \NN$, without more precisions we will always implicitly take $J$ of the form $\intervalleEfo{j_0}{+\infty}$ for $j_0$ large enough, and we will simply write $a_{\underline{i}} \ll b_{\underline{i}}$, $b_{\underline{i}} \gg a_{\underline{i}}$ and $a_{\underline{i}} \asymp b_{\underline{i}}$.
\end{Nota}

\begin{Def}
For each non-zero $\bx\in\RR^3$, we denote by $[\bx]$ its image in $\Proj^2(\RR)$. Commonly, we define the projective distance between two non-zero vectors $\bx$ and $\by$ of $\RR^3$ by
\[
    \dist{[\bx]}{[\by]} =\dist{\bx}{\by} := \frac{\norm{\bx\wedge \by}}{\norm{\bx}\espc\norm{\by}}.
\]
The projective distance satisfies in particular the triangle inequality:
\[
    \dist{\bx}{\bz} \leq \dist{\bx}{\by} + \dist{\by}{\bz}\quad \textrm{for all } \bx,\by,\bz\in \RR^3\setminus\{0\}.
\]
\end{Def}

We also recall the following inequalities, valid for $\bx, \by, \bz\in\RR^3$ (see \cite[Lemma 2.2]{roy2007two}):
\begin{Lem}
We have
\begin{align}
\label{formule technique wedge 1}
& \norm{\psc{\bx}{\bz}\by-\psc{\bx}{\by}\bz} \leq 2\norm{\bx}\espc\norm{\by\wedge\bz},\\
\label{formule technique wedge 2}
& \norm{\by}\espc\norm{\bx\wedge\bz} \leq \norm{\bz}\espc\norm{\bx\wedge\by} + 2 \norm{\bx}\espc\norm{\by\wedge\bz}.
\end{align}
\end{Lem}

\begin{Def}
\label{Def fonction sturmienne}
Let $(s_k)_{k\geq 0}$ be a sequence of positive integers (except for $s_0$) with $s_0 = -1, s_1 = 1$. For $k\geq 0$, we set $t_k = s_0+s_1+\dots+s_k$, and we associate to $(s_k)_k$ a function $\psi$ defined on $\intervalleEfo{0}{+\infty}$ by
\begin{equation*}
\left\{ \begin{array}{ll}
\psi(t_k) = t_{k-1} - 1 \textrm{ for } k\geq 1 \\
\psi(i) = i-1 \textrm{ if $i$ is not among the $t_k$.}
\end{array} \right.
\end{equation*}
The function $\psi$ is called \textsl{Sturmian}.\\
We denote by $\FF$ the set of all Sturmian functions for which the corresponding sequence $(s_k)_k$ is bounded.
\end{Def}

\begin{Rem2}
 Such a sequence is entirely characterised by its associated Sturmian function (the sequence $(t_k)_{k\geq 1}$ being exactly the integers $n$ such that $\psi(n)\leq n-2$).\\
 The Sturmian functions $\psi\in\FF$ belong to a bigger class of functions called \textsl{asymptotically reduced} functions and studied by Fischler in \cite{fischler2006palindromic} and \cite{fischler2007palindromic} (cf \cite[Definition 2.1]{fischler2007palindromic}). Generalizing the construction of Roy and the use of his bracket (cf \cite{roy2003approximation} and\cite{roy2004approximation}) Fischler shows that for each word $w$ with \textsl{abundant palindromic prefixes}, we may associate a function $\psi$ asymptotically reduced (if $w$ is a \textsl{Sturmian} word, the associated function $\psi$ is Sturmian) and that we may build a real number $\xi$ whose certain diophantine exponents are intimately connected to $\psi$.
\end{Rem2}

\begin{Nota}
For the rest of the paper, we set a sequence $(s_k)_{k\geq 0}$ of positive integers (except for $k=0$, with $s_0 = -1, s_1 = 1$), which is not necessarily bounded. We set $t_k = s_0+\dots + s_k$ and we denote by $\psi$ the associated Sturmian function. We also set $\phi = [0;s_2,s_3,\dots]$.
\end{Nota}

\section{$\psi$-Sturmian sequences in $\GL(\CC)$}
\label{section suites psi-sturmiennes de GL2(C)}

The goal of this section is to generalize Roy's definition of Fibonacci sequence in $\GL(\CC)$ (see \cite[§3]{roy2007two}) and to identify in Proposition \ref{Propriete sur y_k} some useful properties of the sequence $(\by_i)_i$ (see Definition \ref{Def (y_i)}). This sequence will give ``good'' solutions of Problem $E'_{\lambda,X}$. Several properties have combinatorial analogues for characteristic Sturmian words (see \cite[§5]{bugeaud2005exponentsSturmian}).

\begin{Def}
A $\psi$\textsl{-Sturmian} sequence in a monoid is a sequence $(\bw_i)_{i\geq 0}$ such that
\[
    \bw_{i+1} = \bw_i^{s_{i+1}}\bw_{i-1} \quad (i\geq 1).
\]
Clearly, such a sequence is entirely determined by its first two elements $\bw_0$ and $\bw_1$. We say that a $\psi$-Sturmian sequence is a \textsl{Sturmian sequence of angle} (or slope) $\phi$ (where $\phi = [0;s_2,s_3,\dots]$).
\end{Def}

\begin{Exe}
 In the special case $\phi = [0;1,1,1,\dots] = 1/\gamma$ (then we have \mbox{$(t_k)_{k\geq 0} = (k-1)_{k\geq0}$} and $\psi(n) = n-2$ for all $n$), a Sturmian sequence of angle $1/\gamma$ in a monoid is a Fibonacci sequence according to Roy (cf at the beginning of §~$3$ of \cite{roy2007two}).
\end{Exe}

The following Proposition is a direct consequence of \cite[Proposition $3.1$]{roy2007two}.

\begin{Prop}[Roy, 2007]
\label{Prop Roy ouvert Zariski}
There exists a non-empty Zariski open subset $\mathcal{U}$ of $\GL(\CC)^2$ such that for all $(\bw_0,\bw_1)\in \mathcal{U}$, there is $N\in \GL(\CC)$ such that $\bw_1N, \bw_0\transpose{N}$ and $\bw_1\bw_0\transpose{N}$ are symmetric. When  $\bw_0$ and $\bw_1$ have integer coefficients, we may take $N$ with integer coefficients.
\end{Prop}

\begin{Prop} \
\label{Prop w_kw_k-1 = w_k-1wk}
\begin{enumerate}
\item \label{Prop 1. w_kw_k-1 = w_k-1wk}
Let $\bw_0, \bw_1, N \in \GL(\CC)$ be such that $\bw_1N, \bw_0\transpose{N}$ and $\bw_1\bw_0\transpose{N}$ are symmetric. Then $\bw_1\bw_0\transpose{N} = \bw_0\bw_1N$.
\item
\label{Prop 2. w_kw_k-1 = w_k-1wk}
Let $(\bw_i)_{i\geq 0}$ be a $\psi$-Sturmian sequence of $\GL(\CC)$ and let $N$ be $ \GL(\CC)$. We set $N_k = N$ if $k$ is even, $N_k = \transpose{N}$ if $k$ is odd. Suppose that \mbox{$\bw_1\bw_0\transpose{N} = \bw_0\bw_1N$}.  Then for all $k\geq 1$ we have
 \begin{equation}
 \label{Eq w_kw_k-1 = w_k-1wk}
    \bw_{k-1}\bw_k N_{k+1} = \bw_k\bw_{k-1}N_k.
 \end{equation}
\end{enumerate}
\end{Prop}

\begin{Dem}
For \ref{Prop 1. w_kw_k-1 = w_k-1wk}: we have the equality $\bw_1\bw_0\transpose{N} = \big(\bw_1N\big)N^{-1}\big(\bw_0\transpose{N}\big)$. By taking the transpose of these expressions and using the fact that $\bw_1N, \bw_0\transpose{N}$ and $\bw_1\bw_0\transpose{N}$ are symmetric, we find
\[
    \bw_1\bw_0\transpose{N} = \big(\bw_0\transpose{N}\big)\transpose{N^{-1}}\big(\bw_1N\big) = \bw_0\bw_1N.
\]
For \ref{Prop 2. w_kw_k-1 = w_k-1wk} (one may compare this Proposition with Lemma $5.1$ of \cite{bugeaud2005exponentsSturmian}): one may adapt the arguments of the Proposition $1$ of \cite{allouche2001transcendence} to our context.\\
The case $k=1$ is true by hypothesis. Suppose that for $k\geq 1$ one has $\bw_{k-1}\bw_k N_{k+1} = \bw_k\bw_{k-1}N_k$. Then, one has
\begin{align*}
    \bw_k\bw_{k+1}N_{k+2} = \bw_k\big(\bw_k^{s_{k+1}}\bw_{k-1}N_{k+2}\big) & = \bw_k^{s_{k+1}}\big(\bw_k\bw_{k-1}N_k\big) \\
     & = \bw_k^{s_{k+1}}\big(\bw_{k-1}\bw_k N_{k+1}\big)\\
     & = \bw_{k+1}\bw_kN_{k+1},
\end{align*}
and we can conclude by induction.

\end{Dem}

\begin{Def}
\label{Def (y_i)}
Let $(\bw_i)_{i\geq 0}$ be a $\psi$-Sturmian sequence in $\GL(\CC)$ and let $N\in\GL(\CC)$ be such that $\bw_1\bw_0\transpose{N} = \bw_0\bw_1N$. We set $N_k = N$ if $k$ is even, $N_k = \transpose{N}$ if $k$ is odd. We define a sequence $(\by_j)_{j\geq -2}$ in the following way: we set $\by_{-2} = \bw_0\transpose{N}$, $\by_{-1} = \bw_1N$, and for $1\leq k$, $0\leq l < s_{k+1}$, we set
\begin{align}
\label{Eq Def y_i}
    \by_{t_k + l} = \bw_k^{l+1}\bw_{k-1}N_k.
\end{align}
(note that \eqref{Eq Def y_i} remains valid for $l=s_{k+1}$ by Eq. \eqref{Eq w_kw_k-1 = w_k-1wk}).\\
In particular, we have
\begin{equation}
\label{Eq y_psi(k) = w_k-1}
    \by_{\psi(t_k)} = \bw_{k-1}N_k \quad (k\geq 1).
\end{equation}
\end{Def}

\begin{Prop}
\label{Propriete sur y_k}
We have the following properties:
\begin{enumerate}
\item $\by_{t_k+l} = \bw_k^{l+1} \by_{\psi(t_k)}$ for $k\geq 1$ and $0\leq l \leq s_{k+1}$
\label{Propriete 1 sur y_k}
\item $\by_{j+1} = \bw_k\by_j$, for $k\geq 0$ and $t_k\leq j < t_{k+1}$
\label{Propriete 2 sur y_k}
\item $\by_{j} = \bw_k\by_{\psi(j)}$, for $k\geq 1$, $t_k\leq j < t_{k+1}$
\label{Propriete 4 sur y_k}
\item For all $j \geq 0$, we have
\label{Propriete 3 sur y_k}
\begin{equation}
    \label{Eq y_j+1 = crochet}
    \by_{j+1} = \by_j\by_{\psi(j)}^{-1}\by_j
\end{equation}
In particular, if $\bw_1N, \bw_0\transpose{N}$ and $\bw_1\bw_0\transpose{N}$ are symmetric, then for each \mbox{$j\geq -2$}, $\by_j$ is symmetric.
\end{enumerate}
\end{Prop}

\begin{Dem} The property \ref{Propriete 1 sur y_k} (resp. \ref{Propriete 2 sur y_k}) results from the definition of $\by_m$ and from \eqref{Eq y_psi(k) = w_k-1}.\\
The property \ref{Propriete 4 sur y_k} is implied by \ref{Propriete 2 sur y_k} if $j>t_k$ (since the we have $\psi(j) = j-1$) and by \ref{Propriete 1 sur y_k} si $j=t_k$.\\
Finally \eqref{Eq y_j+1 = crochet} is implied by \ref{Propriete 2 sur y_k} and \ref{Propriete 4 sur y_k}.
Furthermore, the property that $\bw_0\transpose{N}$, $\bw_1N$ and $\bw_1\bw_0\transpose{N}$ are symmetric is equivalent by definition to $\by_{-2}$, $\by_{-1}$ and $\by_0$ are symmetric: we can conclude by induction using \eqref{Eq y_j+1 = crochet}.

\end{Dem}

\begin{Rem2}
Equation \eqref{Eq y_j+1 = crochet} may be rewritten as
\begin{equation}
    \label{Eq y_j+1 = crochet 2}
    \det(\by_{\psi(j)})\by_{j+1} = [\by_j,\by_j,\by_{\psi(j)}],
\end{equation}
where $[\cdot,\cdot,\cdot]$ is Roy's bracket introduced in \cite{roy2004approximation}. This bracket is used in a crucial way in \cite{fischler2007palindromic} by Fischler (see \cite[§2]{fischler2007palindromic} and the first remark below Theorem $4.1$ that may be linked to \eqref{Eq y_j+1 = crochet 2}). Equation \eqref{Eq y_j+1 = crochet} reflects common combinatorial properties of Sturmian words. One may for instance compare it with Eq. $(4)$ of \cite{fischler2007palindromic}.
\end{Rem2}

\section{Admissible $\psi$-Sturmian sequences}
\label{section suites psi-sturmiennes admissibles}

In this section we establish several arithmetical properties of admissible $\psi$-Sturmian sequences in the monoid $\MM$ (see Definition \ref{Def monoide suite admissible}), including the generalized recurrence found by Bugeaud and Laurent (cf Lemma $7.1$ of \cite{bugeaud2005exponentsSturmian} and its proof). We also generalize in Definition \ref{Def (z_i)_i} the construction of the sequence $(\bz_i)_i$ of \cite[Proposition $4.1$]{roy2007two}. This sequence will give ``good'' solutions of Problem $E_{\omega,X}$ and plays a crucial role to build a $3$-system associated to a proper $\psi$-Sturmian number and it allows us to calculate standard diophantine exponents $\hlambda_2,\lambda_2,\homega_2$ and $\omega_2$. The main result of this section is Proposition \ref{Prop arithm Tr, det, y_k, z_k}.\\
The following definition is directly inspired by Definition $3.2$ of \cite{roy2007two}.

\begin{Def}
\label{Def monoide suite admissible}
Let $\MM = \Mat(\ZZ)\cap\GL(\CC)$ denote the monoid of $2\times2$ integer matrices with non-zero determinant.
We say that a $\psi$-Sturmian sequence in $\MM$ is \textsl{admissible} if there is a matrix $N\in\MM$ such that $\bw_1N, \bw_0\transpose{N}$ and $\bw_1\bw_0\transpose{N}$ are symmetric. In this case, matrices of the sequence $(\by_i)_{i\geq-2}$ associated by the Definition \ref{Def (y_i)} are symmetric.
\end{Def}

The examples of admissible sequences used by Roy in \cite{roy2007two} lead to the construction of admissible $\psi$-Sturmian sequences. Details are gathered in Section \ref{Section Exemple de Roy}. Following Roy's approach (cf §4 of \cite{roy2007two}) we identify $\RR^3$ (resp. $\ZZ^3$) with the space of $2\times 2$ symmetric matrices with real (resp. integer) coefficients under the map
\[
   \bx = (x_0,x_1,x_2) \longrightarrow \left(
    \begin{array}{cc}
        x_0 & x_1 \\
        x_1 & x_2
    \end{array}
    \right)  .
\]
If we write $\bx = (x_0,x_1,x_2)\in\RR^3$, then set $\det(\bx) = x_0x_2-x_1^2$ and $\Tr(\bx) = x_0+x_2$. Likewise, given three symmetric matrices $\bx, \by, \bz$, let $\bx\wedge\by$, $\psc{\bx}{\by}$ and $\det(\bx,\by,\bz)$ denote respectively the inner product, the scalar product and the determinant of corresponding vectors of $\RR^3$. We define the \textsl{content} of a matrix $\bw\in\Mat(\ZZ)$ or of a point $\by\in\ZZ^3$ as the greatest common divisor of their coefficients. We say that such a matrix or point is \textsl{primitive} if its content is $1$.

\begin{Nota}
In this section, $(\bw_i)_{i\geq 0}$ is an admissible $\psi$-Sturmian sequence. We denote by $N\in\MM$ and $(\by_j)_{j\geq -2}$ the matrices defined in Definition \ref{Def (y_i)}.\\
Moreover let $J$ be the matrix $J = \left(\begin{array}{cc}
        0 & 1 \\
        -1 & 0
    \end{array}
    \right).$
\end{Nota}

\begin{Def}
\label{Def (z_i)_i}
We define the sequence $(\bz_j)_{j\geq -1}$ in the following way:
\begin{equation}
\label{Eq Def z_i}
    \bz_{t_k+l} = \frac{1}{\det(\bw_k)}\by_{\psi(t_{k+1})}\wedge \by_{t_k+l},\quad 0\leq k\;\textrm{ and } 0\leq l< s_{k+1}.
\end{equation}
\end{Def}

Proposition \ref{Prop arithm Tr, det, y_k, z_k} below and Corollary \ref{Cor arithm Tr, det, y_i, z_j} extend results of Proposition $4.1$ and of Corollary $4.2$ in \cite{roy2007two}.

\begin{Prop} \
\label{Prop arithm Tr, det, y_k, z_k}
\begin{enumerate}
\item  Recurrence of sequences $(\Tr(\bw_i))_{i}$ and $(\det(\bw_i))_{i}$:
\begin{align}
\label{Eq recurrence Tr et det}
\Tr(\bw_{k}^l\bw_{k-1}) = \Tr(\bw_{k})\Tr(\bw_{k}^{l-1}\bw_{k-1})-\det(\bw_k)\Tr(\bw_{k}^{l-2}\bw_{k-1})
\end{align}
where $k\geq 1$ and $l\geq 2$. In particular, we have
\begin{align}
\label{Eq Tr mod det 1}
\Tr(\bw_{k}^l\bw_{k-1}) &\equiv \Tr(\bw_{k})^{l-1}\Tr(\bw_{k}\bw_{k-1}) \mod \det(\bw_k)
\end{align}
\item  Recurrence of the sequence $(\by_i)_{i}$:
\begin{equation}
\label{Eq recurrence y_i}
    \by_{t_k+l+1} = \Tr(\bw_k)\by_{t_k+l} - \det(\bw_k)\by_{\psi(t_k+l)}, \quad 0\leq l < s_{k+1},
\end{equation}
with $k\geq 1$. In particular, the case $k-1$ applied with $l=s_k -1$ gives
\begin{equation}
\label{Eq recurrence y_i bis}
    \by_{t_k} = \Tr(\bw_{k-1})\by_{t_k-1} - \det(\bw_{k-1})\by_{\psi(t_k-1)}, \quad k\geq 2,
\end{equation}
\item  Recurrence of the sequence $(\bz_i)_{i}$:
\begin{align}
\label{Eq recurrence z_i}
\left\{ \begin{array}{ll}
\bz_{t_k+l+1} = \Tr(\bw_k)\bz_{t_{k}+l} - \by_{\psi(t_{k+1})}\wedge \by_{\psi(t_{k}+l)} \quad (0 \leq l < s_{k+1}-1)\\
\bz_{t_{k+1}} = \Tr(\bw_{k-1})\bz_{t_k-1} - \by_{\psi(t_{k})}\wedge \by_{\psi(t_{k}-1)} \quad (k \geq 2),
\end{array} \right.
\end{align}
the first relation remaining valid for $k\geq 1$.
\item For $k\geq 0$ we have
\begin{align}
\label{Eq det(y,y,y)}
\det(\by_{t_k-1}, \by_{t_k}, \by_{t_k+1}) = -\det(\bw_k)\det(\by_{t_k})\Tr(JN_{k+1}),
\end{align}
in particular, if $N$ is not symmetric (which is equivalent to $\Tr(JN) \neq 0$), then $\by_{t_k-1}$, $\by_{t_k}$ and $\by_{t_k+1}$ are linearly independent for all $k\geq 0$.
\item For $k\geq 0$ we have
\begin{align}
\label{Eq wedge des z_j}
\bz_{t_{k+1}}\wedge \bz_{t_k+l} = \det(N_k)\Tr(JN_{k+1})\by_{t_k+l}\quad (0\leq l < s_{k+1}).
\end{align}
\end{enumerate}
\end{Prop}

\begin{Dem}\
Let $k\geq 1$. Following Roy's idea and using Cayley-Hamilton Theorem for $\bw_k$, one finds
\begin{equation}
\label{Eq Cayley-Hamilton}
    \bw_k^2 = \Tr(\bw_k)\bw_k-\det(\bw_k)\Id,
\end{equation}
where $\Id$ is the identity matrix of $\Mat(\RR)$. Then, we prove \eqref{Eq recurrence Tr et det} by multiplying each side of the previous equality on the left by $\bw_k^{l-2}$ and on the right by $\bw_{k-1}$ (with $l\geq 2$), then by taking the trace. One gets Equation \eqref{Eq Tr mod det 1} by induction by taking \eqref{Eq recurrence Tr et det} modulo $\det(\bw_k)$.\\
For recurrence \eqref{Eq recurrence y_i}, one starts with \eqref{Eq Cayley-Hamilton}, and by multiplying each side of this equality on the left by $\bw_k^{l}$ (with $k\geq 1$, $0\leq l < s_{k+1}$) and on the right by $\bw_{k-1}N_k$, one finds
\[
    \bw_{k}^{l+2}\bw_{k-1}N_k = \Tr(\bw_{k})\bw_{k}^{l+1}\bw_{k-1}N_k-\det(\bw_k)\bw_{k}^{l}\bw_{k-1}N_k,
\]
which by the definition of the $\by_i$ may be rewritten as
\[
    \by_{t_k+l+1} = \Tr(\bw_{k})\by_{t_k+l}-\det(\bw_k)\bw_{k}^{l}\by_{\psi(t_k)}.
\]
To conclude it suffices to note that for each $0<l<s_{k+1}$, we have $\by_{\psi(t_k+l)} = \by_{t_k+l-1} = \bw_{k}^{l}\by_{\psi(t_k)}$. Also note that \eqref{Eq recurrence y_i} is equivalent to $(19)$ of \cite{bugeaud2005exponentsSturmian}.\\
We find the first relation of \eqref{Eq recurrence z_i} by using \eqref{Eq recurrence y_i} together with the definition of $\bz_{t_k+l+1}$ and $\bz_{t_k+l}$. For the second relation, using \eqref{Eq recurrence y_i}, we obtain by induction on $l$, $0\leq l < s_{k+1}$,
\begin{equation}
\label{Eq inter sur prod vect y_l}
    \by_{t_k+l}\wedge\by_{t_k+l+1} = \det(\bw_k)^{l}\by_{t_k}\wedge\by_{t_k+1} = \det(\bw_k)^{l+1}\by_{\psi(t_k)}\wedge\by_{t_k}.
\end{equation}
Similarly, we may show that for each point $\bx\in\RR^3$,
\begin{equation}
\label{Eq inter sur prod vect y_l bis}
    \det(\bx,\by_{t_k+l},\by_{t_k+l+1}) = \det(\bw_k)^{l}\det(\bx,\by_{t_k},\by_{t_k+1}).
\end{equation}
Equation \eqref{Eq inter sur prod vect y_l} implies that
\begin{align*}
\det(\bw_{k+1})\bz_{t_{k+1}} = \by_{t_{k}+s_{k+1}-1}\wedge \by_{t_{k}+s_{k+1}} = \det(\bw_k)^{s_{k+1}}\by_{\psi(t_{k})}\wedge \by_{t_{k}}.
\end{align*}
If $k\geq 2$, Equation \eqref{Eq recurrence y_i bis} provides the relation
\begin{align*}
    \by_{\psi(t_{k})}\wedge \by_{t_{k}} =  \Tr(\bw_{k-1})\det(\bw_{k-1})\bz_{t_{k}-1} - \det(\bw_{k-1})\by_{\psi(t_{k})}\wedge \by_{\psi(t_k-1)}.
\end{align*}
Finally by using $\det(\bw_{k+1}) = \det(\bw_k)^{s_{k+1}}\det(\bw_{k-1})$ we obtain
\begin{align*}
    \bz_{t_{k+1}} &=  \frac{\det(\bw_k)^{s_{k+1}}\det(\bw_{k-1})}{\det(\bw_{k+1})}\times\big(\Tr(\bw_{k-1})\bz_{t_{k}-1} - \by_{\psi(t_{k})}\wedge \by_{\psi(t_k-1)} \big) \\
    & = \Tr(\bw_{k-1})\bz_{t_{k}-1} - \by_{\psi(t_{k})}\wedge \by_{\psi(t_k-1)}.
\end{align*}
In order to prove \eqref{Eq det(y,y,y)}, we use the formula $\det(\bx,\by,\bz) = \Tr(J\bx J\by J\bz)$ (valid for all $\bx,\by,\bz\in\RR^3$; see formula $(2.1)$ of \cite{roy2004approximation}).
If $k\geq 1$, by noticing that $J\bx J \bx = -\det(\bx) \Id$ for all $\bx\in\RR^3$ (cf \cite{roy2004approximation}) we get
\begin{align*}
\det(\by_{t_k}, \by_{t_k+1}, \by_{t_k-1})
& = -\det(\by_{t_k+1})\Tr(J \by_{t_k} \by_{t_k+1}^{-1}\by_{\psi(t_{k+1})}) \\
& = -\det(\bw_k\by_{t_k})\Tr\big(J \by_{t_k} (\bw_k\by_{t_k})^{-1}(\bw_kN_{k+1})\big) \\
& = -\det(\bw_k)\det(\by_{t_k})\Tr(J N_{k+1}).
\end{align*}
We must eventually show \eqref{Eq wedge des z_j}. Let $k,l$ be integers, with $k\geq 0$ and $0\leq l < s_{k+1}$. Equation \eqref{Eq inter sur prod vect y_l} (for the two cases $s_{k+1}-1$ and $l$) provides the identity
\[
    \by_{t_{k+1}-1}\wedge \by_{t_{k+1}} = \det(\bw_k)^{s_{k+1}-l-1}\by_{t_k+l}\wedge \by_{t_k+l+1},
\]
and by using identity $(\bx\wedge\by)\wedge(\by\wedge\bz)=\det(\bx,\by,\bz)\by$ (valid for all points $\bx,\by,\bz$ of $\RR^3$), we get
\begin{align*}
\det(\bw_{k+1})\det(\bw_k)\bz_{t_{k+1}}\wedge \bz_{t_k+l} &= (\by_{t_{k+1}-1}\wedge \by_{t_{k+1}})\wedge(\by_{t_k-1}\wedge \by_{t_k+l}) \\
&= -\det(\bw_k)^{s_{k+1}-l-1}\det(\by_{t_k-1},\by_{t_k+l},\by_{t_k+l+1})\by_{t_k+l}\\
&= -\det(\bw_k)^{s_{k+1}-1}\det(\by_{t_k-1},\by_{t_k},\by_{t_k+1})\by_{t_k+l},
\end{align*}
using \eqref{Eq inter sur prod vect y_l bis} with ($\bx=\by_{t_k-1}$) for the last equality.
Finally we end the proof for \eqref{Eq wedge des z_j} by using \eqref{Eq det(y,y,y)} and identity $\det(\bw_{k+1}) = \det(\bw_k)^{s_{k+1}}\det(\bw_{k-1})$.

\end{Dem}

\begin{Cor}
\label{Cor arithm Tr, det, y_i, z_j}
Suppose that $\Tr(\bw_1^l\bw_0)$ and $\det(\bw_1^l\bw_0)$ are relatively prime for $i=0,1,\dots,s_2+1$, as well as $\Tr(\bw_1)$ and $\det(\bw_1)$. Then
\begin{enumerate}
\item \label{Cor point 2 arithm Tr, det, y_i, z_j} for all $k\geq 1$, $0\leq l\leq s_{k+1}+1$, $\Tr(\bw_{k}^l\bw_{k-1})$ and $\det(\bw_{k}^l\bw_{k-1})$ are relatively prime,
\item \label{Cor point 3 arithm Tr, det, y_i, z_j} for all $k\geq 1$, $0\leq l\leq s_{k+1}+1$, the matrix $\bw_{k}^l\bw_{k-1}$ is primitive,
\item \label{Cor point 4 arithm Tr, det, y_i, z_j} the content of $\by_{j}$ divides $\det(N)$ ($j \geq -2$),
\item \label{Cor point 5 arithm Tr, det, y_i, z_j} the point $\det(\bw_2)\bz_j$ belongs to $\ZZ^3$ ($j\geq -1$), and its content divides $\det(\bw_2)^2\det(N)^2\Tr(JN)$.
\end{enumerate}

\end{Cor}

\begin{Dem}
We adapt the structure of \cite[proof of Corollary 4.2]{roy2007two} to our context.\\
We show \ref{Cor point 2 arithm Tr, det, y_i, z_j} by induction. Let $\PP(k)$ be the property `` $\Tr(\bw_k^{l}\bw_{k-1})$ and $\det(\bw_k^{l}\bw_{k-1})$ are relatively prime for $l=0,\dots,s_{k+1}+1$, as well as $\Tr(\bw_k)$ and $\det(\bw_k)$ ''. By hypothesis $\PP(1)$ is true. Suppose that $\PP(k_0)$ is true for a fixed $k_0\geq 1$. Then equations provided by $\PP(k_0)$ for cases $l=s_{k_0+1}$ and $l=s_{k_0+1}+1$ ensure that $\Tr(\bw_{k_0+1})$ and $\det(\bw_{k_0+1})$ are relatively prime, as well as $\Tr(\bw_{k_0}\bw_{k_0+1})$ = $\Tr(\bw_{k_0+1}\bw_{k_0})$ and $\det(\bw_{k_0+1}\bw_{k_0})$. By \eqref{Eq Tr mod det 1}, we have
\[
\Tr(\bw_{k_0+1}^{l}\bw_{k_0})\equiv \Tr(\bw_{k_0+1})^{l-1}\Tr(\bw_{k_0+1}\bw_{k_0}) \mod \det(\bw_{k_0+1})\quad (l\geq 2).
\]
Now it suffices to notice that for each $l\geq 2$, $\det(\bw_{k_0+1}^{l}\bw_{k_0})$ has the same prime factors than
$\det(\bw_{k_0+1})$ (and than $\det(\bw_2)$), therefore the previous equation ensures that $\Tr(\bw_{k_0+1}^{l}\bw_{k_0})$ and $\det(\bw_{k_0+1}^{l}\bw_{k_0})$ are relatively prime. Thus $\PP(k_0+1)$ is true, which completes the induction step.\\
The assertion \ref{Cor point 3 arithm Tr, det, y_i, z_j} results directly from \ref{Cor point 2 arithm Tr, det, y_i, z_j} since the content of a matrix $\bw$ always divides its trace and its determinant.\\
Let $k,l$ be integers with $k\geq 0$ and $0\leq l < s_{k+1}$. By definition, we have $\by_{t_k+l} = \bw_k^{l+1}\bw_{k-1}N_k$, thus
\[
    \by_{t_k+l}\Adj(N_k) = \det(N)\bw_k^{l+1}\bw_{k-1},
\]
where $\Adj(N_k)$ denotes the transpose of the comatrix of $N_k$. In particular, since $\bw_k^{l+1}\bw_{k-1}$ is primitive, the content of $\by_{t_k+l}$ divides $\det(N)$ (the same argument enables us to deal with cases $y_{-1}$ and $y_{-2}$) which prove \ref{Cor point 4 arithm Tr, det, y_i, z_j}.\\
Finally we have to prove \ref{Cor point 5 arithm Tr, det, y_i, z_j}. The fact that $\det(\bw_2)\bz_{t_k+l}$ belongs to $\ZZ^3$ for $k=0,1,2$ and $0\leq l < s_{k+1}$ is immediate by the definition of $\bz_{t_k+l}$ and because $\det(\bw_1)$ and $\det(\bw_0)$ divide $\det(\bw_2)$. The recurrences of \eqref{Eq recurrence z_i} imply that $\det(\bw_2)\bz_j\in\ZZ^3$ for each $j\geq 0$. Moreover, the content of $\det(\bw_2)\bz_{t_k+l}$ divides that of $\det(\bw_2)^2\bz_{t_{k+1}}\wedge\bz_{t_k+l}$, thus by \eqref{Eq wedge des z_j} and \ref{Cor point 4 arithm Tr, det, y_i, z_j}, it divides $\det(\bw_2)^2\det(N)^2\Tr(JN)$.

\end{Dem}

\section{$\psi$-Sturmian sequences with multiplicative growth}
\label{section suites psi-sturm croissance multi}

In this section we show that under a multiplicative growth hypothesis, there exists a real number $\delta\geq0$ such that $|\det(\omega_k)|\asymp\norm{\omega_k}^{\delta}$ (see Proposition \ref{Prop existence delta sur le det}). This number $\delta$ will play a crucial role to determine the Diophantine exponents associated with $\xi$ and it is needed to define proper $\psi$-number.

\begin{Def}
\label{Def croissance multiplicative}
We denote by $\norm{\bw}$ the norm of a matrix $\bw\in\Mat(\RR)$, defined as the maximum of the absolute values of its coefficients.\\
Let $(\bw_i)_{i\geq 0}$ be a $\psi$-Sturmian sequence in $\GL(\RR)$. We state that $(\bw_i)_i$ has a \textsl{multiplicative growth} if there exist two positive constants $c_1$ and $c_2$ such that for each integer $k\geq 1$ and for each $1\leq l \leq s_{k+1}+1$, we have
\begin{align}
\label{Eq Def croissance multiplicative}
    c_1\norm{\bw_k}\times\norm{\bw_k^{l-1}\bw_{k-1}} \leq \norm{\bw_k^l\bw_{k-1}} \leq c_2\norm{\bw_k}\times\norm{\bw_k^{l-1}\bw_{k-1}},
\end{align}
which could be more concisely rewritten as
\[
    \norm{\bw_k^l\bw_{k-1}} \asymp \norm{\bw_k}\times\norm{\bw_k^{l-1}\bw_{k-1}}
\]
for $k \geq 1$ and $1\leq l \leq s_{k+1}+1$.
\end{Def}

The next Lemma is implied by the proof of \cite[Lemma 5.1]{roy2007two}.

\begin{Lem}
\label{Lemme croissance multiplicative}
Let $\bw_0, \bw_1\in\GL(\RR)$ be two matrices of the form $\left(
    \begin{array}{cc}
         a & b \\
        c & d
    \end{array}
    \right)$ with $1\leq a \leq\min\{b,c\}\leq\max\{b,c\}\leq d$. Then, the $\psi$-Sturmian sequence defined by $\bw_0$ and $\bw_1$ has multiplicative growth and we may take $c_1=1$ and $c_2=2$ in \eqref{Eq Def croissance multiplicative}.
\end{Lem}

\begin{Lem}[Bugeaud-Laurent (2005)]
\label{Lemme Bugeaud, Laurent mu_k}
Let $(X_n)_{n\geq 0}$ be a sequence of positive real numbers which tends to infinity and $(s_n)_{n\geq 1}$  be a sequence of positive integers. Suppose that there exists $c\geq 1$ such that
\[
    c^{-s_{n+1}}X_n^{s_{n+1}}X_{n-1} \leq X_{n+1} \leq c^{s_{n+1}}X_n^{s_{n+1}}X_{n-1} \quad (n\geq 1).
\]
For $n\geq 0$, we define $\mu_n$ by the formula $X_{n+1} = X_n^{\mu_n}$. Then, as $n$ tends to infinity we have
\[
    \mu_n = [s_{n+1}; s_n,\dots,s_1](1+\PeO(1)).
\]

\end{Lem}

\begin{Dem}
Cf \cite[Lemma 4.1]{bugeaud2005exponentsSturmian}.

\end{Dem}

The next lemma is a technical result that we will use in order to get asymptotical information on $\psi$-Sturmian sequences with multiplicative growth. It will be used in particular for the proof of Proposition \ref{Prop existence delta général}.

\begin{Lem}
\label{Lemme X_n perso}
Let $(X_n)_{n\geq 0}$ be a sequence of positive real numbers which satisfies the hypothesis of Lemma \ref{Lemme Bugeaud, Laurent mu_k}. Then
\begin{enumerate}
\item \label{enum 0 Lemme X_n perso} The sequence $(\log(X_n))_n$ strictly increases for $n$ large enough.
\item \label{enum 1 Lemme X_n perso} We have $\liminf_{n\rightarrow+\infty}\frac{\log(X_{n+2})}{\log(X_n)} > 1.$
\item \label{enum 2 Lemme X_n perso} There exists $\lambda > 1$ such that $\log(X_n)\gg \lambda^n$.
\item \label{enum 3 Lemme X_n perso} The sequence $(1/\log(X_n))_n$ is a summable sequence and as $n$ tends to infinity
\[
    \sum_{k\geq n}\frac{1}{\log(X_k)} \asymp \frac{1}{\log(X_n)}.
\]
\end{enumerate}
\end{Lem}

\begin{Dem}
First we prove \ref{enum 0 Lemme X_n perso}. Since the sequence $(\log(X_n))_n$ tends to $+\infty$, there exists $N\geq 1$ such that for each $n\geq N$, $\log(X_{n-1})+\log(c) > 0$, where $c$ is the constant involved in Lemma \ref{Lemme Bugeaud, Laurent mu_k}. For all $n\geq N$, by the growth hypothesis on $(X_n)_n$ we have
\begin{align*}
    \log(X_{n+1}) & \geq s_{n+1}\big(\log(X_n)+\log(c)\big)+\log(X_{n-1}) > \log(X_n),
\end{align*}
with the choice of $N$. Thus the sequence $(\log(X_n))_n$ strictly increases starting index $N$.\\
For the assertion \ref{enum 1 Lemme X_n perso} we define $\mu_n$ by the formula $\log(X_{n+1}) = \mu_n\log(X_n)$. According to Lemma \ref{Lemme Bugeaud, Laurent mu_k} we have $\mu_n = (1+o(1))[s_{n+1};s_n,\dots,s_1]$. Then, we conclude by noticing that

\begin{align*}
\frac{\log(X_{n+2})}{\log(X_n)} = \mu_{n+1}\mu_n \sim [s_{n+2};s_{n+1},\dots,s_1][s_{n+1};s_n,\dots,s_1] \geq \frac{3}{2}.
\end{align*}
We may deduce from this the existence of $\lambda>1$ such that
\begin{equation}
\label{eq inter Lemme X_n perso}
    \liminf_{n\rightarrow+\infty}\log(X_{n+2})/\log(X_n) > \lambda^2 > 1,
\end{equation}
and such a real $\lambda$ satisfies $\log(X_{n}) \gg (\lambda^2)^{n/2} = \lambda^n$, hence \ref{enum 2 Lemme X_n perso}. In particular, the sequence $(1/\log(X_n))_n$ is summable. Set
\[
    R_n = \sum_{k\geq n}\frac{1}{\log(X_k)} < +\infty.
\]
For $n$ large enough we have the trivial lower bound $R_n\geq 1/\log(X_n)$. On the other hand, for $n$ large enough and choosing $\lambda$ satisfying \eqref{eq inter Lemme X_n perso} we have
\begin{align*}
    R_n \leq \frac{1}{\log(X_n)} + \frac{1}{\log(X_{n+1})} + \sum_{k\geq n}\frac{1}{\lambda^2\log(X_{k})} \leq \frac{2}{\log(X_n)} + \frac{1}{\lambda^2}R_n,
\end{align*}
hence $(1-\frac{1}{\lambda^2})R_n \leq 2/\log(X_n)$, which concludes the proof of the assertion \ref{enum 3 Lemme X_n perso}.

\end{Dem}

The following property is crucial to build the $3$-system of Section \ref{Section 3-système (partiel) d'un nombre psi sturmien} (cf Figure \ref{figure 3system_partiel2.png}) and to compute the exponents $\hlambda_2,\lambda_2,\homega_2$ and $\omega_2$ associated to a $\psi$-Sturmian number defined in Proposition \ref{Prop existence y = (1,xi,xi^2)}.

\begin{Prop}
\label{Prop existence delta général}
Let $(\bw_k)_{k\geq 0}$ be a $\psi$-Sturmian sequence in $\GL(\RR)$ with multiplicative growth such that $(\norm{\bw_k})_k$ tends to infinity, and let $(D_k)_{k\geq 0}$ be a sequence of real numbers $\geq 1$ satisfying the recurrence
\begin{equation}
\label{eq rec suite réels sturmienne}
    D_{k+1} = D_k^{s_{k+1}}D_{k-1}\quad\textrm{for all $k\geq 1$.}
\end{equation}
Then, there exists a real number $\delta\geq 0$ such that
\[
    D_k \asymp \norm{\bw_k}^\delta.
\]
In particular there exists a sequence $(\hW_k)_{k\geq 0}$ satisfying \eqref{eq rec suite réels sturmienne} and such that
\[
    \hW_k \asymp \norm{\bw_k}.
\]
\end{Prop}

\begin{Dem}
For each $k\geq 0$, we write $W_k = \norm{\bw_k}$ and we choose a sequence $(D_k)_{k\geq 0}$ of real numbers $\geq 1$ satisfying the recurrence \eqref{eq rec suite réels sturmienne}. Omitting a finite number of initial terms if necessary, we may suppose that $W_k \geq 2$ for each $k$. By multiplicative growth, the sequence $(W_k)_k$ satisfies the hypothesis of Lemma \ref{Lemme X_n perso}. First, we prove that
\begin{equation}
\label{eq inter 1 Lemme existence delta général}
    \log(D_k) = \GrO(\log(W_k)).
\end{equation}
By multiplicative growth, there is a constant $M>0$ such that
\begin{equation}
\label{eq inter 3 Lemme existence delta général}
    s_{k+1}(\log(W_k)-M)+\log(W_{k-1})\leq \log(W_{k+1}),
\end{equation}
for each $k$ large enough. According to the assertion \ref{enum 2 Lemme X_n perso} of Lemma \ref{Lemme X_n perso}, there exists $k_0\geq 1$ such that
\begin{equation}
\label{eq inter 4 Lemme existence delta général}
    \big(1-\frac{1}{k}\big)\log(W_k)\leq \big(1-\frac{1}{k+1}\big)\big(\log(W_k)-M\big),
\end{equation}
for all $k\geq k_0$. Now, let $C>0$ be such that for $k = k_0-1,k_0$
\begin{equation}
\label{eq inter 2 Lemme existence delta général}
    C\log(D_k)\leq \big(1-\frac{1}{k}\big)\log(W_k).
\end{equation}
By induction we show that \eqref{eq inter 2 Lemme existence delta général} remains valid for each $k\geq k_0$. Suppose that it is true for $k-1$ and $k$. Then
\begin{align*}
    C\log(D_{k+1}) & = s_{k+1}C\log(D_k)+C\log(D_{k-1}) \\
    & \leq s_{k+1}\big(1-\frac{1}{k}\big)\log(W_k)+\big(1-\frac{1}{k-1}\big)\log(W_{k-1})\\
    & \leq s_{k+1}\big(1-\frac{1}{k+1}\big)\big(\log(W_k)-M\big)+\big(1-\frac{1}{k+1}\big)\log(W_{k-1})\\
    & \leq \big(1-\frac{1}{k+1}\big)\log(W_{k+1});
\end{align*}
indeed, the first inequality is given by the recurrence hypothesis (for $k$ and $k-1$). The second one is provided by \eqref{eq inter 4 Lemme existence delta général} and the third one by \eqref{eq inter 3 Lemme existence delta général}. Therefore \eqref{eq inter 2 Lemme existence delta général} is true for each $k\geq k_0$; this proves in particular Equation \eqref{eq inter 1 Lemme existence delta général}.\\
Now, we define $\ee_k$ by
\[
    \ee_k = \log(W_k)\times\Big(\frac{\log(D_{k+1})}{\log(W_{k+1})} - \frac{\log(D_k)}{\log(W_k)}\Big).
\]
and we shall prove that $\ee_k = \GrO(1)$. For all $k\geq 1$, by using $\log(W_{k+1})=s_{k+1}\log(W_{k})+\log(W_{k-1})+\GrO(s_{k+1})$, we obtain
\begin{align*}
\ee_k &= \frac{\log(W_k)\log(D_{k-1})-\log(D_k)\log(W_{k-1})}{\log(W_{k+1})} + \GrO\Big(\frac{s_{k+1}\log(D_k)}{\log(W_{k+1})}\Big) \\
& = -\frac{\log(W_k)}{\log(W_{k+1})}\ee_{k-1} + \GrO(1),
\end{align*}
since $s_{k+1}\log(D_k)/\log(W_{k+1}) \leq \log(D_{k+1})/\log(W_{k+1}) = \GrO(1)$ by \eqref{eq inter 1 Lemme existence delta général}. This directly implies $\ee_k = \frac{\log(W_{k-1})}{\log(W_{k+1})}\ee_{k-2} + \GrO(1)$. By the assertion \ref{enum 1 Lemme X_n perso} of Lemma \ref{Lemme X_n perso}, there are $\mu > 0$ and $k_0\geq 1$ such that
$\frac{\log(W_k)}{\log(W_{k+2})} < \mu < 1$ for each $k\geq k_0$. We deduce from the previous estimates that there is a constant $C>0$ such that $|\ee_{k+2}| \leq \mu|\ee_{k}| + C$, for all $k$. This implies that $\ee_k = \GrO(1)$ since $\mu<1$.\\
Finally, since $1/\log(W_k)$ is summable by the assertion \ref{enum 3 Lemme X_n perso} of the Lemma \ref{Lemme X_n perso}, the same goes for
\[
    \frac{\ee_k}{\log(W_k)} = \frac{\log(D_{k+1})}{\log(W_{k+1})} - \frac{\log(D_k)}{\log(W_k)},
\]
and we obtain the existence of a real number $\delta\geq 0$ such that $\log(D_k)/\log(W_k)$ tends to $\delta$ as $k$ tends to infinity. Furthermore we have
\begin{align*}
    \Big|\delta - \frac{\log(D_k)}{\log(W_k)}\Big| \leq \Big|\sum_{j\geq k} \frac{\ee_j}{\log(W_j)}\Big|
     \leq \GrO\Big(\sum_{j\geq k} \frac{1}{\log(W_j)}\Big) \leq \GrO\Big(\frac{1}{\log(W_k)}\Big).
\end{align*}
Finally, we have shown that $\log(D_k) = \delta\log(W_k) + \GrO(1)$, which proves the first part of this Lemma.\\
For the second assertion, it suffices to take a sequence $(D_k)_k$ satisfying the recurrence \eqref{eq rec suite réels sturmienne} with first terms $D_0, D_1 > 1$. Such a sequence tends to $+\infty$, and so the real number $\delta\geq 0$ previously defined cannot be zero. Then, it suffices to set $\hW_k:= D_k^{1/\delta}$.

\end{Dem}

\begin{Prop}
\label{Prop existence delta sur le det}
Let $(\bw_k)_{k\geq 0}$ be a $\psi$-Sturmian sequence in $\GL(\RR)$ with multiplicative growth such that $(\norm{\bw_k})_k$ tends to infinity. Then, there exists $0\leq\delta \leq 2$ such that
\[
    |\det(\bw_k)| \asymp \norm{\bw_k}^\delta.
\]
Let $c_1, c_2$ be the implicit constants of the multiplicative growth such that \eqref{Eq Def croissance multiplicative} is satisfied. Suppose that there exists $\alpha, \beta \geq 0$ such that the relation
\begin{align}
\label{Eq encadrement det alpha beta}
    (c_2\norm{\bw_k})^{\alpha} \leq |\det(\bw_k)| \leq (c_1\norm{\bw_k})^{\beta},
\end{align}
is true for $k=0,1$. Then this relation remains true for all $k\geq 0$; in particular we have the estimates
\[
    \alpha\leq \delta\leq \beta.
\]
\end{Prop}

\begin{Dem}
Proposition \ref{Prop existence delta général} applied with the sequence $(D_k)_{k\geq 0}=(|\det(\bw_k)|)_{k\geq 0}$ provides the existence of $\delta \geq 0$ such that $|\det(\bw_k)| \asymp \norm{\bw_k}^\delta.$ Moreover, since $|\det(\bw_k)|\leq 2 \norm{\bw_k}^2$, we always have $\delta\leq 2$.\\
In order to prove the second assertion, we follow Roy's arguments (cf \cite[Proposition 5.3]{roy2007two}) and we proceed by induction, \eqref{Eq encadrement det alpha beta} being true for $k=0,1$. Let us recall that the multiplicative growth gives
\[
(c_1\norm{\bw_{k+1}})^{s_{k+2}}\norm{\bw_{k}} \leq \norm{\bw_{k+2}} \leq (c_2\norm{\bw_{k+1}})^{s_{k+2}}\norm{\bw_{k}} \quad (k\geq 0).
\]
Suppose that \eqref{Eq encadrement det alpha beta} is true for $k=j$ and $k=j+1$ (for an index $j\geq 0$). Then we have
\begin{align*}
    |\det(\bw_{j+2})| = |\det(\bw_{j+1})|^{s_{j+2}}|\det(\bw_{j})| \leq \Big((c_1\norm{\bw_{j+1}})^{s_{j+2}}(c_1\norm{\bw_{j}})\Big)^{\beta}
    \leq (c_1\norm{\bw_{j+2}})^{\beta},
\end{align*}
and similarly $|\det(\bw_{j+2})| \geq (c_2\norm{\bw_{j+2}})^{\alpha}$. Thus \eqref{Eq encadrement det alpha beta} is valid for $k = j+2$, which concludes the induction.

\end{Dem}

\begin{Rem2}
In section \ref{Section Exemple de Roy}, we will investigate Roy's examples (cf \cite[Example 5.4]{roy2007two}) and show that they satisfy \eqref{Eq encadrement det alpha beta} (which corresponds to $(12)$ of \cite{roy2007two}).

\end{Rem2}

\section{Construction of Sturmian type numbers}
\label{section xi construction de Roy}

The main goal of this section is to establish the existence of the $\psi$-number $\xi$ (Proposition \ref{Prop existence y = (1,xi,xi^2)}) and to estimate precisely quantities $\norm{\bz_i}$, $\norm{\by_i}$, $\norm{\by\wedge\by_i}$, $|\by\cdot\bz_i|$ (where $\by=(1,\xi,\xi^2)$) in Proposition \ref{Prop estim normes yi zj}. These estimates will be used to construct the $3$-system $\bP$ of Section \ref{subsection construction partielle 3-système de xi} (see Figure \ref{figure 3system_partiel2.png}). The Diophantine exponents associated with $\xi$ are determined thanks to this $3$-system. The definition of Sturmian type number is given in Definition \ref{Def nombre de type sturmien}.

\begin{Nota}
In this section, $(\bw_i)_{i\geq 0}$ is an admissible $\psi$-Sturmian sequence in $\MM$ with multiplicative growth such that $(\norm{\bw_i})_i$ tends to $+\infty$. We denote by $N$, $(\by_i)_{i\geq -2}$ and $(\bz_i)_{i\geq -1}$ the matrix in $\GL(\RR)$ and the two sequences of symmetric matrices which are associated to $(\bw_i)_{i\geq 0}$ by Definition \ref{Def (y_i)} and by \eqref{Eq Def z_i}. We assume that $\Tr(JN)\neq 0$ (\ie $N$ is not symmetric). Finally $\delta \geq 0$ is the exponent given by Proposition \ref{Prop existence delta sur le det} and satisfying
\[
    |\det(\bw_i)| \asymp \norm{\bw_i}^{\delta}.
\]
\end{Nota}

\begin{Prop}
\label{Prop existence y = (1,xi,xi^2)}
Assume $\delta < 2$. Then there exists $\by\in\RR^3\setminus\{0\}$ such that $\det(\by) = 0$ and
\[
    \norm{\by_i\wedge \by} \ll \frac{|\det(\by_i)|}{\norm{\by_i}}.
\]
Moreover, if $\delta < 1$, the coordinates of $\by$ are linearly independent over $\QQ$ and we may assume $\by = (1,\xi,\xi^2)$ for a real number $\xi$ satisfying $[\QQ(\xi):\QQ] > 2$. Such a number $\xi$ is called a \textsl{$\psi$-Sturmian number}. If $\delta < \sigma/(1+\sigma)$, where $\sigma$ is defined by
\[
    \sigma = \frac{1}{\limsup_{k\rightarrow\infty}[s_{k+1};s_k,\dots,s_1]},
\]
and if the content of $\by_i$ is bounded, then we say that $\xi$ is a \textsl{proper $\psi$-Sturmian number}.

\end{Prop}

\begin{Def}
\label{Def nombre de type sturmien}
We say that a real number $\xi$ is of \textsl{Sturmian type} if there exists $\psi\in\FF$ such that $\xi$ is a proper $\psi$-Sturmian number.
\end{Def}

\begin{Dem}
We adapt the proof of \cite[Proposition $6.1$]{roy2007two} to our context. Let $i\geq 0$ be an integer and let us write $i=t_k+l$ with $k\geq 0$ and $0\leq l < s_{k+1}$. First, we have
\begin{align}
\label{Eq yi^y_i+1 < det(yi) * ||w_k||}
\norm{\by_i\wedge \by_{i+1}} \ll |\det(\by_i)|\times\norm{\bw_k}.
\end{align}
Indeed, as Roy noticed, one has
\begin{equation}
\label{eq technique x^y < xJy}
\norm{\bx\wedge \bz} \leq 2\norm{\bx J\bz}\quad \forall \bx,\bz\in\RR^3,
\end{equation}
since the coefficients of the diagonal of $\bx J \bz$ coincide with the first and third coefficients of
$\bx\wedge \bz$, while the sum of the coefficients of $\bx J \bz$ outside of the diagonal is the middle coefficient of $\bx\wedge \bz$ multiplied by $-1$. But according to Proposition \ref{Propriete sur y_k} \ref{Propriete 2 sur y_k} and using the fact that $\by_{i+1}$ and $\by_i$ are symmetric, we have
\begin{align*}
\norm{\by_i J\by_{i+1}} &= \norm{\by_i J (\by_i\transpose{\,\bw_k})} = \norm{\det(\by_i)J\transpose{\,\bw_k}}  = |\det(\by_i)|\espc\norm{\bw_k}
\end{align*}
(recall the formula $\bx J\bx = \det(\bx)J$ for each point $\bx\in\RR^3$). By using both the previous estimate and \eqref{eq technique x^y < xJy} with $\bx=\by_i$ and $\bz = \by_{i+1}$, we find \eqref{Eq yi^y_i+1 < det(yi) * ||w_k||}, which one may compare with $(16)$ of \cite{roy2007two}.\\
We define $\delta_i = |\det(\by_i)|/\norm{\by_i}^2$. We claim that there exists $c > 0$ such that
\begin{align}
\label{Eq dist(yi,y_i+1)< delta_i}
    \dist{\by_i}{\by_{i+1}} \leq c\delta_i
\end{align}
(one may compare this estimate with (19) of \cite{roy2007two}). Indeed, by \eqref{Eq Def y_i}, by multiplicative growth and with $c_1,c_2>0$ denoting the constants of \eqref{Eq Def croissance multiplicative}, we have
\begin{align*}
\norm{\by_{i+1}} = \norm{\bw_k^{l+2}\bw_{k-1}N_k} & \leq 2\norm{N}\espc\norm{\bw_k^{l+2}\bw_{k-1}} \\
& \leq 2c_2\norm{N}\espc\norm{\bw_k}\espc\norm{\bw_k^{l+1}\bw_{k-1}} \\
& \leq (4c_2\norm{N}\espc\norm{N^{-1}})\espc\norm{\bw_k}\espc\norm{\bw_k^{l+1}\bw_{k-1}N_k} \\
& \;\;\; = (4c_2\norm{N}\espc\norm{N^{-1}})\espc\norm{\bw_k}\espc\norm{\by_i},
\end{align*}
similarly we show that $\norm{\by_{i+1}} \geq \big(4c_1^{-1}\norm{N}\espc\norm{N^{-1}}\big)^{-1}\norm{\bw_k}\espc\norm{\by_i}$, thus
\begin{equation}
\label{Eq estimation ||y_i+1|| = ||y_i||*||w_k||}
    \norm{\by_{i+1}} \asymp\norm{\bw_k}\espc\norm{\by_i}.
\end{equation}
We complete the proof of \eqref{Eq dist(yi,y_i+1)< delta_i} by applying \eqref{Eq yi^y_i+1 < det(yi) * ||w_k||}. Furthermore, the ratio $\delta_{i+1}/\delta_i$ tends to $0$ as $i$ tends to infinity, since
\[
    \frac{\delta_{i+1}}{\delta_i} = \bigg(\frac{\norm{\by_i}}{\norm{\by_{i+1}}}\bigg)^2 \frac{|\det(\by_{i+1})|}{|\det(\by_{i})|} \asymp \frac{|\det(\bw_{k})|}{\norm{\bw_{k}}^2} \asymp \norm{\bw_{k}}^{\delta - 2} \mathop{\longrightarrow}_{k\rightarrow\infty} 0 \quad \textrm{(because $\delta < 2$).}
\]
We follow Roy's argument: by the above, there is an index $i_0 \geq 1$ such that $\delta_{i+1} \leq \delta_i/2$ for each $i\geq i_0$. Then we find for each $j >i>i_0$
\begin{align}
\label{Eq inter dist yi et yj}
    \dist{\by_i}{\by_j} \leq \sum_{m = i}^{j-1} \dist{\by_m}{\by_{m+1}} \leq c\sum_{m = i}^{j-1} \delta_m \leq 2c\delta_i.
\end{align}
Therefore the sequence $([\by_i])_{i\geq 0}$ converges in $\Proj^2(\RR)$ to a point $[\by ]$ for some non-zero $\by\in\RR^3$. Additionally, the ratio $\delta_i = |\det(\by_i)|/\norm{\by_i}^2$ depends only on the class $[\by_i]$ of $\by_i$ in $\Proj^2(\RR)$ and tends to $0$ as $i$ tends to infinity; we deduce by continuity that $|\det(\by)|/\norm{\by}^2 = 0$, \ie that $\det(\by) = 0$. By continuity and by letting $j$ tends to infinity, we also deduce from \eqref{Eq inter dist yi et yj} the estimate
\[
    \norm{\by_i\wedge\by} \ll \frac{|\det(\by_i)|}{\norm{\by_i}}.
\]
Assume now that $\delta < 1$ and let $\bu\in\ZZ^3$ be such that $\psc{\bu}{\by} = 0$. As for $\delta_i$, if we set $\delta_i' = |\det(\by_i)|/\norm{\by_i}$, then we have $\delta_{i+1}'/\delta_i' \asymp |\det(\bw_k)|/\norm{\bw_k} \asymp \norm{\bw_k}^{\delta-1}$, which tends to $0$ as $i$ tends to infinity. Thus the same goes for $\delta_i'$. By \eqref{formule technique wedge 1} we have
\[
    2\norm{\bu}\espc\norm{\by_i\wedge\by} \geq \norm{\psc{\bu}{\by}\by_i-\psc{\bu}{\by_i}\by} = |\psc{\bu}{\by_i}|\espc\norm{\by}
\]
for each $i\geq 0$. Yet, by the above, the left side tends to $0$ and the integer $\psc{\bu}{\by_i}$ is necessarily zero for $i$ large enough. This implies that $\bu = 0$ since by hypothesis $\Tr(JN)\neq 0$ implies that $\by_{t_k-1}$, $\by_{t_k}$ and $\by_{t_k+1}$ are linearly independent for each $k$ (thanks to \eqref{Eq det(y,y,y)} of Proposition \ref{Prop arithm Tr, det, y_k, z_k}). Thus, the coordinates of $\by$ are linearly independent over $\QQ$. In particular, the first coordinate of $\by$ is not zero, and by dividing $\by$ by this coordinate, we may assume that it is equal to $1$. In this case, by denoting $\xi$ the second coordinate of $\by$, the condition $\det(\by) = 0$ implies $\by = (1,\xi,\xi^2)$. The linearly independence over $\QQ$ of these three numbers ensures that $[\QQ(\xi):\QQ] > 2$.

\end{Dem}

\begin{Cor}
\label{Cor croissance ||y_i|| et |Li|}
Assume $\delta < 1$ and let $\by = (1,\xi,\xi^2)$ be the vector given by Proposition~\ref{Prop existence y = (1,xi,xi^2)}. Then
\[
    \norm{\by_i\wedge\by} \asymp \frac{|\det(\by_i)|}{\norm{\by_i}}.
\]
Moreover
\[
    \norm{\by_{i+1}}\asymp \norm{\bw_k}\espc\norm{\by_i} \quad\textrm{and}\quad \norm{\by_{i+1}\wedge\by}\asymp \norm{\bw_k}^{\delta-1}\norm{\by_i\wedge\by},
\]
with $i\geq 0$ and $k$ such that $t_k\leq i< t_{k+1}$.
In particular, as $i$ tends to infinity we have $\norm{\by_i} = o(\norm{\by_{i+1}})$ and \mbox{$\norm{\by_{i+1}\wedge\by} = o(\norm{\by_i\wedge\by})$}.

\end{Cor}

\begin{Dem}
The estimate of $\norm{\by_i\wedge\by}$ is essentially a consequence of Proposition \ref{Prop existence y = (1,xi,xi^2)} which provides the estimate $|\det(\by_i)|\gg \norm{\by_i}\espc\norm{\by_i\wedge\by}$. Write $\by_i = (y_{i,0}, y_{i,1},y_{i,2})$. We conclude in a classical manner that
\begin{align*}
\det(\by_i)
= \left| \begin{array}{cc}
y_{i,0} & y_{i,1} \\
y_{i,1} & y_{i,2}
\end{array} \right|
= \left| \begin{array}{cc}
y_{i,0} & y_{i,1} - \xi y_{i,0} \\
y_{i,1} & y_{i,2} - \xi y_{i,1}
\end{array} \right|
= y_{i,0}(y_{i,2} - \xi y_{i,1}) - y_{i,1}(y_{i,1} - \xi y_{i,0}),
\end{align*}
which gives $|\det(\by_i)| \ll \norm{\by_i}\espc\norm{\by_i\wedge\by}$, hence the announced result.\\
Equation \eqref{Eq estimation ||y_i+1|| = ||y_i||*||w_k||} yields directly  $\norm{\by_{i+1}}\asymp \norm{\bw_k}\espc\norm{\by_i}$, and using the relations $\det(\by_{i+1}) = \det(\bw_k)\det(\by_i)$, $|\det(\bw_k)|\asymp \norm{\bw_k}^{\delta}$ and the first part of the corollary, we find the last estimate: $\norm{\by_{i+1}\wedge\by}\asymp \norm{\bw_k}^{\delta-1}\norm{\by_i\wedge\by}$.

\end{Dem}

\begin{Lem}
\label{Lemme lien norme wedge et inverse}
Assume that the sequence $(\by_i)_{i\geq -2}$ tends in $\Proj^2(\RR)$ to $\by = (1,\xi,\xi^2)$ ($\xi\neq 0$) and that
\begin{equation}
\label{Eq hypo petit o Lemme lien norme wedge et inverse}
    \norm{\by_i} = o(\norm{\by_{i+1}})\quad\textrm{and}\quad \norm{\by_{i+1}\wedge\by} = o(\norm{\by_i\wedge\by})
\end{equation}
as $i$ tends to infinity. For each $i,j$ large enough, $i>j$, we have
\begin{equation}
\label{Eq ||y_i*y_j^-1||}
\frac{\norm{\by_i\wedge\by_j}}{|\det(\by_j)|} \asymp \norm{\by_i\by_j^{-1}} \asymp \norm{\by_i}\frac{\norm{\by_j\wedge\by}}{|\det(\by_j)|} \asymp \frac{\norm{\by_i}}{\norm{\by_j}}.
\end{equation}

\end{Lem}

\begin{Dem}
The last estimate of \eqref{Eq ||y_i*y_j^-1||} is a direct consequence of Corollary \ref{Cor croissance ||y_i|| et |Li|}.
For each $i\geq -2$, write $\by_i = (y_{0,i},y_{1,i},y_{2,i})$ and set
\[
    \Li{i}{1} = y_{0,i}\xi - y_{1,i}\quad \textrm{and}\quad \Li{i}{2} = y_{0,i}\xi^2 - y_{2,i}.
\]
Note that $y_{1,i}\xi - y_{2,i} =  \Li{i}{2} - \xi\Li{i}{1}$ and that $\norm{\by_i\wedge\by}\asymp\max(|\Li{i}{1}|,|\Li{i}{2}|)$. Write $\bw := \det(\by_j)\by_i\by_j^{-1} =
\left(
\begin{array}{cc}
w_{0,0} & w_{0,1} \\
w_{1,0} & w_{1,1}
\end{array}\right)$.
We have
\begin{align}
\bw & = \left(
\begin{array}{cc}
y_{0,i}y_{2,j} - y_{1,i}y_{1,j} & y_{1,i}y_{0,j} - y_{0,i}y_{1,j} \\
y_{1,i}y_{2,j} - y_{2,i}y_{1,j} & y_{2,i}y_{0,j} - y_{1,i}y_{1,j}
\end{array}
\right)\notag\\
& \label{eq inter coeff y_i*y_j^-1} = \left(
\begin{array}{cc}
y_{0,i}\big(\xi\Li{j}{1}-\Li{j}{2}\big) + y_{1,j}\Li{i}{1} & y_{0,i}\Li{j}{1} - y_{0,j}\Li{i}{1} \\
y_{1,i}\big(\xi\Li{j}{1}-\Li{j}{2}\big) + y_{1,j}\big(\Li{i}{2}-\xi\Li{i}{1}\big) & y_{1,i}\Li{j}{1}-y_{0,j}\big(\Li{i}{2}-\xi\Li{i}{1}\big)
\end{array}
\right),
\end{align}
the last equality being obtained by rewriting the coefficients of $\bw$ using determinants, for instance
\begin{align*}
w_{0,0} =
\left|
\begin{array}{cc}
y_{0,i} & y_{1,j} \\
y_{1,i} & y_{2,j}
\end{array}
\right|
= \left|
\begin{array}{cc}
y_{0,i} & y_{1,j} \\
y_{1,i}-\xi y_{0,i} & y_{2,j} -\xi y_{1,j}
\end{array}
\right|
& = \left|
\begin{array}{cc}
y_{0,i} & y_{1,j} \\
-\Li{i}{1} & \xi\Li{j}{1} -\Li{j}{2}
\end{array}
\right|\\
& = y_{0,i}\big(\xi\Li{j}{1}-\Li{j}{2}\big) + y_{1,j}\Li{i}{1}
\end{align*}
(and we proceed in a similar manner to express the three remaining coefficients of the matrix $\bw$). Since $\norm{\by_j}\espc\norm{\by_i\wedge\by} = o(\norm{\by_i}\espc\norm{\by_j\wedge\by})$ as $j$ tends to infinity uniformly with respect to $i$ such that $i>j$, Equation \eqref{eq inter coeff y_i*y_j^-1} directly implies that
\[
    \norm{\bw} \ll \norm{\by_i}\espc\norm{\by_j\wedge\by}
\]
for $i > j$ with $i,j$ large enough.\\
In the view of the form of the coefficient $w_{0,1}$ in \eqref{eq inter coeff y_i*y_j^-1} and the growth hypothesis given by \eqref{Eq hypo petit o Lemme lien norme wedge et inverse}, we may show that if $|\xi|\espc|\Li{j}{1}| \geq |\Li{j}{2}|/2$, then there exists $c'>0$, a constant which depends only on $\xi$, such that for each $j$ large enough, we have
\[
    \norm{\bw} \geq |w_{0,1}| \geq c'\norm{\by_i}\espc\norm{\by_j\wedge\by}.
\]
Similarly, if $|\xi||\Li{j}{1}| \leq |\Li{j}{2}|/2$, considering this time the coefficient $w_{1,0}$ of \eqref{eq inter coeff y_i*y_j^-1} we may show that the existence of $c''>0$, a constant which depends only on $\xi$, such that
\[
    \norm{\bw} \geq |w_{1,0}| \geq  c''\norm{\by_i}\espc\norm{\by_j\wedge\by}.
\]
This ensures $\norm{\bw} \gg \norm{\by_i}\espc\norm{\by_j\wedge\by}$ as $j$ tends to infinity uniformly with respect to $i$ > $j$, and thus $\norm{\bw} \asymp \norm{\by_i}\espc\norm{\by_j\wedge\by}$ under the same conditions. This is equivalent to the mid estimation of \eqref{Eq ||y_i*y_j^-1||}. Moreover according to the previous consideration, either $|w_{0,1}| \geq c'\norm{\by_i}\espc\norm{\by_j\wedge\by}$, or $|w_{1,0}| \geq c''\norm{\by_i}\espc\norm{\by_j\wedge\by}$. Thus $\norm{\bw}\asymp \max(|w_{0,1}|,|w_{1,0}|)$. But we have the identity
\[
    \by_i\wedge\by_j = (w_{1,0},w_{1,1}-w_{0,0},- w_{0,1}),
\]
 thus, since up to a sign $|w_{0,1}|$ and $|w_{1,0}|$ are two coefficients of $\by_i\wedge\by_j$, this ensures that $\norm{\by_i\wedge\by_j} \gg \norm{\bw}$. On the other hand, the previous formula of $\norm{\by_i\wedge\by_j}$ directly provides estimate $\norm{\by_i\wedge\by_j} \ll \norm{\bw} $. Finally $\norm{\by_i\wedge\by_j} \asymp \norm{\bw} = |\det(\by_j)|\espc\norm{\by_i\by_j^{-1}}$ for $j$ large enough, $i>j$, which ends the proof of this lemma.

\end{Dem}

Estimates of the next proposition were established by Roy for the Fibonacci case, although they were formulated in a slightly different way (see Proposition $6.1$ of \cite{roy2007two} and its proof). They are at the center of the construction of a $3$-system representing a $\psi$-Sturmian number (cf Figure \ref{figure 3system_partiel2.png}).

\begin{Prop}
\label{Prop estim normes yi zj}
Assume $\delta < 1$ and denote by $\by = (1,\xi,\xi^2)$ the vector given by Proposition \ref{Prop existence y = (1,xi,xi^2)}. We thus have the following estimates
\begin{enumerate}
\item \label{Enum 1 Prop estim normes yi zj} $\norm{\by_i\wedge\by} \asymp \frac{|\det(\by_i)|}{\norm{\by_i}}$,
\item \label{Enum 5 Prop estim normes yi zj} $\norm{\by_{i+1}} \asymp \norm{\by_i}^2\norm{\by_{\psi(i)}}^{-1}$,
\item \label{Enum 2 Prop estim normes yi zj} $\norm{\bz_i} \asymp \norm{\by_{\psi(i)}}$,
\item \label{Enum 3 Prop estim normes yi zj} $|\psc{\bz_i}{\by_{i+1}}| \asymp |\det(\by_i)|$,
\item \label{Enum 4 Prop estim normes yi zj} $|\psc{\bz_i}{\by}| \asymp \frac{|\det(\by_i)|}{\norm{\by_{i+1}}}$.
\end{enumerate}
\end{Prop}

\begin{Dem}
The assertion \ref{Enum 1 Prop estim normes yi zj} is provided by Corollary \ref{Cor croissance ||y_i|| et |Li|}.\\
For \ref{Enum 5 Prop estim normes yi zj}, considering \eqref{Eq y_j+1 = crochet} it suffices to apply Lemma \ref{Lemme lien norme wedge et inverse} with $(i+1,i)$ and $(i,\psi(i))$:
\begin{align*}
    \norm{\by_{i+1}}\espc\norm{\by_i}^{-1}\asymp\norm{\by_{i+1}\by_i^{-1}}
    = \norm{\by_{i}\by_{\psi(i)}^{-1}}\asymp \norm{\by_{i}}\espc\norm{\by_{\psi(i)}}^{-1}.
\end{align*}
Now, fix $k\geq 0$ and $0\leq l< s_{k+1}$. For\ref{Enum 2 Prop estim normes yi zj}, by \eqref{Eq ||y_i*y_j^-1||} (by taking the transpose), we have
\begin{align*}
|\det(\bw_k)|\espc\norm{\bz_{t_k+l}}
& \asymp |\det(\by_{\psi(t_{k+1})})|\espc\norm{\by_{\psi(t_{k+1})}^{-1}\by_{t_k+l}} \\
& \asymp |\det(\by_{\psi(t_{k+1})})|\espc\norm{N_{k+1}^{-1}\bw_{k}^l\bw_{k-1}N_k} \\
& \asymp |\det(\bw_k)|\espc\norm{\bw_{k}^l\bw_{k-1}N_k},
\end{align*}
the second estimate is obtained using \eqref{Eq y_psi(k) = w_k-1} and \eqref{Eq Def y_i}.
We can conclude noticing that $\bw_{k}^l\bw_{k-1}N_k = \by_{t_k+l-1} = \by_{\psi(t_k+l)}$ if $l\geq 1$ and $\bw_{k}^l\bw_{k-1}N_k = \by_{\psi(t_k)}$ if $l=0$.\\
Let us show \ref{Enum 3 Prop estim normes yi zj}. By \eqref{Eq inter sur prod vect y_l bis} (with $\bx=\by_{t_k-1}$) and by using identity $\psc{\bx\wedge\by}{\bz} = \det(\bx,\by,\bz)$ (valid for all $\bx,\by,\bz\in\RR^3$) we have by \eqref{Eq det(y,y,y)}
\begin{align*}
|\psc{\bz_{t_k+l}}{\by_{t_k+l+1}}| = |\det(\bw_k)|^{-1}|\det(\by_{t_k-1},\by_{t_k+l},\by_{t_k+l+1})|
&= |\det(\bw_k)|^{l-1}|\det(\by_{t_k-1},\by_{t_k},\by_{t_k+1})|\\
&\asymp |\det(\bw_k)|^{l}\det(\by_{t_k}) \asymp \det(\by_{t_k+l}).
\end{align*}
Finally for \ref{Enum 4 Prop estim normes yi zj} we use Property \eqref{formule technique wedge 1} of the wedge product. By \ref{Enum 1 Prop estim normes yi zj} and \ref{Enum 2 Prop estim normes yi zj} we thus obtain
\begin{align*}
\norm{\psc{\bz_{t_k+l}}{\by}\by_{t_k+l+1} - \psc{\bz_{t_k+l}}{\by_{t_k+l+1}}\by}\ll \norm{\bz_{t_k+l}}\espc\norm{\by_{t_k+l+1}\wedge\by}
&\ll \norm{\by_{\psi(t_k+l)}}\frac{|\det(\by_{t_k+l+1})|}{\norm{\by_{t_k+l+1}}}\\
&\ll \frac{|\det(\by_{t_k+l+1})|}{\norm{\bw_k}^2},\\
\end{align*}
since by multiplicative growth we have $\norm{\by_{t_k+l+1}} \asymp \norm{\bw_k^{l+2}\bw_{k-1}} \asymp \norm{\bw_k}^2\norm{\bw_k^{l}\bw_{k-1}} \asymp \norm{\bw_k}^2\norm{\by_{\psi(t_k+l)}}$. Yet, by \ref{Enum 3 Prop estim normes yi zj} we have $\frac{|\det(\by_{t_k+l+1})|}{\norm{\bw_k}^2} = o\Big(|\det(\by_{t_k+l})|\Big) = o\Big(\norm{\psc{\bz_{t_k+l}}{\by_{t_k+l+1}}\by}\Big)$, thus necessarily, $\norm{\psc{\bz_{t_k+l}}{\by}\by_{t_k+l+1}} \sim \norm{\psc{\bz_{t_k+l}}{\by_{t_k+l+1}}\by} \asymp |\det(\by_{t_k+l})|$ as $k$ tends to infinity.

\end{Dem}

\section{Partial $3$-system representing a $\psi$-Sturmian number}
\label{Section 3-système (partiel) d'un nombre psi sturmien}

Section \ref{subsection geom param des nb} is devoted to some reminders on Schmidt and Summerer's parametric geometry of numbers and on the notion of $3$-system (see Definition \ref{Def n-système}). In Section \ref{subsection construction partielle 3-système de xi}, under the assumption that the sequence $(s_i)_i$ associated with our Sturmian recurrence is bounded (note that all previous results remain valid without this hypothesis), we construct a $3$-system which partially represents a $\psi$-number $\xi$ (see Figure \ref{figure 3system_partiel2.png}). Propositions \ref{Prop Dessin L_j, L_j*} and \ref{Prop 3-sytème partiel} play a central role in this paper, their proof are presented in Section \ref{subsection Prop dessins}. This allows us to determine the six parametric Diophantine exponents associated with $\xi$: see Theorems \ref{Thm exposants geom param nb psi-sturmien} and \ref{Thm exposant psi_2 souligné}, our main result. For two exponents (namely $\pu_2$ and $\po_3$) we need an additional assumption on $\delta$ in order to determine their precise value (by ignoring the \textsl{gray areas} of the $3$-system representing $\xi$). Finally in Section \ref{Section zones d'incertitude} we construct some particular integer points whose trajectory goes through the gray areas. They require strong conditions for the gray areas to be maximal.

\subsection{Parametric geometry of numbers}
\label{subsection geom param des nb}
Let $n\geq 2$ be an integer and let $\bu\in\RR^n$ be a point whose coordinates are linearly independent over $\QQ$ (in the following we will have $n=3$ and $\bu=(1,\xi,\xi^2)$). In this section we quickly present Schmidt and Summerer's tools from the parametric geometry of numbers (cf \cite{Schmidt2009} and \cite{Schmidt2013}) and the main result due to Roy \cite{Roy_juin}, \cite{Roy_octobre} (namely Theorem \ref{Thm Roy conjecture S&S} below which settles a conjecture of Schmidt and Summerer). We recall the notion of $n$-systems in Definition \ref{Def n-système}.

\begin{Nota}
Let $I$ be an unbounded set of non negative real numbers, and let $(K_q)_{q\in I}$, $(C_q)_{q\in I}$ be two families of convex bodies of $\RR^n$ indexed by $I$. We note $K_q \asymp C_q$ and we say that $(K_q)_{q\in I}$ and $(C_q)_{q\in I}$ are \textsl{equivalent} if there exists $c>0$ such that for each index $q$ large enough we have
\[
    \frac{1}{c}K_q\subset C_q \subset cK_q.
\]
In the following, the letter $q$ will always denotes a positive real number. For all points $\bx,\by\in \RR^n$ we denote by $\bx\cdot\by = \psc{\bx}{\by}$ the standard scalar product of $\bx$ and $\by$, and in this section $\norm{\cdot}$ denotes the Euclidean norm associated with the scalar product.
\end{Nota}

\begin{Rem2}
Since all norms are equivalent in finite dimension and since estimates of the type $\asymp$ remain valid up to a multiplicative constant, we may suppose that the norm of  Proposition \ref{Prop estim normes yi zj} is the norm of our choice (although the calculations were made with the norm of Definition \ref{Def croissance multiplicative}).

\end{Rem2}

One of the founding ideas of the parametric geometry of numbers is to consider a family of convex bodies parameterized by a real number $q$ and to study the successive minima associated to this family. The choice of the family may differ according to the context. In this paper we choose to consider the convex bodies family of \cite{Roy_octobre}.

\begin{Def}
\label{Def corps convexe, exposants minima...}
We set
 \[
    \mathcal{C}_{\bu}(e^q) := \{\bx\in\RR^n\;;\; \norm{\bx} \leq 1, |\bx\cdot\bu|\leq e^{-q} \}
 \]
and
 \[
    \mathcal{C}_{\bu}^*(e^q) := \{\bx\in\RR^n\;;\; \norm{\bx} \leq e^{q}, \norm{\bx\wedge\bu}\leq 1 \}.
 \]
For $j=1,\dots,n$, $\lambda_j(q)$ (resp. $\lambda_j^*(q)$) denotes the $j$-th successive minimum of the convex body $\mathcal{C}_{\bu}(e^q)$ (resp. $\mathcal{C}_{\bu}^*(e^q)$ ) with respect to the lattice $\ZZ^n$. We also define
\begin{align*}
&\CL_{j}(q) = \log \lambda_j(q)\quad &\quad &\psi_{j}(q) = \frac{\CL_{j}(q)}{q}\\
&\overline{\psi}_{j} = \limsup_{q\rightarrow\infty}\psi_{j}(q)\quad&\quad &\underline{\psi}_{j} = \liminf_{q\rightarrow\infty}\psi_{j}(q),
\end{align*}
as well as the analogous quantities $\CL_{j}^*(q)$, $\psi_{j}^*(q)$, $\overline{\psi}_{j}^*$, $\underline{\psi}_{j}^*$  associated to $\lambda_j^*(q)$. We group these successive minima $\CL_j$ (resp. $\CL_j^*$) into a single map $\bL_{\bu}=(\CL_1,\dots,\CL_n)$ (resp. $\bL_{\bu}^*=(\CL_1^*,\dots,\CL_n^*)$).
\end{Def}

\begin{Def}
We set for each $N\geq 1$
\[
\Delta_N:=\{(x_1,\dots,x_N)\in\RR^N\;;\; x_1\leq\dots\leq x_N \}
\]
and
\[
\Phi_N:\RR^N\rightarrow\Delta_N
\]
the continuous map which lists the coordinates of a point in monotone non-decreasing order.
\end{Def}

\begin{Def}
we follow \cite[§3]{Schmidt2013} and we define the \emph{combined graph} of a set of real valued functions defined on a interval $I$ to be the union of their graphs in $I\times\RR$. For a map $\bP:[c,+\infty)\rightarrow \Delta_n$ and an interval $I\subset [c,+\infty)$, we also defined the \emph{combined graph of $\bP$ on $I$} to be the combined graph of its components $P_1,\dots,P_n$ restricted to $I$.\\

\end{Def}

\begin{Def}
\label{Def L_y et L_z}
In order to draw the combined graph of the map $\bL_{\bu}$, it is useful to define for each point $\bx\in\RR^n\setminus\{0\}$ the quantity $\lambda_{\bx}(q)$ (resp. $\lambda_{\bx}^*(q)$) to be the smallest real number $\lambda>0$ such that $\bx\in\lambda\mathcal{C}_{\bu}(e^q)$ (resp. $\bx\in\lambda\mathcal{C}_{\bu}^*(e^q)$). Then we set
\[
    \CL_{\bx}(q) = \log(\lambda_{\bx}(q))\quad\textrm{and}\quad \CL_{\bx}^*(q) = \log(\lambda_{\bx}^*(q)).
\]
Roy calls the graph of $\CL_{\bx}$ (or of $\CL_{\bx}^*$) the \textsl{trajectory} of $\bx$.
Locally, the combined graph of $\bL_{\bu}$ is included in the combined graph of a finite set of $\mathrm{L}_{\bx}$, and for each $\bx\neq 0$ we have
\begin{align*}
    &\CL_{\bx}(q) = \max\big(\log(\norm{\bx}), \log(\norm{\bx\cdot\bu})+q\big), \\
    &\CL_{\bx}^*(q) = \max\big(\log(\norm{\bx\wedge\bu}), \log(\norm{\bx})-q\big).
\end{align*}
\end{Def}

Note that
\[
    \CL_1(q) = \min_{\bx\neq 0}\CL_{\bx}(q)\quad \textrm{and}\quad \CL_1^*(q) = \min_{\bx\neq 0}\CL_{\bx}^*(q).
\]

\begin{Def}
Following \cite{davenport1967approximation} and \cite{davenport1969approximation}, we say that a point $\bx\neq 0$ is a \textsl{minimal point} (with respect to the family $\big(\mathcal{C}_{\bu}(e^q)\big)_{q\geq 0}$, resp. to the family $\big(\mathcal{C}_{\bu}^*(e^q)\big)_{q\geq 0}$) if $\bx$ is such that there is $q\geq 0$ for which
\[
    \CL_1(q) = \CL_{\bx}(q)\quad (\textrm{resp. } \CL_1^*(q) = \CL_{\bx}^*(q)).
\]
If there is no ambiguity we may not specify the considered family of convex bodies.
\end{Def}

\begin{Prop}[Mahler]
\label{Prop Mahler}
If $K_q$ denotes the dual convex body of $\mathcal{C}_{\bu}(e^q)$, then $K_q \asymp \mathcal{C}_{\bu}^*(e^q)$, and this implies that for each $j=1,\dots,n$, we have
\[
    \mathrm{L}_{j}(q) = -\mathrm{L}_{n+1-j}^*(q)+\GrO(1)
\]
as $q$ tends to infinity.
\end{Prop}

\begin{Prop}
\label{Prop Lu}
Functions $\mathrm{L}_{j}$ are continuous, piecewise linear with slopes $0$ and $1$, and satisfy the following properties
\begin{enumerate}
\item $\mathrm{L}_{1}(q)\leq\dots\leq \mathrm{L}_{n}(q)$.
\item \label{enum Prop Lu 2} $\mathrm{L}_{1}(q)+\dots+\mathrm{L}_{n}(q) = q + \GrO(1)$ as $q$ tends to infinity.
\end{enumerate}
In a dual manner, functions $\mathrm{L}_{j}^*$ are continuous, piecewise linear with slopes $0$ and $-1$, and satisfy
\begin{enumerate}
\item $\mathrm{L}_{1}^*(q)\leq\dots\leq \mathrm{L}_{n}^*(q)$.
\item $\mathrm{L}_{1}^*(q)+\dots+\mathrm{L}_{n}^*(q) = -q + \GrO(1)$ as $q$ tends to infinity.
\end{enumerate}
\end{Prop}

Schmidt and Summerer describe precisely the behavior of components $\CL_j$ by introducing the model of $(n,0)$-systems in \cite{Schmidt2013}. In \cite{Roy_octobre} Roy gives the following equivalent definition.

\begin{Def}
\label{Def n-système}
Fix a real number $q_0\geq0$. A $n$-system (or $(n,0)$-system) on  $[q_0,+\infty)$ is a continuous piecewise linear map $\bP=(\CP_1,\dots,\CP_n):[q_0,+\infty)\rightarrow\RR^n$  with the following properties:
\begin{enumerate}
\item \label{Def n-système condition 1} For each $q\geq q_0$, we have $0\leq \CP_1(q)\leq\dots\leq \CP_{n}(q)$ and $\CP_1(q)+\dots+\CP_n(q) = q$,
\item \label{Def n-système condition 2} If $H$ is a non-empty open subinterval of $[q_0,+\infty)$ on which $\bP$ is differentiable, then there is an integer $r$ ($1\leq r\leq n$), such that $\CP_r$ has slope $1$ on $H$ while the other components $\CP_j$ of $\bP$ ($j\neq r$) are constant on $H$.
\item \label{Def n-système condition 3} If $q>q_0$ is a point at which $\bP$ is not differentiable and if the integers $r$ and $s$ for which $\CP_r$ has slope $1$ on $(q-\ee,q)$ and $\CP_s$ has slope $1$ on $(q,q+\ee)$ (for  $\ee>0$ small enough) satisfy $r<s$, then we have $\CP_r(q)=\CP_{r+1}(q)=\dots=\CP_s(q)$.
\end{enumerate}
Given a subset $A$ of $\RR$, we call \textsl{interval of $A$} any interval of $\RR$ included in $A$.
Here, the condition `` $\bP$ is piecewise linear '' means that for all bounded intervals $J\subset [q_0,+\infty)$, the intersection of $J$ with the set $\textrm{D}$ of points in $[q_0,+\infty)$ at which $\bP$ is not differentiable is finite, and that the derivative of $\bP$ is locally constant on $[q_0,+\infty)\setminus D$. The slope of a component $\CP_j$ of $\bP$ on a non empty open interval $H$ of $[q_0,+\infty)\setminus D$ is the constant value of its derivative on $H$, or equivalently the slope of its graph over $H$.

\end{Def}

\begin{Def}
\label{Def n-système dual}
Fix a real number $q_0\geq0$. A dual $n$-system on $[q_0,+\infty)$ is a map $\bP:[q_0,+\infty)\rightarrow\RR^n$ such that $-\bP$ is a $n$-system on $[q_0,+\infty)$.
\end{Def}

\begin{Thm}[Schmidt and Summerer, $2013$]
\label{Thm Schmidt et S n-système}
For each non-zero point $\bu\in\RR^n$, there exist $q_0>0$ and a $n$-system $\bP$ on $[q_0,+\infty)$ such that $\norm{\bL_{\bu}-\bP}_{\infty}$ is bounded over $[0,+\infty)$.

\end{Thm}

\begin{Thm}[Roy, $2015$]
\label{Thm Roy conjecture S&S}
For each $n$-system $\bP$ on a interval $[q_0,+\infty)$, there exists a non-zero point $\bu\in\RR^n$ such that $\norm{\bL_{\bu}-\bP}_{\infty}$ is bounded on $[0,+\infty)$.

\end{Thm}

\begin{Def}
\label{Def n-system représentant xi}
Let $\xi$ be a real number and write $\bu=(1,\xi,\dots,\xi^n)$. We say that a $n$-system $\bP$ on $[q_0,+\infty[$ (resp. a dual $n$-system $\bP$ on $[q_0,+\infty[$)\textsl{ represents $\xi$}, or is a \textsl{representant of $\xi$}, if it satisfies that $\norm{\bL_{\bu}-\bP}_{\infty}$ (resp. $\norm{\bL_{\bu}^*-\bP}_{\infty}$) is bounded on $[q_0,+\infty[$.\\
Theorem \ref{Thm Schmidt et S n-système} of Schmidt and Summerer ensures that there always exists a $n$-system which represents $\xi$.

\end{Def}

\subsection{Construction of a partial $3$-system representing $\xi$}
\label{subsection construction partielle 3-système de xi}

In this section, under the assumption that the sequence $(s_i)_i$ associated with our Sturmian recurrence is bounded (note that all previous results remain valid without this hypothesis), we construct a $3$-system which partially represents a $\psi$-Sturmian number $\xi$ (see Figure \ref{figure 3system_partiel2.png}). Propositions \ref{Prop Dessin L_j, L_j*} and \ref{Prop 3-sytème partiel} play a central role in this paper, their proof are presented in Section \ref{subsection Prop dessins}. We deduce from Proposition \ref{Prop 3-sytème partiel} our main result: Theorem \ref{Thm exposants geom param nb psi-sturmien}.

\begin{Nota}
We keep notations (and hypotheses) of Section \ref{section xi construction de Roy}: $(\bw_i)_{i\geq 0}$ is an admissible $\psi$-Sturmian sequence in $\MM$ with multiplicative growth such that $(\norm{\bw_i})_i$ tends to $+\infty$. We denote by $N$, $(\by_i)_{i\geq -2}$ and $(\bz_i)_{i\geq -1}$ the matrix in $\GL(\RR)$ and the two sequences of symmetric matrices which are associated to $(\bw_i)_{i\geq 0}$ by Definition \ref{Def (y_i)} and by \eqref{Eq Def z_i}. We assume that $\Tr(JN)\neq 0$ (\ie $N$ is not symmetric). Finally $\delta \geq 0$ is the exponent given by Proposition \ref{Prop existence delta sur le det} and satisfying
\[
    |\det(\bw_i)| \asymp \norm{\bw_i}^{\delta}.
\]
In this section we additionally assume that $\delta < 1$.\\
Let $\by = (1,\xi,\xi^2)$ denote the vector given by Proposition \ref{Prop existence y = (1,xi,xi^2)}. In the following, $\CL_1, \CL_2, \CL_3$ and $\CL_{\bu} = (\CL_1,\CL_2,\CL_3)$ denote the functions defined in Section \ref{subsection geom param des nb} for the point $\bu=\by\in\RR^3$. Recall that we also assume that the sequence $(s_i)_{i\geq 0}$ associated to $\psi$ is bounded, so that by multiplicative growth, Proposition \ref{Prop existence delta sur le det} directly provides the following estimate
\begin{equation}
    |\det(\by_i)| \asymp \norm{\by_i}^{\delta}.
\end{equation}
We set
\begin{equation}
\label{Eq Def sigma}
    \sigma = \frac{1}{\limsup_{k\rightarrow\infty}[s_{k+1};s_k,\dots,s_1]} = \liminf_{k\rightarrow\infty}\frac{\log(\norm{\bw_{k}})}{\log(\norm{\bw_{k+1}})}.
\end{equation}
Let $(\hW_k)_{k\geq 0}$ be a sequence (provided by Proposition \ref{eq rec suite réels sturmienne}) of positive real numbers satisfying
\begin{enumerate}
\item $\hW_{k+1} = \hW^{s_{k+1}}\hW_{k-1}$ (for $k\geq 1$),
\item $\hW_k \asymp \norm{\bw_k}$.
\end{enumerate}
Let $k_0\geq 1$ be an integer such that for each $k\geq k_0-1$ we have $\hW_k> 1$ (in particular $(\log(\hW_k))_k$  is increasing  for $k\geq k_0$).

\end{Nota}

\begin{Def}
\label{Def Z_i, Y_i, E_i et E_i*}
Let $k,l$ be integers with $0\leq k$ and $0\leq l < s_{k+1}$. We define $\hY_{t_k+l}, \mathcal{E}_{t_k+l}^*, \hZ_{t_k+l}$ and $\mathcal{E}_{t_k+l}$ by
\begin{align*}
\hY_{t_k+l} &= \hW_k^{l+1}\hW_{k-1}\\
\mathcal{E}_{t_k+l}^* &= \big(\hW_k^{l+1}\hW_{k-1}\big)^{\delta-1}\\
\hZ_{t_k+l} &= \hW_k^{l}\hW_{k-1} \\
\mathcal{E}_{t_k+l} &= \hW_k^{(\delta-1)(l+1)-1}\hW_{k-1}^{\delta-1}.
\end{align*}
Note that these formulas for $\hY_{t_k+l}$ and $\mathcal{E}_{t_k+l}^*$ remain valid in the case $l=s_{k+1}$. We also have $\hZ_i = \hY_{\psi(i)}$.
\end{Def}

\begin{Prop}
\label{Prop estim normes VIA W_k de yi zj}
Let $i\geq 0$ be an integer. As $i$ tends to infinity we have:
\begin{enumerate}
\item $\log(\norm{\by_{i}}) = \log(\hY_{i}) + \GrO(1)$.
\item $\log(\norm{\by_{i}\wedge\by}) = \log(\mathcal{E}_{i}^*) + \GrO(1)$.
\item $\log(\norm{\bz_{i}}) = \log(\hZ_{i}) + \GrO(1)$.
\item $\log(|\bz_{i}\cdot\by|) = \log(\mathcal{E}_{i}) + \GrO(1)$.
\end{enumerate}
\end{Prop}

\begin{Dem}
Since $(s_i)_i$ is bounded, we obtain by multiplicative growth and by \eqref{Eq Def y_i}
\[
    \log(\norm{\by_{t_k+l}}) = (l+1)\log(\hW_k)+\log(\hW_{k-1}) + \GrO(1)
\]
for $0\leq k$ tending to infinity, with $0\leq l\leq s_{k+1}$. All assertions may be deduced from the estimates of Proposition \ref{Prop estim normes yi zj}.

\end{Dem}

\begin{Def}
\label{Def L_i et L_i*, q_i, c_i}
Let $k,l$ be integers with $0\leq k$ and $0\leq l < s_{k+1}$. Let us write $i=t_k+l$ and for $q\geq0$ let us set
\begin{align*}
    \hCL_{i}(q) &= \max\big(\log(\hZ_{i}),\log(\mathcal{E}_{i})+q\big), \\
    \hCL_{i}^*(q) &= \max\big(\log(\mathcal{E}_{i}^*),\log(\hY_{i})-q\big),
\end{align*}
and
\begin{align*}
    q_{i} = \log(\hZ_{i})-\log(\mathcal{E}_{i}) = \log(\hY_{i})-\log(\mathcal{E}_{i}^*) = (2-\delta)\log(\hY_{i}).
\end{align*}
We also define
\begin{align*}
    c_i & = \log(\hY_{i+1})-\log(\mathcal{E}_{i}^*) = q_i + \log(\hW_k),
\end{align*}
the intersection point abscissa of $\hCL_i^*$ and $\hCL_{i+1}^*$, which is also the intersection point abscissa of $\hCL_i$ and $\hCL_{\psi^{-1}(i)}$ (recall that by the definition of $\psi$, we have $\psi^{-1}(i) = i+1$ if $i < t_{k+1} - 1$, and $\psi^{-1}(t_{k+1}- 1) = t_{k+2}$). Also note that
\begin{equation}
\label{Eq valeurs aux c_i}
-\hCL_i^*(c_i) = -\hCL_{i+1}^*(c_i) = -\log(\mathcal{E}_{i}^*) \quad\textrm{and}\quad \hCL_i(c_i) = \hCL_{\psi^{-1}(i)}(c_i) = \log(\hZ_{\psi^{-1}(i)}).
\end{equation}
\end{Def}

\begin{figure}[H]
  \centering
  \includegraphics[width=13.5cm]{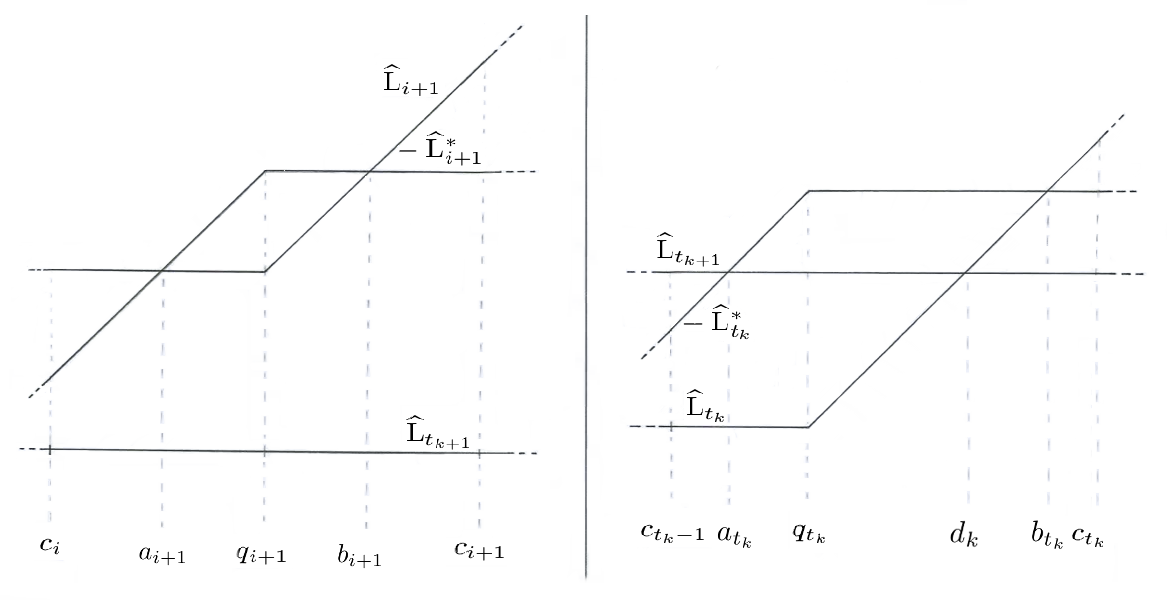}
  \caption{\label{figure 3system_partiel1.png}
  Combined graph of $\hCL_{t_{k+1}},-\hCL_{i+1}^*, \hCL_{i+1}$ on $[c_i,c_{i+1}]$}
\end{figure}

\begin{figure}[H]
  \centering
  \includegraphics[width=14cm]{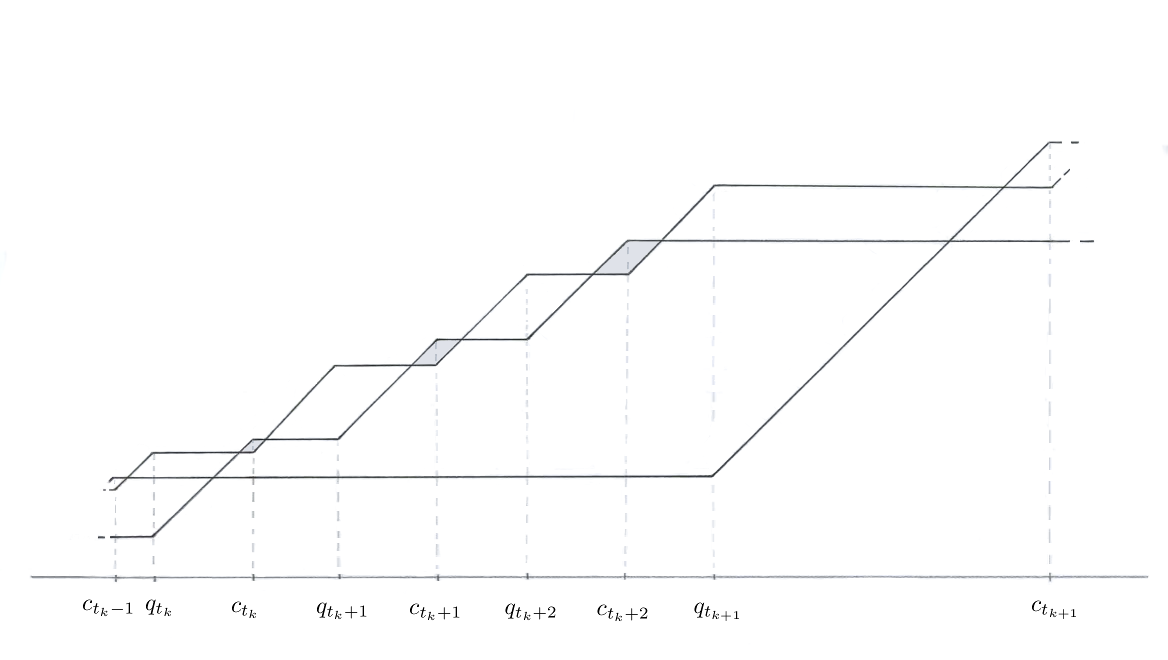}
  \caption{\label{figure 3system_partiel2.png}
  Combined graph of the $3$-system $\bP$}
\end{figure}

\begin{Prop}
\label{Prop Dessin L_j, L_j*}
We have the following properties:
\begin{enumerate}
\item \label{enum 1 Prop Dessin L_j, L_j*} For each $i\geq t_{k_0}$, we have $q_i<c_i<q_{i+1}< c_{i+1}$.
\item \label{enum 2 Prop Dessin L_j, L_j*} There exists a constant $C>0$ such that for any $i$ and for any $q>0$ we have
\[
    |\hCL_{i}(q)-\mathrm{L}_{\bz_i}(q)|\leq C\quad \textrm{and}\quad |\hCL_{i}^*(q)-\mathrm{L}_{\by_i}^*(q)|\leq C,
\]
\ie the graph of $\hCL_{i}$ (resp. $\hCL_{i}^*$) approximates the trajectory of $\bz_i$ (resp. of $\by_i$) within $\GrO(1)$, uniformly with respect to $i$.
\item \label{enum 3 Prop Dessin L_j, L_j*} The function $\hCL_{i}$ is continuous piecewise linear, constant on the interval $[0, q_i]$, and increasing with slope $1$ on the interval $[q_i,+\infty)$.
\item \label{enum 4 Prop Dessin L_j, L_j*} The function $-\hCL_{i}^*$ is continuous piecewise linear, increasing with slope $1$ on the interval $[0, q_i]$, and constant on the interval $[q_i,+\infty)$.
\item \label{enum 5 Prop Dessin L_j, L_j*} Assume that $\delta < \frac{\sigma}{1+\sigma}$. For $i = t_k+l$ with $0\leq l < s_{k+1}-1$ and $k$ large enough, the form of the combined graph of $\hCL_{t_{k+1}},-\hCL_{i+1}^*, \hCL_{i+1}$ on $[c_i,c_{i+1}]$ is similar to the drawing on the left of Figure \ref{figure 3system_partiel1.png}. The form of the combined graph of $\hCL_{t_{k+1}},-\hCL_{t_k}^*, \hCL_{t_k}$ on $[c_{t_k-1},c_{t_k}]$ is similar to the drawing on the right of Figure \ref{figure 3system_partiel1.png}.
\item\label{enum 6 Prop Dessin L_j, L_j*} Assume $\delta < \frac{\sigma}{1+\sigma}$. As $i=t_k+l$ (with $0\leq l <s_{k+1}$) tends to infinity, we have
\begin{equation}
\label{Eq Prop Dessin L_j, L_j* croissance en q_i}
    -\hCL_i^*(q_i) - \max(\hCL_{i}(q_i), \hCL_{t_{k+1}}(q_i)) \longrightarrow +\infty.
\end{equation}
\end{enumerate}
\end{Prop}

\begin{Def}
We define $\bP = \big(\CP_1,\CP_2,\CP_3\big)$ by setting for $t_{k}-1\leq i < t_{k+1}-1$ (with $k_0\leq k$), and $c_i\leq q < c_{i+1}$:
\[
    \bP(q) = \Phi_3\Big(\hCL_{t_{k+1}}(q),-\hCL_{i+1}^*(q),\hCL_{i+1}(q)\Big),
\]
and we denote by $I_{i+1} = [a_{i+1},b_{i+1}] \subset [c_i,c_{i+1}]$ the interval on which $P_3 = -\hCL_{i+1}^*$ (see Figure \ref{figure 3system_partiel1.png}). We also set $I_{i+1}'=[b_i,a_{i+1}]$.
\end{Def}

\begin{Prop}
\label{Prop 3-sytème partiel}
Recall that $\bL_{\bu} = (\CL_1,\CL_2,\CL_3)$ denotes the map of successive minima (see Definition \ref{Def corps convexe, exposants minima...} with $\bu=(1,\xi,\xi^2)$). Assume that $\delta < \frac{\sigma}{1+\sigma}$ holds and that the content of $\by_i$ is bounded (which is indeed the case if the hypotheses of Corollary \ref{Cor arithm Tr, det, y_i, z_j} are satisfied).
Then $\bP$ is a $3$-system on $[q_{t_{k_0}},+\infty)$ and its combined graph is of the type of the one given in Figure \ref{figure 3system_partiel2.png}. Moreover, there exists $C>0$ depending only on $\xi$ such that
\begin{enumerate}
\item \label{enum 1 Prop 3-sytème partiel} $|\CL_1(q) - \CP_1(q)|\leq C$ for any $q$ large enough,
\item \label{enum 2 Prop 3-sytème partiel} $|\CL_2(q) - \CP_2(q)|\leq C$ and $|\CL_3(q) - \CP_3(q)|\leq C$ for each $q\in I_j$ with $j$ large enough,
\item \label{enum 3 Prop 3-sytème partiel} $\CP_2(q) - C \leq \CL_2(q) \leq \CL_3(q) \leq \CP_3(q) + C$ for each $q\in I_j'$ with $j$ large enough.
\end{enumerate}
Roughly speaking we may say that up to a bounded difference the combined graph of $\bL_{\bu}$ coincides with that of $\bP$ outside the shaded areas on Figure \ref{figure 3system_partiel2.png}, and for each $j$ large enough, the combined graph of $\CL_2$ and $\CL_3$ on the interval $I_j'$ is included -- within $\GrO(1)$ -- in the corresponding shaded area (we also know that in this ``\textsl{gray}'' area the graphs of $\CL_2$ and $\CL_3$ will zigzag around a line with slope $1/2$ but here it is not important for our purpose).
.
\end{Prop}

The following proposition states a classical result in parametric geometry of numbers and establishes a relation between standard diophantine exponents $\hlambda_2(\xi) ,\lambda_2(\xi) ,\homega_2(\xi), \omega_2(\xi)$ and the diophantine exponents $\underline{\psi_1},\overline{\psi_1}, \underline{\psi_3},\overline{\psi_3}$ (attached to the vector $\bu=(1,\xi,\xi^2)$) which come from parametric geometry of numbers. Cf \cite{Schmidt2009} and \cite{Roy_juin}; note that $\lambda_2(\xi),\hlambda_2(\xi),\omega_2(\xi),\homega_2(\xi)$ are respectively denote by $\lambda(\bu),\hlambda(\bu),\tau(\bu),\widehat{\tau}(\bu)$ in \cite{Roy_juin}.

\begin{Prop}
\label{Prop dico exposants}
We have
\begin{align}
\label{Eq dico exposants}
\big(\underline{\psi_1},\overline{\psi_1}, \underline{\psi_3},\overline{\psi_3}\big) = \Big(\frac{1}{\omega_2(\xi)+1},\frac{1}{\homega_2(\xi)+1}, \frac{\hlambda_2(\xi)}{\hlambda_2(\xi)+1}, \frac{\lambda_2(\xi)}{\lambda_2(\xi)+1} \Big).
\end{align}
\end{Prop}

In particular, Jarn\'ik's relation \eqref{Eq Jarnik} may be rewritten as
\begin{equation}
\label{Eq Jarnik version geom param}
    2\underline{\psi}_3(\xi)+2\overline{\psi}_1(\xi)-3\underline{\psi}_3(\xi)\overline{\psi}_1(\xi)-1 = 0.
\end{equation}

\begin{Thm}
\label{Thm exposants geom param nb psi-sturmien}
Assume that $\delta<\frac{\sigma}{1+\sigma}$ and that the content of $\by_i$ is bounded. Then we have
\begin{align*}
&\underline{\psi}_1 = \frac{\sigma}{(2-\delta)(1+\sigma)},\quad \overline{\psi}_1 = \frac{1}{(1-\delta)(1+\sigma)+2}, \\
&\po_2 = \frac{1}{2+\sigma},\\
&\underline{\psi}_3 = \frac{(1-\delta)(1+\sigma)}{1+2(1-\delta)(1+\sigma)},\\
&\frac{1-\delta}{2-\delta} \leq \po_3 \leq \max\Big(\frac{1-\delta}{2-\delta},\frac{1}{2-\delta+\sigma}\Big),
\end{align*}
where $\sigma$ is defined by \eqref{Eq Def sigma} and the considered exponents are as in Definition \ref{Def corps convexe, exposants minima...} (for $\bu=(1,\xi,\xi^2)$). In particular, if $\delta$ satisfies the condition $\delta\leq h(\sigma)$ with
\begin{equation}
\label{Eq inegalité plus forte delta}
    h(\sigma) = \frac{\sigma}{2}+1-\sqrt{\big(\frac{\sigma}{2}\big)^2+1}\leq\frac{\sigma}{1+\sigma},
\end{equation}
then we have
\[
    \po_3 = \frac{1-\delta}{2-\delta}.
\]
\end{Thm}

Using \eqref{Eq dico exposants} and Roy's examples (see Section \ref{Section Exemple de Roy}), Theorem \ref{Thm exposants geom param nb psi-sturmien} directly implies Theorem \ref{Thm propre exposants intro}.

\begin{Rem}
\label{Rmq sur condition delta < sigma/(1+sigma)}
If the hypotheses of Corollary \ref{Cor arithm Tr, det, y_i, z_j} are fulfilled, then the content of $\by_i$ is bounded. Also note that if $\delta > \frac{\sigma}{1+\sigma}$, then $\bP$ (see Figure \ref{figure 3system_partiel2.png}) is no longer a $3$-system (see the proof of the assertion \ref{enum 5 Prop Dessin L_j, L_j*} of Proposition \ref{Prop Dessin L_j, L_j*}). Actually, for $\delta > \frac{\sigma}{1+\sigma}$ the expression which gives the value of $\hlambda_2(\xi)$ becomes $< 1/2$, and such value of $\hlambda_2(\xi)$ is clearly forbidden. Finally, we think that our condition $\sigma \leq h(\sigma)$ is not the best possible and may be improved for particular $\psi$-Sturmian numbers. See Section \ref{Section zones d'incertitude} for some clues supporting this conjecture.\\
For $\pu_2$ (which is the only parametric exponent not given by Theorem \ref{Thm exposants geom param nb psi-sturmien}) the situation is more complicated and we need to introduce new quantities, see Theorem \ref{Thm exposant psi_2 souligné}.
\end{Rem}

\begin{Dem}[of Theorem \ref{Thm exposants geom param nb psi-sturmien}]
To the end of calculating the exponents of Theorem \ref{Thm exposants geom param nb psi-sturmien}, we use Proposition \ref{Prop 3-sytème partiel} together with Proposition \ref{Prop Dessin L_j, L_j*} (cf Figures \ref{figure 3system_partiel1.png} and \ref{figure 3system_partiel2.png} to which we will often refer).\\
For the exponent $\underline{\psi}_1$ it suffices to calculate
\[
    \underline{\psi}_1 = \liminf_{q\rightarrow+\infty} \frac{\CL_1(q)}{q} = \liminf_{q\rightarrow+\infty} \frac{\CP_1(q)}{q} =  \liminf_{k\rightarrow\infty} \frac{\CP_1(q_{t_k})}{q_{t_k}}.
\]
Yet, in view of Definitions \ref{Def L_i et L_i*, q_i, c_i} and \ref{Def Z_i, Y_i, E_i et E_i*}, we have
\begin{align*}
    \frac{\CP_1(q_{t_k})}{q_{t_k}} = \frac{\log(\hZ_{t_k})}{q_{t_k}} &=\frac{\log(\hW_{k-1})}{(2-\delta)(\log(\hW_k)+\log(\hW_{k-1})} \\
    & =\frac{1}{(2-\delta)(1+\frac{\log(\hW_{k})}{\log(\hW_{k-1})})},
\end{align*}
which decreases with $\frac{\log(\hW_{k})}{\log(\hW_{k-1})}$, and by taking the infimum it follows $\underline{\psi}_1 = \frac{1}{(2-\delta)(1+\frac{1}{\sigma})}$.\\
For $\overline{\psi}_1$ we have
\[
    \overline{\psi}_1 = \limsup_{q\rightarrow+\infty} \frac{\CP_1(q)}{q} =  \limsup_{k\rightarrow\infty} \frac{\CP_1(d_k)}{d_k},
\]
where $d_k = (3-\delta)\log(\hW_k)+(1-\delta)\log(\hW_{k-1})$ is the intersection point abscissa of $\hCL_{t_k}$ and $\hCL_{t_{k+1}}$ (cf Figure \ref{figure 3system_partiel1.png}). Then we have
\begin{equation}
\label{eq inter calcul psi_1 surligne}
    \frac{\CP_1(d_k)}{d_k} = \frac{\log(\hZ_{t_{k+1}})}{d_k} = \frac{1}{3-\delta +(1-\delta)\frac{\log(\hW_{k-1})}{\log(\hW_k)}},
\end{equation}
and by taking the supremum we find $\overline{\psi}_1 = \frac{1}{3-\delta +(1-\delta)\sigma}$.\\
For $\underline{\psi}_3$, it may be shown that it suffices to calculate $\underline{\psi}_3  = \liminf_{k\rightarrow\infty}\frac{\CP_3(b_{t_k})}{b_{t_k}}$. We may also directly use Jarn\'ik's identity \eqref{Eq Jarnik version geom param} with the previous expression of $\overline{\psi}_1$.\\
The calculation of $\overline{\psi}_3$ is a bit more delicate because we have to show that the shaded areas of Figure \ref{figure 3system_partiel2.png} may be ignored under the condition \eqref{Eq inegalité plus forte delta}. In the view of the form of the $3$-system $\bP$ (see Figure \ref{figure 3system_partiel2.png}), we have
\begin{align*}
    \limsup_{q\rightarrow+\infty} \frac{\CP_3(q)}{q} &=  \max\Big(\limsup_{i\rightarrow\infty} \frac{\CP_3(c_i)}{c_i},\limsup_{i\rightarrow\infty}\frac{\CP_3(q_{i})}{q_{i}}\Big).
\end{align*}
On the one hand we have
\begin{align*}
    \frac{\CP_3(q_{i})}{q_{i}} = \frac{-\hCL_i^*(q_i)}{q_i} &= \frac{-\log(\mathcal{E}_i^*)}{q_i} = \frac{1-\delta}{2-\delta}.
\end{align*}
On the other hand we claim that the maximum value of $\CP_3(c_i)/c_i$ for $i=t_k+l$ with $0\leq l \leq s_{k+1}-1$ is reached for $l=s_{k+1}-1$. Indeed, in the view of Figure \ref{figure 3system_partiel1.png}, for $i=t_k+s_{k+1}-1=t_{k+1}-1$, we have
\begin{align*}
    \frac{\CP_3(c_{t_{k+1}-1})}{c_{t_{k+1}-1}} = \frac{\hCL_{t_{k+2}}(c_{t_{k+1}-1})}{c_{t_{k+1}-1}}
    = \frac{\log(\hZ_{t_{k+2}})}{q_{t_{k+1}-1}+\log(\hW_k)} = \frac{\log(\hW_{k+1})}{(2-\delta)\log(\hW_{k+1})+\log(\hW_k)}.
\end{align*}
whereas for $i<t_k+s_{k+1}-1$ we have:
\begin{align*}
\frac{\CP_3(c_i)}{c_i} = \frac{\log(\hZ_{i+1})}{q_i+\log(\hW_k)} = \frac{(l+1)\log(\hW_k)+\log(\hW_{k-1})}{(2-\delta)\big((l+1)\log(\hW_k)+\log(\hW_{k-1})\big)+\log(\hW_k)},
\end{align*}
and the right hand side of the last inequality increases with $l$, therefore it reaches its maximal value for $l=s_{k+1}-1$. Finally we have
\begin{align*}
    \limsup_{i\rightarrow\infty} \frac{\CP_3(c_i)}{c_i} = \limsup_{k\rightarrow\infty}  \frac{\CP_3(c_{t_{k+1}-1})}{c_{t_{k+1}-1}} &= \limsup_{k\rightarrow\infty}\frac{\log(\hW_{k+1})}{(2-\delta)\log(\hW_{k+1})+\log(\hW_k)}\\
    &= \frac{1}{2-\delta+\sigma},
\end{align*}
thus
\[
    \limsup_{q\rightarrow+\infty} \frac{\CP_3(q)}{q} = \max\Big( \frac{1}{2-\delta+\sigma}, \frac{1-\delta}{2-\delta}\Big).
\]
Now, note that if $\delta\leq h(\sigma)$, then $\delta^2-(\sigma+2)\delta+\sigma \geq 0$, so that we have $\max\Big( \frac{1}{2-\delta+\sigma}, \frac{1-\delta}{2-\delta}\Big) = \frac{1-\delta}{2-\delta}$. This ends the proof of the fact that under the condition \eqref{Eq inegalité plus forte delta} we have
\[
 \limsup_{q\rightarrow+\infty} \frac{\CP_3(q)}{q} = \limsup_{k\rightarrow\infty} \frac{\CP_3(q_{t_k})}{q_{t_k}} = \frac{1-\delta}{2-\delta},
\]
and since $\CP_3(q)\leq \CL_3(q)+\GrO(1)$ and $\CP_3(q_{t_k}) = \CL_3(q_{t_k})+\GrO(1)$ by Proposition \ref{Prop 3-sytème partiel}, we deduce that under this condition $\overline{\psi}_3 = \frac{1-\delta}{2-\delta}$.\\
For $\po_2$ it suffices to notice that
\[
    \po_2 = \limsup_{k\rightarrow+\infty}\frac{\CL_2(a_{t_k})}{a_{t_k}}.
\]
This follows from the fact that $\po_2\leq\po_3\leq 1/2$ and $\CL_2$ is (within $\GrO(1)$) always under the line with slope $1/2$ passing through the diagonals of the $s_{k+1}$ shaded areas located between $q_{t_k}$ and $q_{t_{k+1}}$. See Figure \ref{figure 3system_partiel2.png} (it is a consequence of Proposition \ref{Prop 3-sytème partiel} and of assertion \ref{enum Prop Lu 2} of Proposition \ref{Prop Lu}).\\
Recall that $a_{t_k}$ is the intersection point abscissa of $\hCL_{t_{k+1}}$ and $-\hCL_{t_k}^*$ (cf figure \ref{figure 3system_partiel1.png}), therefore $a_{t_k} = \log(\hZ_{t_{k+1}})+\log(\hY_{t_k}) = 2\log(\hW_k)+\log(\hW_{k-1})$. Moreover we have $\CP_2(a_{t_k}) = \hCL_{t_{k+1}}(a_{t_k}) = \log(\hZ_{t_{k+1}})$.
It thus yields
\begin{align*}
\po_2 = \limsup_{k\rightarrow+\infty}\frac{\CL_2(a_{t_k})}{a_{t_k}} = \limsup_{k\rightarrow+\infty}\frac{\CP_2(a_{t_k})}{a_{t_k}} &= \limsup_{k\rightarrow+\infty}\frac{1}{2+\frac{\log(\hW_{k-1})}{\log(\hW_k)}} = \frac{1}{2+\sigma}.
\end{align*}

\end{Dem}

\begin{Thm}
\label{Thm exposant psi_2 souligné}
Assume that $\delta < \frac{\sigma}{1+\sigma}$ and that the content of $\by_i$ is bounded. Let us define
\[
     \tau = \limsup_{k\rightarrow\infty } \frac{1}{[s_{k};s_{k-1},\dots,s_1]}\quad \textrm{and} \quad \sigma' = \liminf_{\substack{k\rightarrow\infty \\ s_{k+1}>1}} \frac{1}{[s_{k};s_{k-1},\dots,s_1]},
\]
(with the convention that $\sigma'=+\infty$ if $s_{k+1}=1$ for each $k$ large enough). Let us set
\[
    \theta(\delta) = \min\Big(\frac{1+\sigma'}{(2-\delta)(2+\sigma')},\frac{1}{2+(1-\delta)(1+\tau)}\Big).
\]
Then we have
\[
    \min\Big(\frac{(1-\delta)(1+\sigma)}{(2-\delta)(1+\sigma)+1} ,\theta(\delta)\Big) \leq \pu_2 \leq \theta(\delta).
\]
In particular, there exists a constant $c > 0$, which depends only on the Sturmian function $\psi$, such that if $\delta \leq c$, then
\[
    \pu_2 = \theta(\delta).
\]
\end{Thm}

\begin{Dem}
We claim that a point at which the function $\CP_2(q)/q$ reaches its minimum over the interval $[q_{t_k},q_{t_{k+1}}]$ belongs to $\{d_k,c_{t_k}\}$ if $s_{k+1}=1$ and belongs to $\{d_k,c_{t_k},q_{t_k+1}\}$ if $s_{k+1}>1$ (see Figures \ref{figure 3system_partiel1.png} and \ref{figure 3system_partiel2.png}).\\
Indeed, if $q$ is a minimum of $\CP_2(q)/q$ over $[q_{t_k},q_{t_{k+1}}]$, then $q$ is necessarily either $d_k$, or one of the $c_{t_k+l}$ (with $0\leq l <s_{k+1}$), or one of the $q_{t_k+l}$ (with $0<l<s_{k+1}$).\\
For points $c_{t_k+l}$ (with $0\leq l < s_{k+1}$) coming from the gray areas, we have
\begin{align*}
 \frac{\CP_2(c_{t_k+l})}{c_{t_k+l}} =  -\frac{\log(\mathcal{E}_{t_k+l}^*)}{q_{t_k+l}+\log(\hW_k)} = \frac{(1-\delta)\big((l+1)\log(\hW_k)+\log(\hW_{k-1})\big)}{(2-\delta)\big((l+1)\log(\hW_k)+\log(\hW_{k-1})\big)+\log(\hW_k)},
\end{align*}
which increases with $l$, thus being always greater than (or equal to)
\begin{align}
\label{eq inter 1 psi_2 souligné}
    \frac{\CP_2(c_{t_k})}{c_{t_k}} = \frac{(1-\delta)\Big(1+\frac{\log(\hW_{k-1})}{\log(\hW_k)}\Big)}{(2-\delta)\Big(1+\frac{\log(\hW_{k-1})}{\log(\hW_k)}\Big)+1}.
\end{align}
For points $q_{t_k+l}$, with $0<l<s_{k+1}$ and $k$ such that $s_{k+1}>1$ (there are infinitely many such $k$ if and only if $\xi$ is not of the Fibonacci type), we have
\begin{align*}
\frac{\CP_2(q_{t_k+l})}{q_{t_k+l}} = \frac{\log(\hZ_{t_k+l})}{q_{t_k+l}} &= \frac{(l+1)\log(\hW_k)+\log(\hW_{k-1})-\log(\hW_k)}{(2-\delta)\big((l+1)\log(\hW_k)+\log(\hW_{k-1})\big)}\\
& = \frac{1}{2-\delta}\Big(1- \frac{\log(\hW_k)}{(l+1)\log(\hW_k)+\log(\hW_{k-1})}\Big)
\end{align*}
which increases with $l$, thus reaching its minimal value for $l=1$. This yields
\begin{equation}
\label{eq inter 2 psi_2 souligné}
    \frac{\CP_2(q_{t_k+l})}{q_{t_k+l}} \geq \frac{\CP_2(q_{t_k+1})}{q_{t_k+1}} = \frac{1}{2-\delta}\Big(1- \frac{1}{2+\frac{\log(\hW_{k-1})}{\log(\hW_{k})}}\Big).
\end{equation}
For points $d_k$, by \eqref{eq inter calcul psi_1 surligne}, we have
\begin{equation}
\label{eq inter 3 psi_2 souligné}
    \frac{\CP_2(d_k)}{d_k} = \frac{\CP_1(d_k)}{d_k} = \frac{1}{3-\delta +(1-\delta)\frac{\log(\hW_{k-1})}{\log(\hW_k)}}.
\end{equation}
We may conclude by taking the infimum of \eqref{eq inter 1 psi_2 souligné}, \eqref{eq inter 2 psi_2 souligné} and \eqref{eq inter 3 psi_2 souligné}.\\
If $\delta=0$, then we have
\[
    \frac{(1-\delta)(1+\sigma)}{(2-\delta)(1+\sigma)+1} = \frac{1+\sigma}{3+2\sigma} > \theta(0) = \min\Big(\frac{1+\sigma'}{4+2\sigma'},\frac{1}{3+\tau}\Big),
\]
by using $\sigma',\tau \geq \sigma$. Thus, by continuity, there exists $c>0$, which depends only on $\sigma,\sigma',\tau$, such that
\[
    \min\Big(\frac{(1-\delta)(1+\sigma)}{(2-\delta)(1+\sigma)+1} ,\theta(\delta)\Big) = \theta(\delta) \quad (0\leq \delta\leq c).
\]

\end{Dem}

\subsection{Discussing the gray areas}
\label{Section zones d'incertitude}

\begin{figure}[H]
  \centering
  \includegraphics[width=14.5cm]{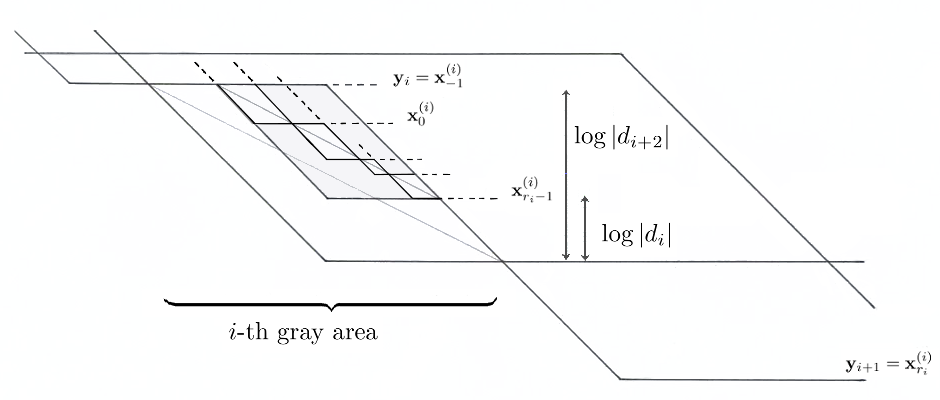}
  \caption{\label{figure zone incertitude Fibo}
  Combined graph of the $\bx_m^{(i)}$ in the $i$-th gray area}
\end{figure}

Here we discuss the condition $\delta \leq h(\sigma)$ of Theorems \ref{Thm exposants geom param nb psi-sturmien} and \ref{Thm propre exposants intro}, which may possibly be improved, at least for some class of $\psi$-Sturmian numbers, as suggested by the constructions presented below. Indeed we construct some particular integer points $\bx^{(i)}_m$ (see Definition \ref{Def x_m^(i)}) whose trajectory goes through the gray areas (see Figure \ref{figure zone incertitude Fibo}). Proposition \ref{Prop condition cas extreme gray area} implies that they require strong conditions for the gray areas to be maximal.\\\\
The question is: what does the combined graph of $\bL_{\bu}$ look like in the gray areas of the $3$-system $\bP$? (see Figure \ref{figure 3system_partiel2.png}). The optimality of condition $\delta\leq h(\sigma)$ (with $h(\sigma)$ as defined in Theorem \ref{Thm exposants geom param nb psi-sturmien}) is equivalent to the existence of $\psi$-Sturmian numbers for which the combined graph of the corresponding successive minima map $\bL_{\bu}$ coincides with the combined graph of $\bP$ at an infinite number of suitable gray areas within $o(q)$ (in other words, at the level of these suitable gray areas, the picture of the combined graph of $\bL_{\bu}$ matches precisely the outlines of the corresponding gray area within $o(q)$).\\\\
We consider in this section the dual $3$-system $-\bP$ (and the map $\bL_{\bu}^*$) rather than $\bP$ (and $\bL_{\bu}$). We will now construct some points $\bx_m^{(i)}$ whose trajectory intersects with the $i$-th gray area and which demand at least two conditions (believed to be hard to satisfy) on the sequence $(\bw_i)_i$, so that the combined graph of $\bL_{\bu}^*$ matches precisely the outlines of the $i$-th gray area. Proposition \ref{Prop condition cas extreme gray area} formulates this idea.\\\\
For an easier approach, let us focus on the Fibonacci case (although our constructions may be used to deal with the general case). In this setting, we have $s_k=1$ for each $k\geq 1$, thus $t_k=k-1$ for each $k\geq 0$. Let us fix an admissible $\psi$-Sturmian sequence $(\bw_i)_i$ with multiplicative growth and such that $(\norm{\bw_i})_i$ tends to infinity. Then, by definition $(\bw_i)_i$ satisfies
\[
    \bw_{i+1} = \bw_i\bw_{i-1}\quad (i\geq 1),
\]
and its associated sequences $(\by_i)_i$, $(\bz_i)_i$ are defined by
\[
\left\{ \begin{array}{l}
\by_i = \bw_{i+2}N_i \quad \textrm{for $i\geq -2$},\\
\bz_i = \frac{1}{\det(\bw_{i+1})}\by_{i-1}\wedge\by_i \quad \textrm{for $i\geq -1$}.
\end{array}
\right.
\]
Note that up to an index shifting, these sequences are exactly the same as Roy's in \cite{roy2007two}. For each $i\geq 0$, set
\[
    t_i = \Tr(\bw_i)\quad\textrm{and}\quad d_i=\det(\bw_i).
\]
Assume that $\delta < \frac{\sigma}{1+\sigma}=\frac{1}{\gamma^2}$, where $\delta$ is given by Proposition \ref{Prop existence delta sur le det} (recall that in the Fibonacci case we have $\sigma = \frac{1}{\gamma}$, where $\gamma$ denotes the golden ratio). Denote by $\by=(1,\xi,\xi^2)$ the vector given by Proposition \ref{Prop existence y = (1,xi,xi^2)}.
By virtue of Proposition \ref{Prop arithm Tr, det, y_k, z_k}, the sequence $(\by_i)_i$ satisfies the following recurrence

\begin{equation}
\label{Eq inter x_m^(i) Rec des y_i}
\by_{i+1} = t_{i+1}\by_{i}-d_{i+1}\by_{i-2}\quad (i\geq 0).
\end{equation}

\begin{Def}
\label{Def x_m^(i)}
Let $\big(p_m^{(i)}/q_m^{(i)}\big)_{-1\leq m\leq r_i}$ be the (finite) convergents of the rational $t_{i+1}/d_{i+1}$. By convention we set $p_{-1}^{(i)}=1$ and $q_{-1}^{(i)}=0$. We denote by $a_0^{(i)}=\lfloor t_{i+1}/d_{i+1}\rfloor,a_1^{(i)},\dots$ the partial quotients of $t_{i+1}/d_{i+1}$. For $-1\leq m\leq r_i$ we define the integer point
\begin{equation}
\label{eq def x_m^i}
    \bx_m^{(i)} = p_m^{(i)}\by_{i}-q_m^{(i)}\by_{i-2}\in\ZZ^3.
\end{equation}
\end{Def}
Note that this sequence begins with $\bx_{-1}^{(i)}=\by_{i}$ and ends with $\bx_{r_i}^{(i)}=\by_{i+1}$. It is known that $(p_m^{(i)})_m$, $(q_m^{(i)})_m$ satisfy the recurrence $u_{m+1} = a_{m+1}^{(i)}u_m+u_{m-1}$, $(m\geq 0)$ (see for instance \cite[Chapter I]{schmidt1996diophantine}). Therefore the sequence $(\bx_m^{(i)})_m$ satisfies
\begin{align}
\label{Eq recurrence x_m^i}
\bx_{m+1}^{(i)} = a_{m+1}^{(i)}\bx_m^{(i)}+\bx_{m-1}^{(i)}\quad (m\geq 0).
\end{align}
By \eqref{Eq inter x_m^(i) Rec des y_i} and using the identity $d_{i+2}=d_{i+1}d_i$, Equation \eqref{eq def x_m^i} may be rewritten as
\begin{equation}
\label{Eq inter x_m^(i) 0}
     \bx_m^{(i)} = \frac{\alpha_m^{(i)}\by_i + \beta_m^{(i)}\by_{i+1}}{d_{i+2}},
\end{equation}
with $\alpha_m^{(i)} = d_i\big(d_{i+1}p_m^{(i)}-t_{i+1}q_m^{(i)}\big)$ and $\beta_m^{(i)} = d_iq_m^{(i)}$. Furthermore, the theory of continued fractions ensures that for $m<r_i$
\begin{equation}
\label{Eq inter x_m^(i) 1}
    |\alpha_m^{(i)}| \asymp \frac{|d_{i+2}|}{q_{m+1}^{(i)}} \ll |d_{i+2}|.
\end{equation}
On the other hand we have
\begin{equation}
\label{Eq inter x_m^(i) 2}
    |\beta_m^{(i)}| = |d_iq_m^{(i)}| \leq |d_{i+2}|,
\end{equation}
since $q_m^{(i)} \leq |d_{i+1}|$. Thus, in view of the growth of $(\norm{\by_i})_i$ we deduce from \eqref{Eq inter x_m^(i) 0} and the previous estimates that
\[
    \norm{\bx_m^{(i)}} \asymp \frac{|\beta_m^{(i)}|}{|d_{i+2}|}\norm{\by_{i+1}}\ll \norm{\by_{i+1}} \quad \textrm{and} \quad \norm{\bx_m^{(i)}\wedge\by} \asymp \frac{|\alpha_m^{(i)}|}{|d_{i+2}|}\norm{\by_{i}\wedge\by} \ll\norm{\by_{i}\wedge\by},
\]
which precisely means that the trajectory of $\bx_m^{(i)}$ goes through the gray area associated with $\by_i$ and $\by_{i+1}$. In fact, thanks to \eqref{Eq inter x_m^(i) 1} and \eqref{Eq inter x_m^(i) 2}, we may prove that the trajectories of $\bx_m^{(i)}$, $m=-1,\dots,r_i$, will zigzag around a line with slope $-1/2$ (which is the vertical translate of the decreasing gray area diagonal by $\frac{1}{2}\log|d_i|$); see Figure \ref{figure zone incertitude Fibo} below.

Note that the graphs of $\CL_1^*$ and $\CL_2^*$ will zigzag around the decreasing $i$-th gray area diagonal. Therefore, the combined graph of the points $\bx_m^{(i)}$ does not match that of $\CL_1^*$ and $\CL_2^*$.

\begin{Prop}
\label{Proposition contenu des x_m^i}
Let $c_m^{(i)}$ denote the content of $\bx_m^{(i)}$. Then
\begin{enumerate}
\item The product $c_m^{(i)}c_{m+1}^{(i)}$ divides $d_{i}$ \label{Prop contenu x_m^i item 1}
\item The $\gcd$ of $c_m^{(i)}$ and $c_{m+1}^{(i)}$ divides the content of $\by_{i}$. \label{Prop contenu x_m^i item 2}
\end{enumerate}
\end{Prop}

\begin{Rem2}
In particular, if there exists a prime $p$ such that $d_i = p^{\alpha_i}$ for each $i$ (which is possible by taking $a=p$ with Roy's examples; see Section \ref{Section Exemple de Roy}), then at least half of the $\bx_m^{(i)}$ are primitive up to a bounded factor.
\end{Rem2}

\begin{Dem}
The recurrence \eqref{Eq inter x_m^(i) Rec des y_i} together with the identity $d_{i+2}=d_{i+1}d_i$ imply $\by_{i-2}\wedge\by_i = d_i\bz_{i+1}$. Also note that $\left| \begin{array}{cc}
p_m^{(i)} & p_{m+1}^{(i)} \\
q_m^{(i)} & q_{m+1}^{(i)}
\end{array} \right| = \pm 1$ by the theory of continued fractions. Hence $\bx_m^{(i)}\wedge\bx_{m+1}^{(i)} = \pm d_{i}\bz_{i+1}$, which proves assertion \ref{Prop contenu x_m^i item 1}.\\
Assertion \ref{Prop contenu x_m^i item 2} is a direct consequence of \eqref{Eq recurrence x_m^i}. If $d$ divides both the contents of $\bx_m^{(i)}$ and $\bx_{m+1}^{(i)}$, then it divides the content of all the $\bx_m^{(i)}$, especially that of $x_{-1}^{(i)} = \by_i$.

\end{Dem}

The following proposition shows that any minimal point $\bv$ which goes near the $i$-th gray area bottom outline is necessarily proportional to a $\bx_m^{(i)}$. Note that $\{\by_i,\by_{i+1}\}$ forms a basis of the linear subspace of $\RR^3$ generated by all the minimal points whose trajectory intersects with the $i$-th gray area.

\begin{Prop}
\label{Proposition point minimal proportionnel à x_m^i}
Let $\bv$ be a minimal point whose trajectory intersects with the $i$-th gray area. Let us write
\[
    \bv=\frac{\alpha\by_i+\beta\by_{i+1}}{d_{i+2}}.
\]
Let us set $\lambda = \det(\bw_2)^2\det(N)^2|\Tr(JN)|$ and suppose that $|\alpha\beta| < \frac{|d_{i+1}|}{2\lambda}$. Then $\bv$ is proportional to a point $\bx_m^{(i)}$ (with $-1\leq m \leq r_i$).

\end{Prop}

\begin{Dem}
We have $\alpha\bz_{i+1} = \bv\wedge\by_{i+1}\in\ZZ^3$ and $\beta\bz_{i+1} = \bv\wedge\by_{i}\in\ZZ^3$. According to assertion \ref{Cor point 5 arithm Tr, det, y_i, z_j} of Corollary \ref{Cor arithm Tr, det, y_i, z_j}, we have $\lambda\alpha,\lambda\beta\in\ZZ$. Let us write $\bv=(a\by_i-b\by_{i-2})/d_{i}$. Similarly, by considering $-b\bz_{i+1} = \bv\wedge\by_i\in\ZZ^3$ and $a\bz_{i+1} = \bv\wedge\by_{i-2}\in\ZZ^3$, we may show that $\lambda a,\lambda b\in\ZZ$. Also note that, following from \eqref{Eq inter x_m^(i) Rec des y_i}, we have
\[
    \bv = \frac{(ad_{i+1}-bt_{i+1})\by_i+b\by_{i+2}}{d_{i+2}},
\]
thus $\alpha = ad_{i+1}-bt_{i+1}$ and $\beta = b$.
Since $\bv$ is a minimal point we have $1\ll |\alpha\beta| \ll |d_{i+2}|$. Now assume that $|\alpha\beta| < \frac{|d_{i+1}|}{2\lambda}$. Also suppose that $b\neq 0$ (otherwise $\bv$ is proportional to $\by_i$). Under these hypothesis we have
\begin{align*}
    \Big|\frac{a'}{b'}-\frac{t_{i+1}}{d_{i+1}}\Big| < \frac{1}{2b'^2},
\end{align*}
where $a'=\lambda a,b'=\lambda b\in\ZZ$. This proves that $a'/b'$ is a convergent $p_m^{(i)}/q_m^{(i)}$ of $t_{i+1}/d_{i+1}$, and $\bv$ is proportional to $\bx_m^{(i)}$.

\end{Dem}

\begin{Prop}
\label{Prop condition cas extreme gray area}
Suppose that up to a bounded difference the combined graph of $\bL_{\bu}^*$ matches the outlines of the $i$-th gray area (see Figure \ref{figure zone incertitude Fibo}). Then
\begin{enumerate}
\item There exists $0\leq m < r_i$ such that $\bx_m^{(i)}/c_m^{(i)}$ is a minimal point and its trajectory matches the $i$-th gray area bottom outline up to a bounded difference (where $c_m^{(i)}$ denotes the content of $\bx_m^{(i)}$ ). Moreover the $(m+1)$-th partial quotient $a_{m+1}^{(i)}$ satisfies the estimate $a_{m+1}^{(i)}\asymp |d_{i+1}|$.
\item The integer $c_m^{(i)}$ satisfies the estimate $c_m^{(i)}\asymp |d_{i}|$ (it is thus maximal).
\end{enumerate}
\end{Prop}

\begin{Dem}
Suppose that the combined graph of $\bL_{\bu}^*$ matches precisely the outlines of the $i$-th gray area up to a bounded difference. This is equivalent to the existence of a minimal point $\bv$ whose trajectory matches precisely the outlines of the $i$-th gray area up to a bounded difference. If we write
\[
    \bv=\frac{\alpha\by_i+\beta\by_{i+1}}{d_{i+2}},
\]
then we have $\alpha \asymp \beta \asymp 1$. According to Proposition \ref{Proposition point minimal proportionnel à x_m^i} there exists $0\leq m < r_i$ such that $\bv$ is proportional to $\bx_m^{(i)}$ if $i$ is large enough. More precisely, $\bv = \pm \bx_m^{(i)}/c_m^{(i)}$, where $c_m^{(i)}$ denotes the content of $\bx_m^{(i)}$. In view of the trajectory of $\bx_m^{(i)}$ (see Figure \ref{figure zone incertitude Fibo}) and by Proposition \ref{Proposition contenu des x_m^i}, the content of $\bx_m^{(i)}$ has to be $\asymp |d_i|$. Moreover by using \eqref{Eq inter x_m^(i) 0} we may deduce that
\[
\big|\alpha_m^{(i)}/d_i\big| = \big|d_{i+1}p_m^{(i)}-t_{i+1}q_m^{(i)}\big| \asymp |\alpha| \asymp 1,
\]
and
\[
    \big|\beta_m^{(i)}/d_i\big| = q_m^{(i)}\asymp |\beta|\asymp 1.
\]
Finally by using the classical estimate
\[
    \Big|\frac{p_m^{(i)}}{q_m^{(i)}}-\frac{t_{i+1}}{d_{i+1}}\Big| \asymp \frac{1}{\big(q_m^{(i)}\big)^2a_{m+1}^{(i)}},
\]
we may conclude that $a_{m+1}^{(i)}\asymp |d_{i+1}|$.

\end{Dem}

Numerical tests with Roy's matrices (see Section \ref{Section Exemple de Roy}) seem to indicate that the partial quotients $a_m^{(i)}$ are small, just as the content of $\bx_m^{(i)}$ (therefore the combined graph of $\bL_{\bu}$ does not match the outlines of the gray area). A more precise study of these quantities could lead us to a better comprehension of the combined graph of $\bL_{\bu}$ in the gray areas.

\subsection{Proof of Propositions \ref{Prop Dessin L_j, L_j*} and  \ref{Prop 3-sytème partiel}}
\label{subsection Prop dessins}

In this section we essentially prove Propositions \ref{Prop Dessin L_j, L_j*} and  \ref{Prop 3-sytème partiel} of Section \ref{subsection construction partielle 3-système de xi}.

\begin{Lem}
We have the following inequality
\begin{equation}
\label{Eq inegalite sigma [s_k+1;...]}
    \limsup_{k\rightarrow+\infty} \frac{1}{[s_{k+1};s_k,\dots,s_1]} \leq \frac{1}{1+\sigma},
\end{equation}
which implies the inequality
\begin{equation}
\label{Eq inter 0 Prop Dessin L_j, L_j*}
    \frac{\sigma}{1+\sigma} \leq 1-\limsup_{k\rightarrow+\infty} \frac{\log(\hW_{k-1})}{\log(\hW_k)}.
\end{equation}
\end{Lem}
Note that \eqref{Eq inegalite sigma [s_k+1;...]} and \eqref{Eq inter 0 Prop Dessin L_j, L_j*} are equalities if the sequence $(s_k)_k$ is equal to the constant sequence $(1,1,1,\dots)$ (Fibonacci case).

\begin{Dem}
By the definition of $\sigma$ we have
\begin{align*}
    \frac{1}{[s_{k+1};s_k,\dots,s_1]} & \leq  \frac{1}{1 +\frac{1}{[s_k;s_{k-1},\dots,s_1]}} \leq  \frac{1}{1 +\sigma}+o(1),
\end{align*}
and we find the first inequality of our Lemma by taking the supremum of the left hand side. In the case $(s_n)_n = (1,1,\dots)$ a quick computation shows that equality is reached by using the formula
\[
    \limsup_{k\rightarrow+\infty} \frac{1}{[s_{k+1};s_k,\dots,s_1]} = \liminf_{k\rightarrow+\infty} \frac{1}{[s_{k+1};s_k,\dots,s_1]} = \frac{1}{\gamma}
\]
(recall that $\gamma$ denotes the golden ratio).\\
The second inequality may be deduced from the first one using Lemma \ref{Lemme Bugeaud, Laurent mu_k}.

\end{Dem}

\begin{Dem}[of Proposition \ref{Prop Dessin L_j, L_j*}]$\newline$
Let $k\geq k_0$ be an integer and let us write $i=t_k+l$ with $0\leq l < s_{k+1}$.\\
Since $\log(\hW_k) > 0$ and $c_i=q_i+\log(\hW_k)$ and $q_{i+1} = q_i+(2-\delta)\log(\hW_k)$ (with $\delta < 1$ by hypothesis), we have the estimates $q_i<c_i<q_{i+1}$, proving immediately \ref{enum 1 Prop Dessin L_j, L_j*}. The assertion \ref{enum 2 Prop Dessin L_j, L_j*} is a direct consequence of estimates of Proposition \ref{Prop estim normes VIA W_k de yi zj} and of the definition of functions $\CL_{\bz_i}$ and $\CL_{\by_i}^*$. Assertions \ref{enum 3 Prop Dessin L_j, L_j*} and \ref{enum 4 Prop Dessin L_j, L_j*} directly result from the definitions of functions $\hCL_i$ and $\hCL_i^*$.\\
Now let us prove \ref{enum 5 Prop Dessin L_j, L_j*}. First assume that $0\leq l<s_{k+1}-1$ and let us show that the form of the combined graph of $\hCL_{t_{k+1}},-\hCL_{i+1}^*, \hCL_{i+1}$ on $[c_i,c_{i+1}]$ is similar to the drawing on the left of Figure \ref{figure 3system_partiel1.png}. By \ref{enum 3 Prop Dessin L_j, L_j*} and \ref{enum 4 Prop Dessin L_j, L_j*} applied with $i+1$ and since the inequality $c_{i+1} < q_{t_{k+1}}$ implies that $\hCL_{t_{k+1}}$ is constant on $[c_i,c_{i+1}]$, it suffices to prove that
\begin{align}
   \hCL_{i+1}(q_{i+1}) & < -\hCL_{i+1}^*(q_{i+1}), \label{Eq inter 1 Prop Dessin L_j, L_j*}\\
   -\hCL_{i+1}^*(c_{i}) & \leq \hCL_{i+1}(c_{i}) \quad \textrm{and} \quad -\hCL_{i+1}^*(c_{i+1}) \leq \hCL_{i+1}(c_{i+1}), \label{Eq inter 2 Prop Dessin L_j, L_j*}\\
   \hCL_{t_{k+1}}(c_i) & < -\hCL_{i+1}^*(c_i). \label{Eq inter 3 Prop Dessin L_j, L_j*}
\end{align}
We claim that the inequality \eqref{Eq inter 1 Prop Dessin L_j, L_j*} holds when $\delta < \sigma$ (implied by our hypothesis $\delta < \frac{\sigma}{1+\sigma}$). Inequalities \eqref{Eq inter 2 Prop Dessin L_j, L_j*} are satisfied under the single hypothesis $\delta \geq 0$ (with equality if and only if $\delta = 0$). The last inequality \eqref{Eq inter 3 Prop Dessin L_j, L_j*} holds under the hypothesis $\delta < \frac{\sigma}{1+\sigma}$.\\
In order to prove \eqref{Eq inter 1 Prop Dessin L_j, L_j*} we consider the following equivalences:
\begin{align*}
 & -\hCL_{t_{k}+l+1}^*(q_{t_{k}+l+1}) > \hCL_{t_{k}+l+1}(q_{t_{k}+l+1}) \\
 & \Leftrightarrow (1-\delta)\big((l+2)\log(\hW_{k}) + \log(\hW_{k-1})\big) > (l+1)\log(\hW_{k}) + \log(\hW_{k-1}) \\
 & \Leftrightarrow \delta < \frac{\log(\hW_{k})}{(l+2)\log(\hW_{k})+\log(\hW_{k-1})}.
\end{align*}
The right hand side of the last inequality decreases with $l$, thus it is always greater than (or equal to) $\log(\hW_{k})/\log(\hW_{k+1})$ because $l \leq s_{k+1}-2$. Yet according to Lemma \ref{Lemme Bugeaud, Laurent mu_k} and by the definition of $\sigma$, we have for any $\ee > 0$
\[
    \frac{\log(\hW_{k})}{\log(\hW_{k+1})} = \frac{1}{[s_{k+1};s_{k},\dots,s_1](1+o(1))} \geq \sigma-\ee
\]
for any $k$ large enough, proving \eqref{Eq inter 1 Prop Dessin L_j, L_j*} by the equivalences we have just showed (and by using $\delta < \sigma$). Moreover, one may deduce from the calculations that $\frac{-\hCL_{i+1}^*(q_{i+1})}{\hCL_{i+1}(q_{i+1})} > 1+\ee$ for $\ee>0$ small enough, which proves that
\begin{equation}
    \label{Eq inter Prop Dessin L_j, L_j* 0 bis}
    -\hCL_{i}^*(q_{i}) - \hCL_{i}(q_{i}) \longrightarrow +\infty
\end{equation}
as $i$ tends to infinity with $i = t_k+l$ and $0<l<s_{k+1}$.\\
For \eqref{Eq inter 2 Prop Dessin L_j, L_j*}, it suffices to note that for each $t_k \leq i < t_{k+1}-1$ we have
\begin{align}
 -\hCL_{i+1}^*(c_{i}) & = (1-\delta)\hCL_{i+1}(c_{i}), \label{Eq inter 5 Prop Dessin L_j, L_j*} \\
 -\hCL_{i+1}^*(c_{i+1}) & = (1-\delta) \hCL_{i+1}(c_{i+1}), \label{Eq inter 6 Prop Dessin L_j, L_j*}
\end{align}
the equality \eqref{Eq inter 6 Prop Dessin L_j, L_j*} remaining valid for $i=t_k-1$.\\
Finally, let us prove that the hypothesis $\delta < \frac{\sigma}{1+\sigma}$ implies \eqref{Eq inter 3 Prop Dessin L_j, L_j*}. Since $c_i$ is the intersection point abscissa of $\hCL_i^*$ and $\hCL_{i+1}^*$, we have
\begin{align*}
    -\hCL_{i+1}^*(c_i) > \hCL_{t_{k+1}}(c_i) &\Leftrightarrow -\hCL_{i}^*(c_i) > \hCL_{t_{k+1}}(c_i) \\
    &\Leftrightarrow -\log(\mathcal{E}_{i}^*) > \log(\hW_k) \\
    &\Leftrightarrow 1-\frac{1}{(l+1)+\frac{\log(\hW_{k-1})}{\log(\hW_k)}} > \delta.
\end{align*}
The left hand side of the last equality decreases with $l$, therefore it reaches its minimum for $l=0$. Yet, by the definition of $\sigma$, for any $\ee > 0$ small enough we have for $k$ large enough
\[
    1-\frac{1}{1+\frac{\log(\hW_{k-1})}{\log(\hW_k)}}> 1-\frac{1}{1+\sigma} - \ee > \delta,
\]
thus allowing \eqref{Eq inter 3 Prop Dessin L_j, L_j*} to hold according to the previous equivalences. This ends the proof that the form of the combined graph of $\hCL_{t_{k+1}},-\hCL_{i+1}^*, \hCL_{i+1}$ on $[c_i,c_{i+1}]$ is similar to the drawing on the left of Figure \ref{figure 3system_partiel1.png} for the case $t_k\leq i < t_{k+1}-1$.\\
Let us determine the form of the combined graph of $\hCL_{t_{k+1}},-\hCL_{t_k}^*, \hCL_{t_k}$ on $[c_{t_k-1},c_{t_k}]$. Thanks to \ref{enum 3 Prop Dessin L_j, L_j*} and \ref{enum 4 Prop Dessin L_j, L_j*}, if we prove that
\begin{equation}
\label{Eq inter 7 Prop Dessin L_j, L_j*}
    \hCL_{t_k}(c_{t_k-1}) \leq -\hCL_{t_k}^*(c_{t_k-1}) \leq \hCL_{t_{k+1}}(c_{t_k-1})
\end{equation}
and
\begin{equation}
\label{Eq inter 8 Prop Dessin L_j, L_j*}
    \hCL_{t_{k+1}}(c_{t_k})\leq -\hCL_{t_k}^*(c_{t_k}) \leq \hCL_{t_k}(c_{t_k})
\end{equation}
hold, then it shows that the form of the combined graph is similar to the drawing on the right of Figure \ref{figure 3system_partiel1.png}. If $t_k+1<t_{k+1}$, since $-\hCL_{t_k}^*(c_{t_k}) = -\hCL_{t_{k}+1}^*(c_{t_k})$ and $\hCL_{t_{k}}(c_{t_k})=\hCL_{t_{k}+1}(c_{t_k})$, then inequalities \eqref{Eq inter 8 Prop Dessin L_j, L_j*} come from the form of the combined graph on $[c_i,c_{i+1}]$ for $i=t_k$, which has already been proved. If $t_k+1 = t_{k+1}$, since $-\hCL_{t_k}^*(c_{t_k}) = -\hCL_{t_{k}+1}^*(c_{t_k})$ and $\hCL_{t_{k}}(c_{t_k})=\hCL_{t_{k+2}}(c_{t_k})$, then inequalities \eqref{Eq inter 8 Prop Dessin L_j, L_j*} are in fact inequalities \eqref{Eq inter 7 Prop Dessin L_j, L_j*} for $k+1$. Let us prove \eqref{Eq inter 7 Prop Dessin L_j, L_j*}. We have
\begin{equation}
\label{Eq inter 9 Prop Dessin L_j, L_j*}
     -\hCL_{t_k}^*(c_{t_k-1}) = -\log(\hY_{t_k})+\log(\hY_{t_k})-\log(\mathcal{E}_{t_k-1}^*) = (1-\delta)\log(\hW_k),
\end{equation}
from which we directly deduce $ -\hCL_{t_k}^*(c_{t_k-1})\leq \hCL_{t_{k+1}}(c_{t_k-1}) (= \log(\hW_k))$, with equality if and only if $\delta=0$. Furthermore, we have $\hCL_{t_k}(c_{t_k-1}) = \log(\hW_{k-1})$, therefore the relation \eqref{Eq inter 9 Prop Dessin L_j, L_j*} shows that \eqref{Eq inter 7 Prop Dessin L_j, L_j*} results from the upper bound
\[
    \delta < 1-\frac{\log(\hW_{k-1})}{\log(\hW_{k})},
\]
which is provided by the hypothesis $\delta < \frac{\sigma}{1+\sigma}$ and by inequality \eqref{Eq inter 0 Prop Dessin L_j, L_j*}. This ends the proof of \eqref{Eq inter 7 Prop Dessin L_j, L_j*}, therefore the combined graph of  $\hCL_{t_{k+1}},-\hCL_{t_k}^*, \hCL_{t_k}$ on $[c_{t_k-1},c_{t_k}]$ has indeed the announced form. Moreover, inequalities \eqref{Eq inter 7 Prop Dessin L_j, L_j*} and \eqref{Eq inter 8 Prop Dessin L_j, L_j*} are strict inequalities if $\delta > 0$.\\
Finally let us prove the assertion \ref{enum 6 Prop Dessin L_j, L_j*}. According to the assertion \ref{enum 5 Prop Dessin L_j, L_j*} (cf Figure \ref{figure 3system_partiel1.png}) for $i=t_k+l$ with $0 < l < s_{k+1}$, we have
\[
    \max\big(\hCL_{i}(q_i), \hCL_{t_{k+1}}(q_i)\big) = \hCL_{i}(q_i),
\]
 and by \eqref{Eq inter Prop Dessin L_j, L_j* 0 bis}, we have the estimate \eqref{Eq Prop Dessin L_j, L_j* croissance en q_i} if $i$ is of the form $i = t_k+l$ with $0 < l < s_{k+1}$. Let us fix $\ee > 0$. For $i=t_k$ we have
\[
    \max\big(\hCL_{i}(q_i), \hCL_{t_{k+1}}(q_i)\big) = \hCL_{t_{k+1}}(q_i),
\]
and for $k$ large enough
\begin{align*}
\frac{-\hCL_{t_k}^*(q_{t_k})}{\hCL_{t_{k+1}}(q_{t_k})} &= \frac{(1-\delta)(\log(\hW_k) +\log(\hW_{k-1}) )}{\log(\hW_k)} \geq (1-\delta)(1+\sigma-\ee).
\end{align*}
The hypothesis $\delta < \frac{\sigma}{1+\sigma}$ is equivalent to $(1-\delta)(1+\sigma) > 1$, thus by choosing $\ee$ small enough, we obtain the existence of a constant $\mu$ such that $-\hCL_{t_k}^*(q_{t_k})/\hCL_{t_{k+1}}(q_{t_k})> \mu > 1$ for each $k$ large enough. This proves \eqref{Eq Prop Dessin L_j, L_j* croissance en q_i} for $i$ of the form $i = t_k$, which ends the proof of \ref{enum 6 Prop Dessin L_j, L_j*}.

\end{Dem}

\begin{Dem}[of Proposition \ref{Prop 3-sytème partiel}]$\newline$
First let us prove that $\bP$ is a $3$-system. We start with the continuity of $\bP$. We are reduced to verifying that $\bP$ is continuous at abscissas $c_i$. Recall that for $t_k\leq i < t_{k+1}-1$, we have
\[
    -\hCL_i^*(c_i) = -\hCL_{i+1}^*(c_i)\quad\textrm{and}\quad \hCL_i(c_i) = \hCL_{i+1}(c_i).
\]
This directly ensures the continuity of $\bP$ on the whole interval $[c_{t_k-1}, c_{t_{k+1}-1})$. Taking $i=t_{k+1}-2$, at point $c_{i+1}$ we have
\[
    -\hCL_{i+1}^*(c_{i+1}) = -\hCL_{t_{k+1}}^*(c_{i+1})\quad\textrm{and}\quad \hCL_{i+1}(c_{i+1}) = \hCL_{t_{k+2}}(c_{i+1}).
\]
Yet, by the definition of $\bP$, we have
\[
    \bP(c_{i+1}) = \Phi_3\Big(\hCL_{t_{k+1}}(c_{i+1}),\hCL_{t_{k+2}}(c_{i+1}),-\hCL_{t_{k+1}}^*(c_{i+1})\Big),
\]
finally leading to $\bP$ being continuous at this point. So the map $\bP:[q_{t_{k_0}},+\infty)$ is continuous and its components $\CP_i$ are piecewise linear with slopes $1$ or $0$. In view of assertion \ref{enum 5 Prop Dessin L_j, L_j*} of Proposition \ref{Prop Dessin L_j, L_j*} and of the combined graph of Figure \ref{figure 3system_partiel1.png}, the form of the combined graph of $\bP$ is that of Figure~\ref{figure 3system_partiel2.png}.\\
Let $k,i$ be integers with $k_0\leq k$ and $t_{k}-1\leq i < t_{k+1}-1$. Recall that $c_i< q_{i+1} < c_{i+1}$. Furthermore, we know that $\hCL_{i+1}$ is constant on $[c_i,q_{i+1}]$ and monotone increasing with slope $1$ on $[q_{i+1},c_{i+1}]$, while $-\hCL_{i+1}^*$ is monotone decreasing with slope $1$ on $[c_i,q_{i+1}]$, and constant on $[q_{i+1},c_{i+1}]$. Since $c_{i+1}<q_{t_{k+1}}$, $\hCL_{t_{k+1}}$ is constant on $[c_i,c_{i+1}]$. We may infer from this that $\bP$ satisfies the point \ref{Def n-système condition 2} of Definition \ref{Def n-système}. In particular, $\CP_1+\CP_2+\CP_3$ is affine with slope $1$.\\
Since by the definition of $\bP$ we have $\CP_1\leq\CP_2\leq \CP_3$, showing that $\bP$ fulfills the point \ref{Def n-système condition 1} of Definition \ref{Def n-système} requires to prove that $\CP_1(q)+\CP_2(q)+\CP_3(q) = q$ for $q\geq q_{t_{k_0}}$. And since $\CP_1+\CP_2+\CP_3$ is affine with slope $1$ from the above, it suffices to verify that it is fulfilled at a particular point $q$. Yet for $k\geq k_0$, we have
\begin{align*}
    \CP_1(q_{t_k})+\CP_2(q_{t_k})+\CP_3(q_{t_k})
    & = \hCL_{t_k}(q_{t_k})+\hCL_{t_{k+1}}(q_{t_k})-\hCL_{t_k}^*(q_{t_k})\\
    & = \log(\hW_{k-1}) +  \log(\hW_{k}) - (\delta-1)(\log(\hW_k) + \log(\hW_{k-1}))\\
    & = q_{t_k}.
\end{align*}
Finally, in view of the form of its combined graph given in Figure \ref{figure 3system_partiel2.png} it is clear that $\bP$ fulfills the point \ref{Def n-système condition 3} of Definition \ref{Def n-système}. As we claimed, $\bP$ is therefore a $3$-system.\\
Now let us prove the assertions on $\bL_{\bu}$ and $\bP$. First note that if $i\neq j$ (with $i$ and $j$ large enough) then $\bz_i$ and $\bz_j$ are linearly independent. If not, then there would exist a constant $\lambda$ such that $\CL_{\bz_i} = \CL_{\bz_j} + \lambda$ and it is clear in view of Figure \ref{figure 3system_partiel2.png} that this never happens. Let $i=t_k+l$ be an integer with $0\leq l < s_{k+1}$ and $k$ large enough. By the definition of the successive minima, we deduce by the above and by assertion \ref{enum 2 Prop Dessin L_j, L_j*} of Proposition \ref{Prop Dessin L_j, L_j*} that we have for $q\in I_i=[a_i,b_i]$
\begin{align}
    \CL_1(q)&\leq \min(\CL_{\bz_{i}}(q), \CL_{\bz_{t_{k+1}}}(q)) = \CP_1(q)+\GrO(1), \label{Eq inter Prop 3-sytème partiel 4}\\
    \CL_2(q) &\leq \max(\CL_{\bz_{i}}(q), \CL_{\bz_{t_{k+1}}}(q)) = \CP_2(q)+\GrO(1). \label{Eq inter Prop 3-sytème partiel 1}
\end{align}
On the other hand, we also have $\CL_1^*(q) \leq \CL_{\by_{i}}^*(q)$. Let us prove that $\by_i$ is proportional to a minimal point and that there exists a contant $c>0$ depending only on $\xi$ such that for each $i$ large enough and for each $q\in I_i$, we have:
\begin{equation}
    \label{Eq inter Prop 3-sytème partiel 3}
    |\CL_1^*(q) -\hCL_{i}^*(q)| \leq c,
\end{equation}
\ie $\by_i$ realizes $\CL_1^*$ on $I_i$ within $\GrO(1)$. Indeed, let $q$ belong to $I_i$ and $\bx\neq 0$ be a point (which depends at first sight on $q$) such that $\CL_1^*(q) = \CL_{\bx}^*(q)$, and suppose that $\by_i$ is not proportional to $\bx$. Then, necessarily $\CL_{\bx}^*(q) \leq \CL_{\by_{i}}^*(q)$, and by the definition of the successive minima we obtain $\CL_2^*(q)\leq \CL_{\by_{i}}^*(q)$. Thus, since $\CL_2 = -\CL_2^* + \GrO(1)$ by Proposition \ref{Prop Mahler}, we have the inequality
\[
    -\CL_{\by_{i}}^*(q) + \GrO(1) \leq \CL_2(q).
\]
This inequality combined with inequality \eqref{Eq inter Prop 3-sytème partiel 1} and with assertion \ref{enum 2 Prop Dessin L_j, L_j*} of Proposition \ref{Prop Dessin L_j, L_j*} yields:
\begin{align}
    \label{Eq inter Prop 3-sytème partiel 2}
    \CP_3(q) = -\hCL_i^*(q) &= -\CL_{\by_{i}}^*(q) + \GrO(1) \leq \CP_2(q) + \GrO(1) = \max\big(\hCL_{i}(q), \hCL_{t_{k+1}}(q)\big) + \GrO(1),
\end{align}
for $q\in I_i$. According to assertion \ref{enum 5 Prop Dessin L_j, L_j*} of Proposition \ref{Prop Dessin L_j, L_j*} (cf Figure \ref{figure 3system_partiel1.png}) and by \eqref{Eq Prop Dessin L_j, L_j* croissance en q_i}, there exists $c'>0$ which does not depend on $i$, such that inequality \eqref{Eq inter Prop 3-sytème partiel 2} does not hold if $q\in [a_i+c',b_i-c']$. If $q$ belongs to this interval, we infer that $\by_i$ is proportional to the point $\bx$ (which is a minimal point), and since the content of $\by_i$ is bounded by hypothesis, it follows that $|\hCL_i^*(q) - \CL_{\bx}^*(q)| \leq \GrO(1)$. Then, with the choice of $\bx$, this may be rewritten as
\[
    |\hCL_i^*(q) - \CL_1^*(q)| \leq \GrO(1).
\]
This inequality is valid on $[a_i+c',b_i-c']$ and by $1$-Lipschitz character of $\hCL_i^*$ and $\CL_1^*$ it directly extends to $I_i$, which proves \eqref{Eq inter Prop 3-sytème partiel 3}. By Proposition \ref{Prop Mahler}, this implies
\[
    \CL_3(q)+\GrO(1) \leq -\hCL_i^*(q)=\CP_3(q),
\]
for $q \in I_i$. From \eqref{Eq inter Prop 3-sytème partiel 4}, \eqref{Eq inter Prop 3-sytème partiel 1} and from the previous inequality we deduce
\[
    \CL_j(q) \leq \CP_j(q) + \GrO(1) \quad \textrm{for } j=1,2,3.
\]
Moreover by assertion \ref{enum Prop Lu 2} of Proposition \ref{Prop Lu} and point \ref{Def n-système condition 1} of Definition \ref{Def n-système} of a $3$-system, this implies
\begin{equation}
\label{Eq inter Prop 3-sytème partiel 5}
    \CL_j(q) = \CP_j(q) + \GrO(1) \quad \textrm{for } j=1,2,3,
\end{equation}
for $q\in I_i$. This proves assertion \ref{enum 2 Prop 3-sytème partiel} of Proposition. Assertion \ref{enum 1 Prop 3-sytème partiel} follows from the growth of $\CL_1$ and from \eqref{Eq inter Prop 3-sytème partiel 5}. Indeed, for $i=t_k+l$ with $0< l \leq s_{k+1}$ we have (cf Figures \ref{figure 3system_partiel1.png} and \ref{figure 3system_partiel2.png})
\[
    \CL_1(b_{i-1})+\GrO(1) = \CP_1(b_{i-1}) = \hCL_{t_{k+1}}(b_{i-1}) = \log(\hW_k),
\]
and
\[
    \CL_1(a_{i})+\GrO(1) = \CP_1(a_{i}) = \hCL_{t_{k+1}}(a_{i}) = \log(\hW_k) = \CL_1(b_{i-1})+\GrO(1),
\]
thus
\[
    \CL_1(q)+\GrO(1) = \CP_1(q)
\]
for $q\in [b_i,a_{i+1}] = I_i'$, and since families $(I_j)_j$ and $(I_j')_j$ cover $[M,+\infty)$ for $M$ large enough, this implies \ref{enum 1 Prop 3-sytème partiel}.\\
The assertion \ref{enum 3 Prop 3-sytème partiel} may also be deduced from \eqref{Eq inter Prop 3-sytème partiel 5}. Indeed these inequalities for $q = b_j, a_{j+1}$ imply $\CL_2(q)+\GrO(1) = \CL_3(q)+\GrO(1) = \CP_2(q) = \CP_3(q)$,
and we conclude noticing that $\CL_2,\CL_3$ are continuous piecewise linear with slopes $0$ or $1$.

\end{Dem}

\section{Examples of proper $\psi$-numbers}
\label{Section exemples (deux)}

In Section \ref{Section Exemple de Roy} we show that Roy's matrices \cite[Examples $3.3$]{roy2007two} give $\psi$-Sturmian sequences and proper $\psi$-Sturmian numbers. Thanks to these examples we may prove Corollary \ref{Cor 2 spectre}. In Section \ref{Section Exemple de Bugeaud et Laurent} we present Bugeaud-Laurent's examples. They give proper $\psi$-Sturmian numbers for which the associated $\delta$ given by Proposition \ref{Prop existence delta sur le det} satisfies $\delta=0$.

\subsection{Roy's matrices}
\label{Section Exemple de Roy}

Let $\psi\in\FF$ be a Sturmian function and let $\sigma$ denote the real number defined by \eqref{Eq Def sigma}. In this section we show that Roy's examples (see Examples $3.3$, $4.3$ and $5.4$ of \cite{roy2007two}) provide examples of proper $\psi$-Sturmian numbers. The main result due to Roy's work \cite{roy2007two} is Proposition \ref{Prop densite delta}.

\begin{Exe}[Roy, 2007] \
\label{Exemple 3.3}
Example $3.3$ of \cite{roy2007two} remains valid in our case: if we fix an integer triple $\underline{A} = (a, b, c)$ such that $a\geq2$ and $c\geq b \geq 1$, then the $\psi$-Sturmian sequence $(\bw_i)_{i\geq 0}$ defined by
\[
    \bw_0 = \left(
    \begin{array}{cc}
        1 & b \\
        a & a(b+1)
    \end{array}
    \right),
    \quad
    \bw_1 = \left(
    \begin{array}{cc}
        1 & c \\
        a & a(c+1)
    \end{array}
    \right)
\]
is admissible with the matrix
\[
    \transpose{\,N} = \left(
    \begin{array}{cc}
        -1+a(b+1)(c+1) & -a(b+1) \\
        -a(c+1) & a
    \end{array}
    \right).
\]

\end{Exe}

With these matrices we have the following result.

\begin{Prop}
\label{Prop exemple psi-sturmien from Roy}
The $\psi$-Sturmian sequence $(\bw_i)_i$ has a multiplicative growth (we may choose the constants $c_1=1$ and $c_2=2$ in \eqref{Eq Def croissance multiplicative}) and $(\norm{\bw_i})_i$ tends to infinity. Moreover for each $i\geq 0$ the matrix $\bw_i$ is primitive, and $\Tr(\bw_i)$ and $\det(\bw_i)$ are relatively prime. In particular, $(\bw_i)_i$ gives a $\psi$-Sturmian number $\xi_{\underline{A}}$.
\end{Prop}

\begin{Rem2}
The hypotheses of Corollary \ref{Cor arithm Tr, det, y_i, z_j} are fulfilled. As Roy points out \cite[Example $5.4$]{roy2007two}, since $\det(\bw_0)=\det(\bw_1)=a$, the condition \eqref{Eq encadrement det alpha beta} is true for $i=0,1$ with
\[
    \alpha = \frac{\log(a)}{\log\big(2a(c+1)\big)}\quad\textrm{and}\quad \beta = \frac{\log(a)}{\log\big(2a(b+1)\big)}.
\]
\end{Rem2}

\begin{Dem}
The multiplicative growth with the announced implicit constants is implied by Lemma \ref{Lemme croissance multiplicative}. The fact that $(\norm{\bw_i})_i$ tends to infinity follows from $c_1=1$ and from $\norm{\bw_0},\norm{\bw_1}>1$ (since $a>1$). Considering the matrices $\bw_i$ modulo $a$, Roy shows \cite[Example $4.3$]{roy2007two} that $\Tr(\bw_i)$ and $\det(\bw_i)$ are relatively prime, which implies that $\bw_i$ is primitive.

\end{Dem}

We denote by $\delta_{\underline{A}}\geq 0$ the real number provided by Proposition \ref{Prop existence delta sur le det}. We have the inequalities
\begin{equation}
\label{Eq inter encadrement delta_A}
    \alpha \leq \delta_{\underline{A}} \leq \beta.
\end{equation}

\begin{Prop}
\label{Prop densite delta}
The set $\{\delta_{\underline{A}}\;|\; \underline{A} = (a,b,c) \textrm{ with $a\geq2$ and $c\geq b \geq 1$} \}$ is dense in $[0,\frac{\sigma}{1+\sigma}]$.
\end{Prop}

\begin{Dem}
Cf \cite[Proof of Corollary $7.2$]{roy2007two}. Fix $\delta\in \big(0,\frac{\sigma}{1+\sigma}\big)$ and $\ee > 0$. Then, there are two integers $k$ and $l$ with $0< l < k$ and
\[
    \delta-\ee \leq \frac{l}{k+2} \leq \frac{l}{k} \leq \delta.
\]
Consider the integer triple $\underline{A} = (a,b,c)$ with $a=2^l, b=2^{k-l}-1$ and $c=2^{k-l}$. With this choice we have $\alpha = l/(k+2)$ and $\beta =  l/k$. Thus by \eqref{Eq inter encadrement delta_A} $\delta_{\underline{A}}$ satisfies $|\delta-\delta_{\underline{A}}| \leq \ee$.

\end{Dem}

\subsection{Bugeaud and Laurent's examples}
\label{Section Exemple de Bugeaud et Laurent}

In this section we briefly present Bugeaud and Laurent's examples and show that their numbers $\xi_{\phi'}$ \cite{bugeaud2005exponentsSturmian} are (proper) $\psi$-Sturmian numbers for which the associated $\delta$ satisfies $\delta = 0$.\\

Let $\underline{s}'=(s_k')_{k\geq 1}$ be a sequence of positive integers. Fix $a$ and $b$ two distinct positive integers. We define inductively a sequence of words $(m_k)_{k\geq 0}$ by the formulas
\[
    m_0 = b,\quad m_1=b^{s_1'-1}a\quad\textrm{and}\quad m_{k+1}=m_k^{s_{k+1}'}m_{k-1}\quad(k\geq 1).
\]
This sequence converges to the infinite word
\[
    m_{\phi} = \lim_{k\rightarrow+\infty}m_k = b^{s_1'-1}a\dots,
\]
which is usually called the \textsl{Sturmian characteristic word of angle} (or of \textsl{slope}) $\phi := [0; s_1', s_2', s_3', \dots]$ constructed on the alphabet $\{a, b\}$. Denote by $\xi_{\phi'} = [0; m_{\phi'}]$ the real number whose partial quotients are successively $0$, and the letters of the infinite word $m_{\phi'}$.\\

Now, set $\underline{s}=(s_k)_{k\geq 0}$ be the sequence defined by $s_0=-1$, $s_1=1$ and $s_k=s_{k}'$ ($k\geq 2$). We denote by $\psi$ the Sturmian function associated with $\underline{s}$ and by $\sigma$ the real number defined by \eqref{Eq Def sigma}. We define a $\psi$-Sturmian sequence $(\bw_i)_{i\geq 0}$ whose first terms are
\[
\bw_0 =
\left(
\begin{array}{cc}
b & 1 \\
1 & 0
\end{array}
\right)
\quad\textrm{and}\quad
\bw_1 =
\left(
\begin{array}{cc}
b & 1 \\
1 & 0
\end{array}
\right)^{s_1'-1}
\left(
\begin{array}{cc}
a & 1 \\
1 & 0
\end{array}
\right).
\]
Then, thanks to Lemma $5.3$ of \cite{bugeaud2005exponentsSturmian}, the sequence $(\bw_i)_i$ is admissible with the matrix
\[
N =
\left[
\left(
\begin{array}{cc}
a & 1 \\
1 & 0
\end{array}
\right)
\left(
\begin{array}{cc}
b & 1 \\
1 & 0
\end{array}
\right)
\right]^{-1}.
\]

\begin{Prop}
\label{Prop exemple psi-sturmien from Bugeaud et Laurent}
The sequence $(\bw_i)_i$ satisfies the hypotheses of Section \ref{section xi construction de Roy} and the associated (proper) $\psi$-Sturmian number is $\xi_{\phi'}$. Moreover, the associated $\delta$ (given by Proposition \ref{Prop existence delta sur le det}) is equal to $0$.
\end{Prop}

\begin{Dem} First note that $(\bw_i)_i$ has multplicative growth (each $\bw_i$ belongs to $\MM$), the sequence $(\norm{\bw_i})_i$ tends to infinity and $|\det(\bw_i)|=1$ for each $i$ (so that the associated $\delta$ is equal to $0$).
We denote by $(\by_i)_i$ the sequence of symmetric matrices associated with $(\bw_i)_i$ and by $l_1,l_2,\dots$ the letters of the word $m_{\phi'}$. Let $p_n/q_n$ denote the $n$-th convergent of $\xi_{\phi'}$. The classical properties of continued fraction expansions give (see for instance \cite[Chapter I]{schmidt1996diophantine})
\[
\left(
\begin{array}{cc}
q_n & q_{n-1} \\
p_n & p_{n-1}
\end{array}
\right)
=
\left(
\begin{array}{cc}
l_1 & 1 \\
1 & 0
\end{array}
\right)\dots
\left(
\begin{array}{cc}
l_n & 1 \\
1 & 0
\end{array}
\right).
\]
Since each matrix $\by_i$ is of the last form for a particular $n$ (by Lemma $5.3$ of \cite{bugeaud2005exponentsSturmian}), this implies directly that the $\psi$-number given by Proposition \ref{Prop existence y = (1,xi,xi^2)} is equal to $\xi_{\phi'}$.

\end{Dem}

\section{Proofs of Theorem \ref{Thm propre exposants intro} and Corollary \ref{Thm densité omega_2 sturm dans [c,+infty]}}
\label{section preuve du corollaire omega_2 densité}

Let $\psi\in\Sturm$ be a Sturmian function. We denote by $\sigma$ the real number defined by \eqref{Eq Def sigma}. In this section we prove Theorem \ref{Thm propre exposants intro} and Corollary \ref{Thm densité omega_2 sturm dans [c,+infty]} stated in the introduction.\\

\begin{Dem}[of Theorem \ref{Thm propre exposants intro}]
$\newline$
Let us set
\[
\Delta_{\psi} = \{\delta_{\underline{A}}\;|\;  \underline{A} = (a,b,c) \textrm{ with $a \geq 2$, $c \geq b \geq 1$}\textrm{ and }\delta_{\underline{A}} < \frac{\sigma}{1+\sigma}\}\cup\{0\}.
\]
Proposition \ref{Prop densite delta} implies that $\Delta_{\psi}$ is dense in $[0,\sigma/(1+\sigma)]$. By virtue of Propositions \ref{Prop exemple psi-sturmien from Roy} and \ref{Prop exemple psi-sturmien from Bugeaud et Laurent}, for each $\delta\in\Delta_{\psi}$, there exists a proper $\psi$-Sturmian number $\xi$ such that $\delta$ is the number associated to $\xi$ by Proposition \ref{Prop existence delta sur le det}. We conclude by Theorem \ref{Thm exposants geom param nb psi-sturmien} together with Proposition \ref{Prop dico exposants}.

\end{Dem}

\begin{Dem}[of Corollary \ref{Thm densité omega_2 sturm dans [c,+infty]}]
$\newline$
In order to establish his results, Cassaigne \cite{cassaigne1999limit} shows that
\begin{equation}
\label{Thm intro spec omega2 eq inter -1}
    \Sturm = \{1/[\bb]\;|\;\bb\in(\ZZ_{\geq 1})^{\NN}\textrm{ such that for each $k\geq 0$ we have } [\bb]\geq[\mathrm{T^k\bb}]\},
\end{equation}
where $[\bb]$ is the continued fraction $[b_0;b_1,\dots]$ and $\mathrm{T}$ is the right shift operator ($[\mathrm{T^k\bb}]$ is the real number $[b_k;b_{k+1},\dots]$).\\
According to Theorem \ref{Thm propre exposants intro}, for each $\sigma\in\Sturm$ the set of $\omega_2(\xi)$ with $\xi$ a proper $\psi$-Sturmian number and $\sigma(\psi)=\sigma$ is dense in the interval $\Big[2/\sigma,1+2/\sigma\Big]$. It is therefore sufficient to show that
\[
    \bigcup_{\sigma\in\Sturm}\Big[\frac{2}{\sigma},\frac{2}{\sigma}+1\Big] = [1+\sqrt{5},2+\sqrt 5]\cup[2+2\sqrt 2,3+2\sqrt 3]\cup[3+\sqrt{13}+\infty).
\]
For $n\geq 1$ and $1\leq a \leq n$ we define the sequence $\bu_{a,n} = (s_k)_{k\geq0 }$ by
\[
\left\{ \begin{array}{l}
s_{2k} = n,\\
s_{2k+1} = a,
\end{array}
\right.
\]
for each $k\geq 0$. It is clear that for each  $k\geq 0$ we have $[\bu_{a,n}]\geq[\mathrm{T^k\bu_{a,n}}]$, thus $1/[\bu_{a,n}]\in\Sturm$, and
\[
    [\bu_{a,n}] = [n;\overline{a,n}] =\frac{an+\sqrt{(an)^2+4an}}{2a}.
\]
Then we define
\[
    \delta_{a,n} = 2[\bu_{a,n}] = n+n\sqrt{1+\frac{4}{an}}.
\]
Let us study the union of intervals $[\delta_{a,n},\delta_{a,n}+1]$ for $n\geq 1$ and $1\leq a\leq n$.\\
We remark that
\begin{equation}
\label{Thm intro spec omega2 eq inter 1}
0\leq \delta_{a,n}-\delta_{a+1,n} \leq 1 \quad (n\geq 2 \textrm{ and } 1\leq a < n),
\end{equation}
and
\begin{equation}
\label{Thm intro spec omega2 eq inter 2}
|\delta_{n+1,n+1}-\delta_{1,n}| < 1 \quad (n\geq 3).
\end{equation}
We may deduce from \eqref{Thm intro spec omega2 eq inter 1} and from \eqref{Thm intro spec omega2 eq inter 2} that
\begin{align*}
\label{Thm intro spec omega2 eq inter 3}
\bigcup_{\substack{1\leq n \\ 1\leq a \leq n }} [\delta_{a,n},\delta_{a,n}+1] &= [\delta_{1,1},\delta_{1,1}+1]\cup[\delta_{2,2},\delta_{1,2}+1]\cup[\delta_{3,3},+\infty) \\
&= [1+\sqrt{5},2+\sqrt 5]\cup[2+2\sqrt 2,3+2\sqrt 3]\cup[3+\sqrt{13},+\infty).
\end{align*}
As pointed out by Bugeaud and Laurent in \cite{bugeaud2005exponentsSturmian}, the two greatest values of $\sigma\in\Sturm$ are $(-1+\sqrt5)/2$ and $-1+\sqrt 2$. This implies that the two smallest values of $2/\sigma$ for $\sigma\in\Sturm$ are $\delta_{1,1}=1+\sqrt 5 = 3.23\dots$ and $\delta_{2,2}=2+2\sqrt 2 = 4.82\dots$. To conclude, it is enough to show that there does not exist any $\sigma\in\Sturm$ such that $\delta_{1,2}<2/\sigma < \delta_{3,3}$. By \eqref{Thm intro spec omega2 eq inter -1}, it suffices to show that if $\bb\in(\ZZ_{\geq 1})^{\NN}$ is such that $[\bb]\geq[\mathrm{T}^k\bb]$ for each $k\geq 0$ and if
\[
    [2;\overline{1,2}] \leq [b_0;b_1,\dots] \leq [3;\overline{3}],
\]
then, either $[b_0;b_1,\dots] = [2;\overline{1,2}]$ or $[b_0;b_1,\dots]  = [3;\overline{3}]$.\\
Note that if $\ba=(a_k)_{k\geq 0}$ and $\ba'=(a_k')_{k\geq 0}$ are two distinct sequences of positive integers and if $k$ is the smallest index for which $a_k\neq a_k'$, then $[\ba] > [\ba']$ if and only if $a_k > a_k'$ if $k$ is even or $a_k<a_k'$ if $k$ is odd.\\
Set $x=[b_0;b_1,\dots]$. The integer part of $x$ is $b_0$ and is equal to $2$ or $3$. Since $x\geq [\mathrm{T}^k\bb]$ for each $k\geq 0$, we have $1\leq b_k\leq b_0$ for all $k$. Suppose that $b_0 = 2$. Then, with the previous remark we may show by induction that $x = [2;\overline{1,2}]$. Similarly, if $b_0=3$ we show that $x=[3;\overline{3}]$.\\

\end{Dem}

\bibliographystyle{abbrv}

\Addresses

\end{document}